\documentclass[twocolumn]{autart}

\usepackage{amssymb}
\usepackage{amsmath}
\usepackage{amsfonts}
\usepackage{graphicx}
\usepackage{epsfig}
\usepackage{cite}
\usepackage{url}
\usepackage{picins}
\usepackage{xcolor}
\usepackage{tikz}
\allowdisplaybreaks[4]

\newtheorem{lemma*}{Lemma}

\DeclareMathOperator{\esssup}{ess\,sup}

\begin{document}

\begin{frontmatter}

\title{Backstepping for Partial Differential Equations: A Survey\thanksref{t1}\thanksref{t2}} % Title, preferably not more 
                                                % than 10 words.

\thanks[t1]{This paper was not presented at any IFAC 
meeting. Corresponding author R~Vazquez. Tel. +34954488148. }
\thanks[t2]{R. Vazquez acknowledges support of grant PID2023-147623OB-I00 funded by MICIU/AEI/10.13039/501100011033 and by ``ERDF A way of making Europe''. The work of M. Krstic was funded by AFOSR grant FA9550-23-1-0535  and NSF grant ECCS-2151525. The work of J. Auriol was funded by Agence Nationale de la Recherche (ANR) via the grant PANOPLY ANR-23-CE48-0001-01}

%\thanks[footnoteinfo]{Funding. This paper was not presented at any IFAC meeting. Submitted to \emph{Automatica} 2024. Corresponding author R.~Vazquez. Tel. +34-954488148.}

\author[vazquez]{Rafael Vazquez}\ead{rvazquez1@us.es},
\author[auriol]{Jean Auriol}\ead{jean.auriol@centralesupelec.fr},
\author[bribiesca]{Federico Bribiesca-Argomedo}\ead{federico.bribiesca-argomedo@insa-lyon.fr},
\author[krstic]{Miroslav Krstic}\ead{krstic@ucsd.edu}%
\address[vazquez]
{Departamento de Ingenier\'ia Aeroespacial, Universidad de Sevilla, 41092 Seville, Spain}
\address[auriol]{Universit\'e Paris-Saclay, CNRS, CentraleSup\'elec, Laboratoire des Signaux et Syst\`emes, 91190, Gif-sur-Yvette, France}
\address[bribiesca]{INSA Lyon, Universite Claude Bernard Lyon 1, Ecole Centrale de Lyon, CNRS, Amp\`ere, UMR5005, 69621
Villeurbanne, France}
\address[krstic]
{Department of Mechanical and Aerospace Engineering, University of California San Diego, La Jolla, CA 92093-0411}%

% \begin{keyword}
% Distributed parameter systems;
% Stabilization;
% Feedback Linearization;
% Partial differential equations;
% Lyapunov Function;
% Boundary Conditions;
% Others;
% \end{keyword}
\begin{abstract}
Systems modeled by partial differential equations (PDEs) are at least
as ubiquitous as those by nature finite-dimensional and modeled by ordinary differential equations (ODEs). And
yet, systematic and readily usable methodologies,
for such a significant portion of real systems, have been historically
scarce. Around the year 2000, the backstepping approach
to PDE control began to offer not only a less abstract alternative to
PDE control techniques replicating optimal and spectrum
assignment techniques of the 1960s, but also enabled the methodologies
of adaptive and nonlinear control, matured in the
1980s and 1990s, to be extended from ODEs to PDEs, allowing feedback
synthesis for systems that
are uncertain, nonlinear, and infinite-dimensional. The PDE
backstepping literature has since grown to hundreds of papers and
nearly a dozen books.
This survey aims to facilitate the entry into this thriving area of
overwhelming size and topical
diversity. Designs of controllers and observers, for parabolic,
hyperbolic, and other classes of PDEs, in one or more dimensions, with
nonlinear, adaptive,
sampled-data, and event-triggered extensions, are covered in the
survey. The lifeblood of control are technology and physics.
The survey places a particular emphasis on applications that have
motivated the development of the theory and which have
benefited from the theory and designs: flows, flexible structures,
materials, thermal and chemically reacting dynamics, energy (from oil
drilling to batteries and magnetic confinement fusion), and vehicles.
 \end{abstract}
\end{frontmatter}

%\tableofcontents 

\section{Introduction}
Partial differential equations (PDEs) are at the heart of numerous scientific and engineering disciplines, providing a powerful framework for modeling and understanding complex systems that exhibit variations across space and time. They are ubiquitous in modeling spatiotemporal phenomena, such as heat transfer, fluid dynamics, structural vibrations, and chemical reactions. Whether it is modeling traffic to alleviate congestion on highways, simulating heat transfer and phase change processes in 3D printing, or predicting the behavior of fluids in various systems, PDEs provide a powerful framework for understanding and controlling complex processes. They also find applications in estimating the state of charge in batteries, addressing challenges in drilling operations, and mitigating thermoacoustic instabilities. This versatility underscores the significance of PDEs in diverse fields, ranging from transportation and manufacturing to energy storage and aerospace engineering. Thus, the study of PDEs remains a vibrant and essential area of research, bridging the gap between abstract mathematical theory and its application in the natural and engineered world. 

More specifically, the control of systems governed by PDEs has applications in many branches of science and engineering.  The ability to effectively control these systems allows for optimizing performance, ensuring safety, and achieving desired behaviors in a wide range of applications.

In this context, backstepping for PDEs has emerged as a powerful and versatile methodology for crafting boundary controllers and observers in the infinite-dimensional setting. Over the past two decades, this constructive approach has been successfully applied to a wide spectrum of PDE systems, as evidenced by the burgeoning number of publications in this domain, comprising not only several hundreds of papers but also almost a dozen of books, see e.g. \cite{vaz2008,kog2020_e,BekiarisLiberis2013,krstic5,meurer2,adaptive,wang2022pde,Yu2023,anfinsen2019adaptive}, some with over a thousand citations. The only \mbox{(then-)general} textbook was published in 2008 \cite{krs2008}, and in the last 15 years, the field has at least quintupled in size and breadth. The proliferation of research has reached a point where it is increasingly challenging for novices to navigate and digest the extensive literature. 

This survey aims to provide a unified and accessible exposition of the backstepping methodology, elucidating its theoretical underpinnings, key extensions, practical applications, and implementation insights. We endeavor to highlight the broad contexts in which backstepping has proven effective, spanning parabolic, hyperbolic, and other classes of PDEs, while also delving into the nuances that have propelled its success. Furthermore, we offer an in-depth analysis of the open challenges and potential research avenues in this continuously evolving field.

\subsection{Origins of PDE backstepping}\label{sect-origins}
The roots of the method can be traced back to the 1970s and 1980s when Volterra-type integral transforms were employed to investigate solutions and controllability properties of PDEs in a handful of lesser-known works \cite{colton,seid1984}. However, it was not until the early 2000s that these transformations were independently rediscovered and leveraged for feedback design and stabilization.

Backstepping, as a control design methodology, first gained prominence in the context of finite-dimensional systems. The origins of the backstepping idea are difficult to trace, because of its simultaneous appearance, often implicit, in the works of Tsinias~\cite{tsinias1989sufficient}, Koditschek~\cite{koditschek1987adaptive}, Byrnes and Isidori~\cite{byrnes1989new}, and Sontag and Sussmann~\cite{sontag1989further}. Backstepping was initially conceived for nonlinear systems with an integrator at the input, such as in robotic systems where the kinematics are separated from the input torque by angular velocity, which integrates in time the torque and other sources of angular acceleration. The recursive application of backstepping led to its generalization to {\em strict-feedback} nonlinear systems, which are of the form $\dot x_i = x_{i+1} + \varphi_i(x_1,\ldots,x_i), i=1,\ldots,n, x_{n+1} =u$. The term (strict)-feedback refers to the ``recirculation" of the states of the integrator chain $\dot x_i = x_{i+1}$ into the $i$-th equation through the nonlinearity $\varphi_i(x_1,\ldots,x_i)$, which is free of the dependencies on the `feedforward’ state variables $x_{i+1}$ through $x_n$. 
%Backstepping, as a control design methodology, first gained prominence in the context of finite-dimensional systems, particularly for a class of nonlinear systems known as strict-feedback systems. 
The book~\cite{KKK}, often regarded as the definitive backstepping reference, has garnered thousands of citations. The popularity of backstepping in the finite-dimensional setting can be attributed to several factors: the ubiquity of strict feedback systems in applications, the constructive nature of the method, and its ability to yield explicit control laws---features seldom found in nonlinear control. Moreover, the approach can be extended to design adaptive controllers \cite{KKK}.

While finite-dimensional backstepping took root in the 1990s, its extension to the infinite-dimensional realm of PDEs emerged a decade later. Single-step ODE backstepping designs were reported in the hyperbolic case in~\cite{coron1998stabilization} and in the parabolic case in~\cite{liu2000backstepping}. The early 2000s then witnessed the efforts of a team at the University of California, San Diego in applying backstepping to PDE systems. These initial works, such as \cite{Boskovic2001}, exhibited some of the hallmarks of PDE backstepping: state transformations, explicit feedback laws, and Lyapunov-based stability analysis. However, they were limited in scope, focusing on specific classes of convection-reaction equations with restricted instability margins. Concurrently, attempts were made to extend backstepping to PDEs via spatial discretization followed by control design \cite{Boskovic2003,Balogh2002}. While the resulting control laws were convergent, these ``discretize-then-design'' approaches suffered from convergence issues as the discretization step tended to zero.

A key moment for PDE backstepping arrived with the work of Smyshlyaev and Krstic \cite{Smyshlyaev2004} after Liu~\cite{liu2003boundary} applied Colton’s method~\cite{colton} to the approach introduced in~\cite{Boskovic2003,Balogh2002}, for the study of the well posedness of the kernel PDE. Focusing on a general 1-D linear reaction-diffusion-advection PDE, they introduced a ``design-then-discretize'' strategy, often called late lumping.  The term \emph{late-lumping} is used in contrast with \emph{early-lumping} approaches, where the controller design is based on a finite-dimensional approximation of the PDE~\cite{morris2010control}.  
  Smyshlyaev and Krstic further extended their methodology to the design of boundary observers in \cite{Smyshlyaev2005}, by exploiting duality principles. These initial contributions by Smyshlyaev and Krstic laid the bedrock for the field of PDE backstepping, and are explored in detail in Section~\ref{Sec_primer}.

\subsection{Alternative methods for PDE Control}

The infinite-dimensional nature of PDE systems poses unique challenges for control design. Unlike ordinary differential equations (ODEs), PDEs exhibit spatially distributed dynamics, which often lead to complex boundary conditions, spatially varying coefficients, and non-local interactions. Classical control techniques developed for finite-dimensional systems are not always directly applicable or may lead to suboptimal performance when applied to PDEs \cite{Lasiecka2000,Curtain2012}.

Over the years, various approaches have been proposed to tackle the control of PDE systems. One prominent framework is the linear quadratic regulator (LQR) theory, which has been extended to infinite-dimensional systems \cite{krener2021optimal,Bensoussan2007}. The LQR approach seeks to minimize a quadratic cost functional subject to the PDE dynamics, leading to optimal feedback controllers. However, the resulting controllers require a solution to an operator (infinite-dimensional) Riccati equation and may be challenging to implement in practice. 
A related classical approach is the Kalman filter and its extensions for state estimation. While backstepping observers provide guaranteed convergence properties through constructive design, Kalman filtering approaches optimize the statistical properties of the estimation error. For instance, a recent work~\cite{afshar2025} extends the Extended Kalman Filter to semilinear infinite-dimensional systems, providing state estimates for nonlinear PDEs under detectability assumptions.

Another traditional approach is the spectral decomposition method, which exploits the eigenfunction expansion of the PDE operator to transform the system into an infinite series of ODEs~\cite{lagnese2012domain,baker2000finite}. This allows for the application of finite-dimensional control techniques to each mode separately. For instance, pole-placement theorems with spectral decomposition of Sturm-Liouville operators can be used to design finite-dimensional output feedback controllers~\cite{lhachemi2022finite}.  However, this approach is limited to systems with simple geometries and boundary conditions, and it may suffer from slow convergence and numerical instabilities. In addition, it is not clear how such methods could be made adaptive when the PDE's parameters are unknown since the dimension of the reduced model would vary in real-time, as the parameter estimate varies.

Lyapunov-based methods have also been widely employed for PDE control \cite{mazenc2011strict,Coron2007}. These methods involve constructing a Lyapunov functional that captures the energy or stability properties of the system. The control law is then designed to ensure the decay of the Lyapunov functional along the system trajectories. While Lyapunov-based designs offer a systematic framework, finding suitable Lyapunov functionals and handling boundary conditions can be challenging, especially for complex PDE systems.

In recent years, flatness-based control has emerged as a powerful tool for PDEs \cite{meurer2,Rudolph2003}. Flatness refers to the property of a system where all states and inputs can be parametrized in terms of a finite set of independent variables, called flat outputs. For flat systems, trajectory planning and tracking can be performed in a straightforward manner. However, the applicability of flatness-based control is limited to systems that possess the flatness property, and the identification of flat outputs can be nontrivial.

Model reference adaptive control (MRAC) approaches have also been extended to distributed parameter systems. For example, Böhm et al. \cite{bohm1998model} developed an infinite-dimensional MRAC framework for a broad class of nonlinear distributed parameter systems. While the approach provides theoretical results on well-posedness, stability, and parameter convergence along with finite-dimensional approximation schemes, it requires full-state measurements and fully distributed (in-domain) control inputs, which are seldom found in practical applications.

Similarly, Port-Hamiltonian systems (PHS) formulations have emerged as a powerful framework for the modeling and control of distributed parameter systems~\cite{Jacob-2012-book,vanderSchaft_2002_DPS,duindam2009modeling}. This formalism is well-suited for capturing the dynamics of large-scale (multi)physics systems and leveraging the physical properties of the systems to design boundary controllers~\cite{wu2017port}. The energy shaping method~\cite{Macchelli-2017-TAC_synthPHS}, for instance, enables the alteration of the closed-loop energy function, providing a design approach with a clear physical interpretation. However, PHS methods necessitate a specific configuration and may encounter challenges when employing boundary controllers to directly modify in-domain coupling terms.

%OTHERS? WHO IS MISSING? LOOK AT THE DPS TCs.\color{black}

In summary, the control of PDE systems is a rich and active field, with various approaches tailored to different classes of systems and control objectives. While backstepping has emerged as a prominent methodology, offering a systematic and constructive framework for boundary control design, it is important to recognize the complementary nature of all the different approaches, their tradeoffs, and the potential for hybrid strategies that combine the strengths of multiple methods.

\subsection{Structure of the survey}
%TO BE COMPLETED AT THE END

In this survey, we embark on a comprehensive exploration of the backstepping methodology and its myriad applications. Section \ref{Sec_primer} provides a comprehensive primer on PDE backstepping fundamentals, introducing key concepts and methodologies. Section \ref{Sec_hyp_systems} considers interconnected hyperbolic and parabolic systems, exploring various configurations and their control challenges. Section \ref{Sec_advanced_topics} covers advanced topics such as adaptive control, higher-dimensional domains, nonlinear PDEs, and moving boundary problems. The practical applications of backstepping are showcased in Section \ref{sec-apps}, ranging from traffic control to multi-agent systems or flow control. Next, Section \ref{sec-implementation} addresses implementation considerations, including numerical solutions of kernel equations and machine learning approaches. Finally, Section \ref{sec-conclusions} concludes the survey by discussing open problems and future research directions in PDE backstepping. As we navigate through this extensive body of work, we strive to strike a balance between accessibility for those new to the field and depth for advanced practitioners. By synthesizing the diverse contributions and providing a unified perspective, we aim to facilitate a deeper understanding of the backstepping methodology and its potential for tackling complex control problems in the infinite-dimensional realm.

\section{A primer in backstepping for PDEs} \label{Sec_primer}
Backstepping for PDEs came together as a general methodology in the highly-cited 2004 paper~\cite{Smyshlyaev2004}. 
Dealing with a 1-D linear reaction-diffusion-advection PDE, the method follows a design-first-then-discretize approach, sometimes called late-lumping (as the last step in the implementation is to approximate the controller by a \emph{lumped parameter} finite-dimensional system). It encompasses three fundamental components:
\begin{enumerate}
    \item Selection of a \emph{target system} that exhibits desired properties (like stability, typically validated by a Lyapunov function) while retaining resemblance to the original system;
    \item Utilization of an integral transformation, known as the \emph{backstepping transformation}, to map the original system to the target system within appropriate functional spaces. This transformation's boundedness and invertibility are crucial to guarantee that the original and target systems share equivalent stability properties;
    \item Derivation of \emph{kernel equations} from both the original and target systems PDEs, together with the transformation. The well-posedness of these equations, typically defined on a triangular domain and of Goursat type (hyperbolic boundary problems), is found by converting them into integral equations and applying successive approximations.
\end{enumerate}
These elements are integral to the backstepping methodology and interlinked; a well-chosen target system leads to solvable kernel equations and an invertible transformation. Conversely, a poorly selected target system may result in unsolvable kernel equations or a non-invertible transformation.

Expanding on these concepts, Smyshlyaev and Krstic introduced a dual methodology for designing both collocated and anti-collated observers using boundary measurements in 2005~\cite{Smyshlyaev2005}. Depending on whether the measurements are collocated or anti-collocated, different backstepping transformations are applied.  By integrating the observer with the controller, output-feedback controllers were derived, whose stabilization efficacy was straightforward to demonstrate via the properties of the backstepping transformations.

This section explores the backstepping methodology for PDEs in more detail using the aforementioned reaction-diffusion designs as a basic paradigmatic example. After introducing some notations and basic PDE mathematics, we explore the selection of target systems, transformations, and kernel equations, culminating in the practical application of these concepts to ensure system stability. The discussion also covers the design of observers and output feedback control laws. We finish the section with other basic early backstepping designs. This gradual exposition lays the foundation for understanding how backstepping principles apply to more
complex systems, explored in later sections.

\subsection{Mathematical Preliminaries: Functional Spaces, Norms, Inequalities, and Stability}\label{sect-math}
Unlike finite-dimensional systems governed by ODEs, PDEs involve states that are functions of both time and space, requiring tools from functional analysis for their study and control. This section briefly introduces some foundational mathematical concepts and notations central to PDE control, focusing on 1D spatial domains\footnote{Higher-dimensional domains are discussed in Section~\ref{sec-higher}}, due to their prevalent use in the backstepping literature, and specifically in the interval $(0,1)$ for simplicity. 

Considering functions $f:(0,1) \rightarrow \mathbb{R}$, the most basic functional spaces in PDE theory are the $L^p(0,1)$ spaces, characterized by integrability properties; for simplicity, we just write $L^p$ obviating the domain unless relevant for the discussion. For $1 \leq p < \infty$, the $L^p$ norm of $f$ is defined as $\Vert f \Vert_{L^p} = \left(\int_0^1 |f(x)|^p  dx\right)^{1/p}$, where the integral might diverge for functions not belonging to $L^p$. Accordingly, the space $L^p$ includes all functions for which this norm is finite, namely $L^p = \{ f : \Vert f \Vert_{L^p} < \infty \}$. On the other hand, the $L^\infty(0,1)$ space is defined by the essential supremum norm $L^\infty = \{ f : \Vert f \Vert_{L^\infty} = \esssup_{x \in (0,1)} |f(x)| < \infty \}$, capturing essentially bounded functions. All $L^p$ spaces are Banach spaces (complete normed vector spaces).  For $p\neq q$, the spaces $L^p$ and $L^q$ are distinct since their norms are not equivalent\footnote{Norm equivalence between two norms $\|\cdot\|_a$ and $\|\cdot\|_b$ means there exist constants $c_1, c_2 > 0$ such that $c_1\|f\|_a \leq \|f\|_b \leq c_2\|f\|_a$ for all $f$ in the space. This concept is crucial in the closed-loop stability proofs of backstepping.}; however, within a finite domain, there exists a nesting relation: $L^p \subset L^q$ for $p > q$, which does not apply to unbounded domains.

 Of the $L^p$ spaces, $L^2$ is the most used in backstepping, as it is distinguished by being a Hilbert space (a complete inner product vector space), equipped with the inner product $(f,g)_{L^2} = \int_0^1 f(x)g(x)  dx$, facilitating the use of geometric and projection methods. However, properties of functions, such as continuity, boundedness, or differentiability, are not guaranteed within $L^2$ spaces alone. For this purpose, functional spaces that account for spatial derivatives are necessary; these are known as Sobolev spaces. We define two primary Sobolev spaces, $H^1$ and $H^2$,  with norms defined as 
$\Vert f \Vert_{H^1}^2= \Vert f \Vert_{L^2}^2  + \Vert f_x \Vert_{L^2}^2$ and 
$\Vert f \Vert_{H^2}^2 = \Vert f \Vert_{L^2}^2  + \Vert f_x \Vert_{L^2}^2 +\Vert f_{xx} \Vert_{L^2}^2$.
The spaces $H^1$ and $H^2$ are composed of all functions, with the corresponding norm being finite. They are both Hilbert spaces equipped with inner products. In this survey, we mostly use the $L^2$, $L^\infty$, $H^1$, and $H^2$ spaces and norms.

Norm inequalities such as the Cauchy-Schwarz, Poincaré, and Agmon's play pivotal roles in PDE analysis, facilitating the estimation and bounding of norms. The Cauchy-Schwarz's inequality asserts that $(f,g)_{L^2} \leq \Vert f \Vert_{L^2} \Vert g \Vert_{L^2}$. Poincaré's inequalities relate the $L^2$ and $H^1$ norms, providing bounds like $\Vert f \Vert_{L^2}^2  \leq 2 f^2(0) + 4\Vert f_x \Vert_{L^2}^2 $. Agmon's inequalities offer bounds for the $L^\infty$ norm in terms of the $H^1$ norm and boundary values, such as $\Vert f \Vert_{L^\infty}^2 \leq 2 \Vert f \Vert_{L^2} \Vert f_x \Vert_{L^2} + f^2(0)$. As an example application of these inequalities, for functions $f \in H^1$ with $f(0) = 0$, the $H^1$ norm can be alternatively defined as simply $\Vert f \Vert_{H^1} = \Vert f_x \Vert_{L^2}$, maintaining equivalence with the traditional definition. 

Functional spaces and norms set the stage for defining stability in the infinite-dimensional setting.
Unlike in finite dimension, there is no general unifying PDE stability theory, even for \emph{linear} PDEs. For each particular family of problems, one has to develop different \emph{norm estimates} to show stability properties, and choosing the right norm or norms becomes essential. For instance, if $u(t,\cdot)$ is the state, a basic statement for the origin $u(t,x)\equiv 0$ being exponentially stable in $L^2$ would be
\begin{equation}
\Vert u(t,\cdot) \Vert_{L^2}  \leq C_1 \mathrm{e}^{-C_2t} \Vert u(0,\cdot) \Vert_{L^2},\label{eqn-expest}
\end{equation} 
The constants $C_1\geq1$ and $C_2>0$ in (\ref{eqn-expest}) are respectively referred to as the overshoot coefficient and the decay rate. Similar statements may be posed for other norms, and combined estimates are usually necessary for systems with more than one state. An example of how to derive estimates like (\ref{eqn-expest}) is presented next in Section~\ref{sec-plant}.

\subsection{Reaction-diffusion PDEs and their stability}\label{sec-plant}

Consider the following reaction-diffusion system, a type of PDE frequently encountered in chemical, thermal, and fluid applications:
\begin{equation}
u_t = \epsilon u_{xx} + \lambda(x) u, \quad u(t,0)=0, \quad u(t,1)=U(t). \label{eqn-reactdifint}
\end{equation}
where $u(t,x)$ represents the state, for $x \in [0,1]$ and $t > 0$, and $U(t)$ denotes a possible actuation variable. In~(\ref{eqn-reactdifint}), the diffusion coefficient $\epsilon$ is assumed to be positive, and the reaction term $\lambda(x)$ is considered continuous. We adopt a Dirichlet-type left boundary condition, $u(t,0)=0$, although Neumann ($u_x(t,0)=0$) or Robin ($u_x(t,0)=q u(t,0)$ for some $q$) boundary conditions can also be utilized. The Dirichlet actuation in Eq. (\ref{eqn-reactdifint}) could also alternatively be replaced with Neumann actuation. In this context, the most appropriate spaces\footnote{See~\cite{Evans2010} for details of how a linear 1-D parabolic PDE with Dirichlet boundary conditions is well-posed in $L^2$ and $H^1$. In this survey, the well-posedness of the PDEs under study is assumed, given that it represents a complex subject on its own.} to operate in are $L^2$ or $H^1$. Considering the former, assume initial conditions $u(0,x) = u_0(x)$ with $u_0 \in L^2$.

Consider first the open loop system by setting $U(t)=0$ in  (\ref{eqn-reactdifint}). Let $\bar \lambda$ represent the upper bound of $\lambda$ within the interval $[0,1]$. When $\bar \lambda < \frac{\epsilon}{4}$, the equilibrium state $u(t,x) \equiv 0$ of  (\ref{eqn-reactdifint}) can be shown to be exponentially stable in open loop, and an estimate similar to  (\ref{eqn-expest}) can be derived, as demonstrated next. Consider the following Lyapunov functional, equivalent to the squared $L^2$ norm:
\begin{equation}
L_1(t) = \frac{1}{2} \int_0^1 u^2  dx, \label{eqn-lyap}
\end{equation}
where the dependence of $u(t,x)$ on time and space is omitted for simplicity. Differentiating $L_1(t)$ with respect to time yields:
\begin{equation}
\frac{dL_1(t)}{dt} = \int_0^1 u u_t  dx = \epsilon \int_0^1 u u_{xx} dx + \int_0^1 \lambda u^2 dx, \label{eqn-lyap1}
\end{equation}
which, upon applying integration by parts and utilizing the open loop boundary conditions specified in  (\ref{eqn-reactdifint}) alongside Poincaré's inequality, leads to:
\begin{equation}
\frac{dL_1(t)}{dt} \leq -\frac{\epsilon}{4} \int_0^1 u^2 dx + \int_0^1 \lambda u^2 dx = -\frac{\epsilon - 4\bar\lambda}{2} L_1(t). \label{eqn-lyap4}
\end{equation}
Integration of $L_1(t)$ in  (\ref{eqn-lyap4}) results in:
\begin{equation}
L_1(t) \leq e^{-\left(\frac{\epsilon}{2} - 2\bar \lambda\right)t} L_1(0), \label{eqn-lyapbound}
\end{equation}
which, by leveraging the definition of $L_1(t)$, implies  (\ref{eqn-expest}) with a unity overshoot and decay rate of $\frac{\epsilon}{2} - 2\bar \lambda$. Specifically, a nonpositive $\bar \lambda$ guarantees the exponential stability of $u(t,x) \equiv 0$ for any diffusion coefficient.

\subsection{Backstepping transformations and target systems}\label{sect-transformation} 
The cornerstone of the backstepping method lies in the use of a coordinate transformation, termed the ``backstepping transformation,'' which redefines the states of a system into new states that exhibit desired characteristics. This transformation, alongside the resulting ``target system,'' the PDE that dictates the dynamics of the transformed states, are essential in any backstepping design. Our survey partly revolves around explaining the selection of transformations and target systems to tackle a variety of problems and applications within the realm of PDE control.

To illustrate these concepts, we start by selecting a target system for the reaction-diffusion equation of Section~\ref{sec-plant}. Given the potentially destabilizing reaction term in  (\ref{eqn-reactdifint}), a logical step in designing a target system is to either neutralize this term or transform it into a value that actively contributes to the system's stability, thereby enhancing the rate of convergence. Based on the results from the Lyapunov analysis in Section~\ref{sec-plant}, a judicious choice for a target system could be:
\begin{eqnarray}
w_t &=& \epsilon w_{xx} - cw, \quad w(t,0) = 0, \quad w(t,1) = 0, \label{eqn-target}
\end{eqnarray}
where $c \geq 0$ is a parameter determined by the designer. This selection ensures that the equilibrium $w(t,x) \equiv 0$ of (\ref{eqn-target}) is exponentially stable, with a decay rate of $\epsilon/2 + 2c$, showcasing how a decay rate can be fine-tuned by adjusting the parameters of the target system.

%The art of selecting ``the right'' target system is one of the principal challenges in backstepping. However, our example illustrates fundamental principles that generally hold true, albeit with possible exceptions. Specifically, the intrinsic nature of the original system, in this case, a reaction-diffusion equation, remains unaltered, preserving the diffusion coefficient; the Dirichlet boundary conditions are also maintained in the target system. Finally, the reaction is either canceled or turned negative. Thus, the target system looks very much like the original system, with just enough changes to become stable.

The art of selecting ``the right'' target system is one of the principal challenges in backstepping. Our example reveals several \emph{fundamental principles} that generally hold true, albeit with possible exceptions:
\begin{enumerate}
    \item The intrinsic nature of the original system must be preserved---in this case, a reaction-diffusion equation remains as such, maintaining its diffusion coefficient;
    \item The boundary conditions of the original system (here, Dirichlet conditions) are maintained in the target system;
    \item The destabilizing terms of the original system are transformed appropriately---the reaction term is either canceled or converted to a stabilizing (negative) term.
\end{enumerate}
Thus, the target system emerges as a carefully modified version of the original system---similar enough to maintain essential characteristics, yet different enough to achieve stability.

Continuing our example, we turn our attention to posing a transformation that maps system (\ref{eqn-reactdifint}) onto the selected target system (\ref{eqn-target}). To fully grasp the appropriate structure for such a transformation, it is instructive to revisit the origins of backstepping briefly. Initially developed for finite-dimensional \emph{strict-feedback} systems as explained in Section~\ref{sect-origins},  backstepping addressed configurations where a sequential chain of integrators ends with the control variable in the final equation. These systems display a ``lower-triangular'' structure, wherein the right-hand side terms of successive equations progressively depend on the current and all preceding states, thus establishing a cascading pattern. The successful application of backstepping to lower-triangular ODEs inspired its adaptation to PDEs, notably reaction-diffusion systems as exemplified by (\ref{eqn-reactdifint}), which manifest ``diagonal'' dynamics. In these dynamics, specific terms, such as the reaction term or the time derivative of the state $u_t$, are influenced solely by the local state $u(t,x)$, although the diffusion term extends this dependence to the infinitesimal neighborhood. Moreover, additional terms might emulate a strict-feedback structure, for instance, $g(x)u_x(0, t)$ and $\int_0^x f(x,\xi) u(t,\xi) d\xi$, thereby introducing a notion of ``spatial causality'' where the temporal evolution at any point $x$ is contingent on both the current state and those preceding it, i.e., $u(t,\xi)$ for $\xi \in [0,x]$. It should be noted, however, that this structure does not imply a finite speed of propagation; in a reaction-diffusion system, disturbances anywhere instantaneously influence the entire domain. This is analogous to a chain of integrators of any length, where perturbations instantaneously impact all states.

Inspired by this discussion, a strict-feedback (or spatially causal) transformation for our example is posed as:
\begin{equation}
w(t,x) = u(t,x) - \int_0^x k(x,\xi) u(t,\xi) d\xi, \label{eqn-backstepping}
\end{equation}
This transformation is henceforth designated as the backstepping transformation, with its integral kernel $k(x,\xi)$ termed the backstepping kernel. It is characterized by a Volterra-type integral and is therefore invertible under relatively lenient conditions on the kernel $k(x,\xi)$, a subject elaborated upon in Section~\ref{sec-inv}.

However, this represents merely one potential transformation for mappings between PDE systems. As the former discussion suggests, in the presence of non-strict-feedback terms, such a transformation might not suffice. An alternate ``full'' transformation
\begin{equation}
w(t,x) = u(t,x) - \int_0^1 k_F(x,\xi) u(t,\xi) d\xi, \label{eqn-fullbackstepping}
\end{equation}
could be viable in certain scenarios. This transformation, incorporating a Fredholm-type integral, presents more complex challenges regarding its invertibility~\cite{yoshida1960lectures}. The exploration of backstepping designs employing such transformations and the specific problems they address is further discussed in Section \ref{sec-fredholm}.

\subsection{Kernel equations: derivation and well-posedness}\label{sec-kerneleqs}
Once the backstepping transformation and target system have been chosen, the next step is to find the conditions that the backstepping kernel needs to verify, ensuring that the original system is mapped into the desired target system. These conditions take the form of (usually) hyperbolic PDEs, known as the kernel equations, and their well-posedness and solvability (or lack thereof) are what ultimately determine if the target system and transformation were correctly chosen.  
The way to derive the kernel equations is by substituting the transformation (\ref{eqn-backstepping}) into the target system (\ref{eqn-target}) and eliminating $w$. Then, integrating by parts, we get a weak formulation of a PDE, which can then be written in standard form. In the procedure, the (uncontrolled) boundary condition is also taken into account.

In particular, for our reaction-diffusion example with Dirichlet boundary conditions, the kernel in \eqref{eqn-backstepping} satisfies a hyperbolic PDE in the triangular domain $\mathcal{T}=\{(x,\xi):0\leq\xi\leq x\leq 1\}$, expressed as
\begin{eqnarray}
\epsilon k_{xx}(x,\xi) - \epsilon k_{\xi\xi}(x,\xi)&=&\left(\lambda(\xi)+c\right) k(x,\xi) , \label{eqn-kerneleq}
\end{eqnarray}
subject to boundary conditions
\begin{align}
 k(x,0) &= 0, \label{eqn-kerneleqbc1} \\
k(x,x) &= -\frac{1}{2\epsilon} \int_0^x \left(\lambda(\xi)+c\right) d\xi. \label{eqn-kerneleqbc2}
\end{align}
Equation (\ref{eqn-kerneleq})--(\ref{eqn-kerneleqbc2}) can be proven well-posed by first reducing it to an integral equation. By defining $\delta = x + \xi$ and $\eta = x - \xi$, and introducing $G(\delta,\eta) = k\left(\frac{\delta+\eta}{2},\frac{\delta-\eta}{2}\right)$, we are able to convert (\ref{eqn-kerneleq}) into an integral form:
\begin{eqnarray}
G(\delta,\eta)&=&
-\frac{1}{4\epsilon} \int_\eta^\delta\left(c+ \lambda \left(\frac{\tau}{2} \right)\right) d\tau
\nonumber \\ 
&&\hspace{-0.3cm}+\frac{1}{4\epsilon} \int_\eta^\delta \int_0^\eta \left(c+ \lambda\left(\frac{\tau-s}{2}\right)\right) G(\tau,s) ds d\tau. \label{eqn-intG}
\end{eqnarray}
The method of successive approximations can be used to solve (\ref{eqn-intG}).
Define
\begin{eqnarray}
G_0(\delta,\eta)
=
-\frac{1}{4\epsilon} \int_\eta^\delta \left(c+ \lambda \left(\frac{\tau}{2} \right) \right)d\tau,
\end{eqnarray}
and for $n\geq1$,
\begin{eqnarray}
G_n(\delta,\eta)&=&G_{n-1}(\delta,\eta)
\nonumber \\&&\hspace{-1cm}
+\frac{1}{4\epsilon} \int_\eta^\delta \int_0^\eta \left(c+ \lambda\left(\frac{\tau-s}{2}\right)\right) G_{n-1}(\tau,s) ds d\tau
.\label{eqn-Gn}
\end{eqnarray}
Then, it can be shown that, since $\lambda$ is continuous,
\begin{equation}
G(\delta,\eta)=\lim_{n\rightarrow \infty} G_n(\delta,\eta),
\end{equation}
defines a unique continuous solution for the integral form of the kernel equations, and therefore, there exists a  unique continuous $k(x,\xi)$ solving (\ref{eqn-kerneleq})--(\ref{eqn-kerneleqbc2}) in $\mathcal T$. This procedure also can be used to show that  $k\in\mathcal C^1(\mathcal T)$ given that $\lambda$ is continuous.

\subsection{Backstepping Feedback Gains and Control Laws}\label{sec-control}
Once the kernel equations are known to be well-posed, they define a valid backstepping transformation into the target system. However, the feedback law for the boundary controller has not been defined yet. This control law is derived via the backstepping transformation, specifically evaluated at the controlled boundary ($x=1$ in the example) and considering the boundary conditions of both original and target systems, given by equations~(\ref{eqn-reactdifint}) and (\ref{eqn-target}). The resulting expression for the feedback control law is:
\begin{equation}
    U(t) = \int_0^1 k(1,\xi) u(t,\xi)  d\xi. \label{eqn-feedback}
\end{equation}
Note that in the context of Neumann boundary control, analogous derivations yield the boundary control law:
\begin{eqnarray}
    u_x(t,1) &=& U(t) \nonumber \\
    &=& k(1,1)u(t,1) + \int_0^1 k_x(1,\xi) u(t,\xi)  d\xi. \label{eqn-feedbackN}
\end{eqnarray}
The implementation of feedback control laws (\ref{eqn-feedback}) and the Neumann boundary control law (\ref{eqn-feedbackN}) necessitates the precise determination of the boundary value (trace) of the gain kernel $k$ at $x=1$, this is, the backstepping kernel evaluated at $x=1$, or of its derivative $k_x$; it is noteworthy that the values of the kernel elsewhere do not play a role in the implementation of backstepping. The feedback control laws (\ref{eqn-feedback}) and \eqref{eqn-feedbackN} also need the full state of the system; this second requirement may be addressed by the use of observers (Section~\ref{sec-obs}) to design output-feedback control laws (Section~\ref{sec-output}).

\subsection{Invertibility of the transformation and closed loop stability}\label{sec-inv}
Once we have shown how the boundary feedback law (\ref{eqn-feedback}) is derived and computed, we study the question of stability of the closed-loop system (\ref{eqn-reactdifint}) with the proposed control law (\ref{eqn-feedback}). We are using the backstepping transformation~(\ref{eqn-backstepping}) to make the system ``behave'' like the target system (\ref{eqn-target}), whose exponential stability was already established in Section~\ref{sec-plant}. Hence, we need to relate the stability properties of the original $u$ system with those of the target $w$ system. That requires not only the backstepping transformation  (\ref{eqn-backstepping}) that maps $u$ into $w$ (which, for this reason, we will refer to as the ``direct'' transformation), but also its inverse, mapping $w$ into $u$. In the case of our example, due to the triangular form of the backstepping transformation, the theory of Volterra integral equations~\cite{yoshida1960lectures} guarantees that the direct transformation is indeed invertible, with the only condition that the kernel $k$ is at least bounded. Since the integral kernel equation (\ref{eqn-intG}) allowed us to show that $k\in\mathcal C^1(\mathcal T)$, the inverse transformation always exists. 

Define then the inverse transformation as
\begin{equation} 
u(t,x)=w(t,x)+\int_0^x l(x,\xi)w(t,\xi),
\end{equation}
where $l$ is known as the inverse kernel. The inverse and direct kernels can be shown to satisfy the following relationship
\begin{equation}
l(x,\xi)=k(x,\xi)+\int_\xi^x k(x,s) l(s,\xi) ds.\label{eqn-inversedirect}
\end{equation}
Equation (\ref{eqn-inversedirect}) does not depend on the original or target systems and allows to compute $l$ from $k$. An alternative approach is to derive a PDE equation for the $l$ kernel, which is done in a similar fashion to the $k$ kernel equation. Care should be taken in more general designs (e.g., those including extra degrees of freedom in the kernel equations) to make sure the $l$ kernel thus obtained is indeed the inverse of the $k$ kernel. Note, however, that (at least in this design) the inverse transformation is not used in the implementation of backstepping; only its existence and properties are required for the proof of closed-loop stability. In particular, the norm of this inverse transformation appears in the
stability estimate as seen below in (\ref{eqn-ul2est}).

Using the direct and inverse backstepping transformations, respectively (\ref{eqn-backstepping}) and (\ref{eqn-inversedirect}), we can use them to compute an estimate of the $L^2$ norm of $w$ in terms of the $L^2$ norm of $u$, and vice-versa, as follows:
\begin{eqnarray}
\Vert w\Vert_{L^2}^2&\leq&\left(1+\Vert k \Vert_{L^\infty(\mathcal T)}\right)^2 \Vert u\Vert^2_{L^2},\label{eqn-l2estdirect}\\
\Vert u\Vert_{L^2}^2&\leq&\left(1+\Vert l \Vert_{L^\infty(\mathcal T)}\right)^2 \Vert w\Vert^2_{L^2}.\label{eqn-l2estinverse}
\end{eqnarray}
Since $k,l\in\mathcal C^1(\mathcal T)$ and $L^\infty \subset \mathcal C^1(\mathcal T)$ because $\mathcal T$ is a bounded domain, the $L^\infty(\mathcal T)$ norm in (\ref{eqn-l2estdirect})--(\ref{eqn-l2estinverse}) is well-defined for $k$ and $l$.
Using the estimate derived in Section~\ref{sec-plant} for a reaction-diffusion equation, in the case that $\lambda=-c$, we get an estimate of exponential stability of $w$, namely that
\begin{equation}
\Vert w(t) \Vert_{L^2}^2 \leq \mathrm{e}^{-(\epsilon/2+c/4)t} \Vert w(0) \Vert_{L^2}^2,
\end{equation}
hence using (\ref{eqn-l2estdirect})--(\ref{eqn-l2estinverse}) we obtain 
\begin{eqnarray}
\Vert u(t) \Vert_{L^2}^2 &\leq& \left(1+\Vert l \Vert_{L^\infty(\mathcal T)}\right)^2\left(1+\Vert k \Vert_{L^\infty(\mathcal T)}\right)^2 \nonumber \\ && \times  \mathrm{e}^{-(\epsilon/2+c/4)t} \Vert u(0) \Vert_{L^2}^2,\label{eqn-ul2est}
\end{eqnarray}
thus showing exponential stability for $u$. As typical in backstepping, one can select the decay rate, but the overshoot coefficient depends on the sizes of the kernels (and will, in general, grow if the decay rate is selected large). Similarly, $H^1$ exponential stability can be derived with more involved computations.

Beyond stability analysis, the backstepping approach can also offer additional benefits related to optimality. As shown in~\cite{Smyshlyaev2004}, certain backstepping designs solve an inverse optimal stabilization problem, where the controller not only stabilizes the system but also minimizes a meaningful cost functional. This inverse optimality property establishes an important connection between backstepping and optimal control approaches.

        \subsection{Observer design and duality}\label{sec-obs}
Section~\ref{sect-transformation} to Section~\ref{sec-inv} encompass the typical full-state feedback backstepping design, from the selection of the target system and transformation to the final proof of stability. However, these control laws are not implementable in practice since they require knowledge of the state at (almost) every point of the interval $(0,1)$. To overcome this issue, backstepping also provides a design procedure for observers that is dual to the controller design methodology, as illustrated next for our example. 

Consider now that the only available measurement of the state $u$ in the system (\ref{eqn-reactdifint})  is $u_x(t,0)$. This setup is known in the literature as anti-collocated since the measurement is done at the uncontrolled end, in contrast to collocated setups where both measurement and actuation happen at the same boundary. Only the former is detailed here for conciseness. 

To obtain an estimate of the state that can be used in control law (\ref{eqn-feedback}), we postulate an observer as a copy of the system plus injection of output error. Denoting the estimates with a hat, the observer equations are
\begin{eqnarray} 
\hat u_t&=&\epsilon \hat u_{xx}+\lambda(x) \hat u+p_1(x) \left(u_x(t,0)-\hat u_x(t,0) \right)
,\label{eqn-obseq}\\
\hat u(t,0)&=&0, \quad
\hat u(t,1)=U(t),\label{eqn-obseqbc2}
\end{eqnarray}
where $p_1(x)$ is an output injection gain to be determined.

Thanks to linearity, the observer error $\tilde u=u-\hat u$ verifies the following autonomous PDE
\begin{eqnarray} 
\tilde u_t&=&\epsilon \tilde u_{xx}+\lambda(x) \tilde u -p_1(x) \tilde u_x(t,0) 
,\label{eqn-obserrq}\\
\tilde u(t,0)&=&0,\quad
\tilde u(t,1)=0,\label{eqn-obserrqbc2}
\end{eqnarray}
The gain $p_1$ is designed with backstepping to obtain exponential stability of the origin $\hat u\equiv 0$, thus guaranteeing convergence of the estimate.

Following the same methodology of Section~\ref{sect-transformation} we map the $\tilde u$ variable into a target state using the transformation
\begin{eqnarray}
\tilde u(t,x)=\tilde w(t,x)-\int_0^x p(x,\xi) \tilde w(t,\xi).\label{eqn-backobs}
\end{eqnarray}
It maps the (\ref{eqn-obserrq})--(\ref{eqn-obserrqbc2}) system into the  $\tilde w$ target system,
\begin{eqnarray}
\tilde w_t&=&\epsilon \tilde w_{xx}-c \tilde w,\quad
\tilde w(t,0)=0,\quad 
\tilde w(t,1)=0,\label{eqn-obstargeteq}
\end{eqnarray}
where $c \geq 0$ is chosen (possibly with a different value than for the controller) by the designer. As seen in Section~\ref{sec-plant} this selection ensures that the equilibrium $\tilde w(t,x) \equiv 0$ of (\ref{eqn-target}) is exponentially stable, with a decay rate of $\epsilon/2 + 2c$, thus guaranteeing convergence of the estimate $\hat u$ to the real value of $u$.

As in Section~\ref{sect-transformation}, from the original system (\ref{eqn-obserrq})--(\ref{eqn-obserrqbc2}), the target system equations (\ref{eqn-obstargeteq}) and the transformation~(\ref{eqn-backobs}), we obtain the $p$-kernel equations defined in the domain $\mathcal T$ (see Section~\ref{sec-kerneleqs}) by
\begin{eqnarray}
\epsilon p_{\xi \xi}&=&\epsilon p_{xx}+\lambda(x) p(x,\xi), \label{eqn-obskerneleq}\\
p(x,x)&=&\frac{1}{2\epsilon} \int_0^x \lambda(\xi) d\xi,\label{eqn-obskerneleqbc1}\\
p(1,\xi)&=&0.\label{eqn-obskerneleqbc2}
\end{eqnarray}
In addition, the following condition should be satisfied: 
\begin{eqnarray}
p_1(x)&=&-\epsilon p(x,0).
\end{eqnarray}
Hence, the gain $p_1$ of the observer is determined from the observer transformation kernel, as the controller gain was derived from the control transformation kernel.

It remains to show that (\ref{eqn-obskerneleq})--(\ref{eqn-obskerneleqbc2}) is well-posed. To show it, we define new variables $\check x=1-\xi$, $\check \xi=1-x$, $\check p(\check x,\check \xi)=p(x,\xi)$ and $\check\lambda(\xi)=\lambda(x)$. Then, (\ref{eqn-obskerneleq})--(\ref{eqn-obskerneleqbc2}) written in terms of the new variables is
\begin{eqnarray}
\epsilon \check p_{\check x \check x}&=&\epsilon \check p_{\check \xi \check \xi}+\check \lambda(\check \xi) \check p(\check x,\check \xi), \label{eqn-chobskerneleq}\\
\check p(\check x,\check x)&=&\frac{1}{2\epsilon} \int_0^{\check x} \check \lambda(\xi) d\xi,\label{eqn-chobskerneleqbc1}\\
\check p(\check x,0)&=&0.\label{eqn-chobskerneleqbc2}
\end{eqnarray}
Notice that (\ref{eqn-chobskerneleq})--(\ref{eqn-chobskerneleqbc2}) are the same as the kernel equations we obtained for $k$ in Section~\ref{sec-kerneleqs} (see equations (\ref{eqn-kerneleq})--(\ref{eqn-kerneleqbc2})). Hence, using the same method that consists of transforming the system into an integral equation and solving it by successive approximation, we get that, given that $\lambda$ is continuous, there is a unique $p\in \mathcal C^1(\mathcal T)$. As in Section~\ref{sec-inv}, we can define an inverse transformation (with a well-posed kernel) and use it to prove the exponential stability of the origin for the $\tilde u$ system. Hence, we get that the estimate $\hat u$ converges (exponentially) to the state $u$ in the $L^2$ norm (also in the $H^1$ norm).

\subsection{Output feedback designs}\label{sec-output}
        
Even though the observer design has interesting applications by itself, such as state estimation in lithium-ion batteries \cite{moura2013batteries,moura2017batteries}, it can also be used in the construction of output-feedback controllers. In this case, the full-state feedback law of Section~\ref{sec-control} and the observer designed in Section~\ref{sec-obs} can be combined to obtain a stabilizing output feedback as follows:
\begin{eqnarray}
U(t)=\int_0^1 k(1,\xi) \hat u(x,\xi),\label{eqn-coutputfb}
\end{eqnarray}
where $\hat u$ is the observed state obtained from (\ref{eqn-obseq})--(\ref{eqn-obseqbc2}) and the controller and observer gains, $k(1,\xi)$ in (\ref{eqn-coutputfb}) and $p(x,0)$ in (\ref{eqn-obseq}) are computed from the kernels that solve the controller and observer kernel equations, respectively equations (\ref{eqn-kerneleq})--(\ref{eqn-kerneleqbc2}) and (\ref{eqn-obskerneleq})--(\ref{eqn-obskerneleqbc2}).

Notice that the controller has become dynamic since a PDE has to be solved to compute the control law. Thus, in studying the stability of the closed loop, the observer state needs also to be considered. 
For a linear system such as (\ref{eqn-reactdifint}), the separation principle (or ``certainty equivalence'') holds, and it can be shown that the output feedback controller defined by (\ref{eqn-coutputfb}) together with (\ref{eqn-obskerneleq})--(\ref{eqn-obskerneleqbc2}) makes the origin of the combined plant-observer system exponentially stable. However, the result is subtle and needs to be stated for rigor in the $H^2$ space since the observer is driven by $u_x(t,0)$, which is not well-defined for $u\in L^2$ or $u \in H^1$; no further details are provided, as it goes beyond the purpose of this survey, but the interested reader may consult \cite{qi2015multi} for an example of such a result.
\subsection{Extensions to other one-dimensional single-variable systems}
After the initial results for parabolic equations were obtained, it was quickly recognized that the methodology presented in the previous sections could be extended to various other PDE systems. Here, several important extensions focusing on 1D systems described by a single PDE are explored. We first discuss designs for a first-order hyperbolic PDE, which models, for example, traffic flows. We then present the extensions to the wave equation and flexible structure models (shear beam and Euler-Bernoulli beam). We mention how backstepping can be applied to complex-valued PDEs such as the Schrödinger and Ginzburg-Landau equations. We finish with an example where backstepping is used to prove null controllability for a Korteweg-de Vries equation.

\subsubsection{First-order hyperbolic PDEs}
First-order hyperbolic PDEs, such as the transport equation and its variants, describe systems with unidirectional propagation and finite speed of information transfer. They find applications in diverse areas, including traffic flow modeling, chemical reactor design, and transmission line analysis.  Unlike second-order hyperbolic PDEs, which exhibit oscillatory behavior, first-order hyperbolic PDEs display a ``delay-like'' response. The control of these systems focuses on regulating downstream states by manipulating upstream inputs.

In \cite{Krstic2008c}, Krstic and Smyshlyaev extended the backstepping method to solve the stabilization problem for a class of first-order hyperbolic PDEs. The target systems in this design are chosen to be pure delays, ensuring convergence to zero in finite time.

As an example, consider a first-order hyperbolic PDE with a potentially destabilizing integral and boundary coupling term:
\begin{equation}
u_t(t,x) = u_x(t,x) + g(x)u(t,0) + \int_0^x f(x,\xi) u(t,\xi) d\xi \label{eqn-fohyp}
\end{equation}
with the control input $u(t,1)=U(t)$. The backstepping transformation 
\begin{equation}
w(t,x) = u(t,x) - \int_0^x k(x,\xi)u(t,\xi) d\xi
\label{eqn-fohyp_transform}
\end{equation}
is used to map (\ref{eqn-fohyp}) into the delay line target system:
\begin{equation}
w_t(t,x) = w_x(t,x), \quad w(t,1) = 0.
\label{eqn-fohyp_target}
\end{equation}
The solution to (\ref{eqn-fohyp_target}) is $w(t,x) = w(t+x-1,1) = 0$ for $t+x>1$. Hence the target system converges to zero in finite time. Using (\ref{eqn-fohyp}), (\ref{eqn-fohyp_transform}), and (\ref{eqn-fohyp_target}), we obtain the following kernel equations:
\begin{align}
k_x(x,y) + k_y(x,y) &= \int_y^x k(x,\xi) f(\xi,y)d\xi - f(x,y),\\
k(x,0) &= \int_0^x k(x,y)g(y)dy - g(x).
\end{align}
These equations can be shown to be well-posed using the method of successive approximations as in Section~\ref{sec-kerneleqs}. Evaluating the transformation (\ref{eqn-fohyp_transform}) at $x=1$ gives the control law:
\begin{equation}
u(t,1) = \int_0^1 k(1,y)u(t,y) dy.
\end{equation}
The finite-time stability of the closed-loop system is then obtained from the finite-time stability of the target system. 

One of the interesting applications of the design in \cite{Krstic2008c} is the control of finite-dimensional systems with actuator and sensor delays.  A system with input delay can be represented as a cascade of an ODE and a transport PDE (which models the delay). Combining the backstepping design for hyperbolic PDEs with the backstepping design for linear ODEs allows recovering well-known infinite-dimensional controllers for systems with delays, such as the Smith predictor for unstable plants. This application is further discussed in Section~\ref{Sec_ODE_PDE}.

\subsubsection{Wave equation and flexible structures}
The wave equation and flexible structure models such as strings, cables, and beams are commonly described by second-order or fourth-order hyperbolic PDEs, whose dynamics often result in oscillatory behaviors. The control of these systems is crucial in numerous applications, including noise suppression, vibration damping, and structural stabilization.

The extension of backstepping to wave equations and flexible structure models proved initially challenging due to the lack of inherent damping in their dynamics. While transformations to cancel destabilizing terms in the equation could be found,  it was not clear how to ensure the stability of the resulting target systems.  The crucial step in extending backstepping to these systems involved identifying attainable target systems that incorporated damping terms. 

Thus, a problem for a wave equation with anti-damping on the uncontrolled boundary was addressed in \cite{Smyshlyaev2009}. In this scenario, the anti-damping term acted through the boundary condition, making the system unstable with infinitely many eigenvalues in the right-half plane. Another example considers an output-feedback problem with instability at its free end and control on the opposite end~\cite{Krstic2008a}. A similar problem of boundary stabilization for a one-dimensional wave equation with an internal spatially varying antidamping term was tackled in \cite{Smyshlyaev2010}. This antidamping term made the uncontrolled system unstable, with all of its eigenvalues located in the right-half plane. The proposed design used a novel backstepping transformation with a 2x2 structure, mapping the unstable system into a stable target system with a desired decay rate. The control gains were expressed in closed form for systems with constant parameters.  Other backstepping designs on the wave equation include \cite{Zhou2012,wu2020output}.  Although these backstepping controllers may have limited practical applicability due to their poor delay-robustness margins, resulting from the infinite number of unstable poles in open-loop, the underlying theoretical developments have been essential for the development of the backstepping methodology. 

For beams, a pioneering work in this area was~\cite{Krstic2006}, where boundary controllers and observers were designed for the wave equation and the Timoshenko beam model with Kelvin-Voigt damping.  In contrast to the reaction-diffusion equation, where the potentially destabilizing reaction term can be simply eliminated in some cases, the target systems in \cite{Krstic2006} contained damping terms arising from Kelvin-Voigt damping (a type of material damping that introduces a strain-rate term in the constitutive equation). This contribution laid the foundation for subsequent studies on the shear beam~\cite{krstic2008beam} and the Euler–Bernoulli beam~\cite{Smyshlyaev2009}. 

A common feature of these designs is that the damping term appears in the target system, even though the original model lacks inherent damping. The introduction of this damping term in the target system allows moving the closed-loop eigenvalues to the left in the complex plane, achieving an arbitrary decay rate. This ability to prescribe the decay rate represents a significant improvement over classical ``boundary damper'' controllers, which typically provide limited damping for wave or beam models.

Note that some recent results on the wave equation and Timoshenko beams exploit their representation as PDE-ODE systems through a Riemann transformation and are explored in Section~\ref{sec-tim}.
\subsubsection{Complex-valued equations}

While most of the early backstepping designs focused on real-valued PDEs, the methodology readily extends to complex-valued equations, encompassing models such as the Schrödinger and Ginzburg-Landau equations. These equations often arise in quantum mechanics, nonlinear optics, and fluid dynamics, and their control presents unique challenges due to the complex nature of their states and dynamics. 

The earliest of these designs tackles the Ginzburg-Landau equation, a complex-valued PDE that models vortex shedding in fluid flows. In \cite{Aamo2005}, boundary controllers are designed leveraging the backstepping method with complex coefficients, noticing that most of the steps shown in previous sections hold in the complex case. This control design has direct relevance to flow control applications, offering a systematic approach for attenuating vortex shedding and enhancing flow stability.

In \cite{krstic2011boundary}, the problem of boundary control for the linearized Schrödinger equation was tackled with backstepping. The key idea was to formally view the Schrödinger equation as a complex-valued heat equation and apply the backstepping method developed for parabolic PDEs. This resulted in a full-state feedback controller that achieved an arbitrary decay rate. The control gains were expressed in closed form, involving Kelvin functions, highlighting the explicit nature of the backstepping design. 

As an illustrative example, consider the linearized Schrödinger equation:
\begin{equation}
\psi_t(t,x) = -j \psi_{xx}(t,x), \quad \psi_x(t,0) = 0, \quad \psi(t,1) = U(t), \label{eqn-schro}
\end{equation}
where $\psi$ is a complex-valued state and $j$ is the imaginary unit.  
The usual Volterra backstepping transformation is used to map (\ref{eqn-schro}) into the target system
\begin{equation}
\phi_t = -j \phi_{xx} - c \phi, \quad \phi_x(t,0) = 0, \quad \phi(t,1) = 0. \label{eqn-schro_target}
\end{equation}
The target system has eigenvalues with real parts equal to $-c$, and therefore the constant $c$ can be used to set the desired decay rate.
Using (\ref{eqn-schro}) and (\ref{eqn-schro_target}), we obtain the following kernel equations:
\begin{align}
k_{xx}(x,y) - k_{yy}(x,y) &= c j k(x,y),\\
k_y(x,0) &= 0,\\
k(x,x) &= -\frac{cj}{2}x.
\end{align}
This equation can be solved as in the real case in closed form to obtain:
\begin{equation}
    k(x, y) = -cjx \frac{I_1\left(\sqrt{cj(x^2 - y^2)}\right)}{\sqrt{cj(x^2 - y^2)}},
\end{equation}
which can be written in terms of  Kelvin functions (see~\cite{krstic2011boundary}) and then used as usual to obtain the control law; exponential stability of the closed-loop system is as before obtained from the stability of the target system and the invertibility of the transformation.

Finally, it must be noted that complex-valued kernels are always obtained when working in higher dimensional geometries and applying Fourier series or Fourier transform methods; see Section~\ref{sec-higher} for more details.

\subsubsection{Korteweg-de Vries equation}
The vast majority of kernel equations encountered in backstepping designs share common structural properties and can be solved using standard techniques. However, some problems lead to unusual kernel equations while still yielding well-posed backstepping designs. A notable example appears in~\cite{xiang2019null}, a work that considers a linearized Korteweg-de Vries equation. This result stands out for two reasons. First, backstepping is used not to design a controller but as a tool to prove null-controllability. Second, despite presenting highly unusual kernel equations that depart significantly from standard forms, the equations remain well-posed. Such examples, where backstepping succeeds despite atypical structural properties, hint perhaps at a broader theoretical foundation that is not yet fully understood.

    \section{Interconnected %\modifJA{hyperbolic} 
    systems} \label{Sec_hyp_systems}

    In this section, we introduce the rich theory of backstepping for interconnected hyperbolic and parabolic systems, a domain that has seen significant theoretical advancements primarily in the last decade. This burgeoning development is a key driver behind the rapid expansion and growing interest in the field. 
    
    Coupled first-order linear hyperbolic systems, potentially interconnected with ODEs, naturally arise in the mathematical description of transport phenomena with finite propagation speeds. Examples include the transport of matter, sound waves, and information. These systems are essential for describing a wide variety of large complex systems, including wave propagation, traffic network systems, electric transmission lines, hydraulic channels, drilling devices, communication networks, smart structures, and multiscale and multiphysics systems~\cite{bastin2016stability,Yu2023}. The backstepping approach has enabled significant breakthroughs for the stabilization of such systems. To present the wide variety of possible PDE-ODE interconnections analyzed in the literature,  we introduce the following general formulation 
\begin{equation}\label{Sec_4_1_v_pde_orig_simp}
\left\{
\begin{array}{l}
\dot{X}(t)=A_0X(t)+E_0v(t,0), \\
    w_t+\Lambda(x)w_x=\Sigma(x)w,\\
    \dot{Y}(t)=A_1Y(t)+E_1u(t,1),\\
    u(0)=Qv(0)+C_0X,~v(1)=R u(1)+C_1Y,
\end{array} 
\right.
\end{equation}
where $x\in [0,1]$ and $t>0$, $w=(u^\top,v^\top)^\top$. The function $u(t,x) \in \mathbb{R}^n$ represents $n$ rightward propagating states (information that travel from the left boundary of the system to the right boundary) while the function $v(t,x) \in \mathbb{R}^m$ represents $m$ leftward propagating states (information that travel from the right boundary to the left boundary). The functions $X$ and $Y$ are finite-dimensional states belonging to $\mathbb{R}^p$ and $\mathbb{R}^q$. The velocity matrix is $ \Lambda(x)=\text{diag}(\Lambda^+(x),$ $-\Lambda^-(x))$ where $\Lambda^+=\text{diag}(\lambda_1,\cdots,\lambda_n)$ and $\Lambda^-=\text{diag}(\mu_1,\cdots,\mu_m)$ are diagonal submatrices and 
$$
-\mu_m(x)<\cdots<-\mu_1(x)<0<\lambda_1(x)<\cdots< \lambda_n(x),
$$
for all $x\in [0,1]$. If some velocities are equal, then the velocity matrix presents isotachic blocks (these are blocks of states sharing the same transport speed). This specific case is exposed in Section~\ref{Sec_isotachic}. The functions $\mu_j$ and $\lambda_i$ are continuously differentiable. The in-domain coupling term $\Sigma$ is defined as
$\Sigma(x)=\begin{pmatrix}
    \Sigma^{++}(x) & \Sigma^{+-}(x)\\ \Sigma^{-+}(x) & \Sigma^{--}(x)
\end{pmatrix}$, where $\Sigma^{++}(x) \in \mathbb{R}^{n\times n}$, $\Sigma^{+-}(x) \in \mathbb{R}^{n\times m}$, $\Sigma^{-+}(x) \in \mathbb{R}^{m\times n}$ and $\Sigma^{--}(x) \in \mathbb{R}^{m\times m}$ are continuous functions.
Applying the exponential change of coordinates presented in~\cite{Coron2013}, it is possible to consider that the matrices $\Sigma^{++}$ and $\Sigma^{--}$ do not have diagonal terms as there always exist canonical coordinates which allow transforming the system into such a canonical form. Most contributions in the literature consider the velocity matrix $\Lambda$ and the coupling matrix $\Sigma$ as constant to simplify the exposition, although the different computations and results can be easily extended to spatially varying matrices. We adopt this choice in the rest of this section. Finally, the boundary couplings~$Q$ and~$R$ are constant and are respectively called \textbf{distal reflection} (reflection at the unactuated boundary) and \textbf{proximal reflection} (reflection at the actuated boundary). The matrices $A_0 \in \mathbb{R}^{p\times p}, E_0 \in \mathbb{R}^{p\times m}, A_1 \in \mathbb{R}^{q\times q}, E_1 \in \mathbb{R}^{q\times n}, C_0 \in \mathbb{R}^{n\times p}$ and $C_1 \in \mathbb{R}^{m\times q}$ are constant. Depending on the configurations considered in the different subsections, some coupling terms will be set to zero. We will also specify the control input and the available measurement. %\modifJA{To add what is commented below?} 
%The initial condition is denoted~$(X_0,u_0,v_0,Y_0)$, where $X_0 \in \mathbb{R}^{p}, Y_0 \in \mathbb{R}^q$ and where $u_0$ and $v_0$ are usually chosen in $L^2$ or $H^1$ (with appropriate compatibility conditions in this latter case)~\cite{bastin2016stability}.

Most (if not all) of the contributions presented in this section are based on backstepping transformations that generalize the one introduced in ~\eqref{eqn-backstepping}. Denoting $(\bar X, \gamma, \bar Y)$ the target system state (where $\gamma=(\alpha^\top, \beta^\top)^\top$), we have 
% \begin{equation}\label{Sec_4_BS_transf}
% \left\{
% \begin{array}{l}
% \bar X=X-\int_0^1 K_X w(\xi)d\xi +F_XY,\\
% \gamma(x) = u(x) - \int_0^x K(x,\xi) w(\xi) d\xi+F(x)Y, \\
% \bar Y=Y.
% \end{array} 
% \right.
% \end{equation}
% which can be expressed in the compact operator form
\begin{align}\label{Sec_4_BS_transf_matrix_form}
    \begin{pmatrix}
        \bar X(t) \\ \gamma(t,x) \\ \bar Y(t)
    \end{pmatrix}=  \begin{pmatrix}
         \mathbb{I}_{X} & F_w & F_Y\\ 
         0 & \mathbb{V} & F(x) \\
         0 & 0 &  \mathbb{I}_{Y} 
    \end{pmatrix}\begin{pmatrix}
         X(t) \\ w(t,x) \\ Y(t)
    \end{pmatrix},
\end{align}
where $\mathbb{I}_{\lbrace X,Y\rbrace}$ are identity operators of adequate size (that can reasonably be extended to bounded and boundedly invertible operators), $F_w$ is a bounded operator, $F_Y$ is a matrix and $F$ is a bounded operator (in most cases, a family of matrices parametrized by the spatial variable). Finally, the operator $\mathbb{V}$  is chosen as a bounded and \emph{boundedly invertible} integral operator. It is usually a Volterra integral operator of the second kind (as in ~\eqref{eqn-backstepping}). However, as it will appear throughout this section, several configurations require adding affine terms (with a diagonal structure) or even considering Fredholm integral operators. In this latter case, the invertibility of the transformation is not granted anymore and should be carefully proved.  Due to its block-triangular structure and the invertibility of the blocks on the diagonal, the transformation~\eqref{Sec_4_BS_transf_matrix_form} is boundedly invertible. Consequently, as explained in Section~\ref{sec-inv}, the original system and the target system share equivalent stability properties.

   \subsection{Coupled Hyperbolic systems}  \label{Sec_hyp_system}

   \subsubsection{$2\times 2$ hyperbolic systems} \label{Sec_hyp_2_2}
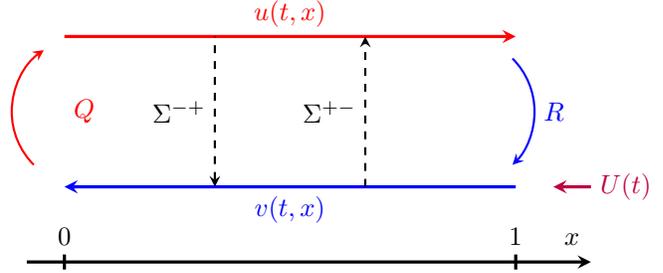
\begin{figure}[htb]%
\begin{center}
\begin{tikzpicture}
\draw [>=stealth,->,red,very thick] (0,0) -- (6,0);
\draw [red] (3,0) node[above]{$u(t,x)$};

\draw [>=stealth,<-,blue,very thick] (0,-2) -- (6,-2);
\draw [blue] (3,-2) node[below]{$v(t,x)$};

\draw [>=stealth,<-,dashed, thick] (2,-2) -- (2,0);
\draw(2,-1) node[left]{$\Sigma^{-+}$};

\draw [>=stealth,->,dashed, thick] (4,-2) -- (4,0);
\draw(4,-1) node[left]{$\Sigma^{+-}$};

\draw [red,>=stealth, thick](-0.7,-1) arc (-180:-135:1);
\draw [red,>=stealth,->, thick](-0.7,-1) arc (-180:-235:1);
\draw [red] (0,-1) node[right]{$Q$};

\draw [blue,>=stealth, thick](6.25,-1) arc (0:45:1);
\draw [blue,>=stealth,->, thick](6.25,-1) arc (0:-45:1);
\draw [blue] (6.25,-1) node[right]{$R$};

\draw [>=stealth,<-,purple,very thick] (6.5,-2) -- (7,-2);
\draw [purple] (7,-2) node[right]{$U(t)$};

%\draw [>=stealth,->,purple,very thick] (-1,0) -- (-0.5,0);
%\draw [purple] (-1,0) node[left]{$U(t)$};

\draw [>=stealth,->,very thick] (-0.5,-3) -- (7,-3);
\draw [very thick] (0,-3.1) -- (0,-2.9);
\draw (0,-2.9) node[above]{0};
\draw [very thick] (6,-3.1) -- (6,-2.9);
\draw (6,-2.9) node[above]{1};
\draw (6.75,-2.9) node[above]{$x$};
\end{tikzpicture}
\caption{Schematic representation of the system~\eqref{Sec_4_1_v_pde_orig_simp} in the absence of ODE states for $n=m=1$}
    \label{fig:orig_classic}
\end{center}
\end{figure}

The key ideas for the backstepping stabilization of coupled hyperbolic systems have been introduced in~\cite{Coron2013,vazquez2011backstepping} for two coupled scalar equations. Consider the system~\eqref{Sec_4_1_v_pde_orig_simp} with $n=m=1$ and without ODEs. This is the more general form for a one-dimensional $2\times 2$ hyperbolic linear system, without including integral or boundary terms~\cite{bastin2016stability}. We assume $Q\ne 0$. A schematic representation of this system is given in Figure~\ref{fig:orig_classic}.  The open-loop stability can be characterized using Lyapunov-based conditions~\cite{bastin2016stability,coron2008dissipative}. Due to the in-domain coupling terms $\Sigma^{+-}$ and $\Sigma^{-+}$ and the boundary coupling terms $Q$ and $R$, the open-loop system may be unstable.

We assume that the control law $U(t)\in \mathbb{R}$ acts at the $x=1$ boundary, i.e. we have $v(t,1)=Ru(t,1)+U(t)$.
%A stabilizing controller is designed in~\cite{Coron2013} using the backstepping approach. 
Given the potentially destabilizing coupling terms in ~\eqref{Sec_4_1_v_pde_orig_simp}, a logical step in designing the backstepping target system is to neutralize these terms. More precisely, in~\cite{Coron2013}, the target system $(\alpha,\beta)$ is considered as a copy of the initial dynamics~\eqref{Sec_4_1_v_pde_orig_simp} but with $\Sigma=0$ and $R=0$. 
%Denoting $\gamma=(\alpha,\beta)^\top$ the target-system state,
% \modifJA{its
%  exponential stability is shown in~\cite{Coron2013} by analyzing the following Lyapunov function equivalent to the squared $L^2$ norm: 
% \begin{align}
%     L(t)=\int_0^1\frac{\mathrm{e}^{-\nu x}}{\lambda(x)}\alpha^2(t,x)+a\frac{\mathrm{e}^{\nu x}}{\mu(x)}\beta^2(t,x) dx, \label{Lyap_hyper}
% \end{align}
% where $a$ and $\nu$ are two positive constants.}
This target system is finite-time stable. Indeed, applying the method of characteristics, we have that for all $t\geq t_F=\frac{1}{\lambda_1}+\frac{1}{\mu_1}$,
\begin{align}
    \alpha(t,x)=\beta(t,x)=0.
\end{align}
The finite-time $t_F$ corresponds to the minimum stabilization time given in~\cite{li2010strong}.

As shown in~\cite{Coron2013}, the original system can be mapped to the desired target system using the backstepping transformation~\eqref{Sec_4_BS_transf_matrix_form} (without the $X$ and $Y$ states and where $\mathbb{V}$ is a Volterra integral operator of the second kind). The kernel $K$ satisfies a set of hyperbolic PDEs that can be rewritten as integral equations using the method of characteristics. Then, applying the method of successive approximations, it is possible to show the existence of a unique fixed point and, consequently, a unique continuous solution to these kernel equations~\cite[Theorem A.1]{Coron2013}. Higher regularity of the solutions can also be shown under adequate regularity of the coefficients appearing in the kernel PDEs~\cite[Theorem A.2]{Coron2013}.

% the following hyperbolic PDE
% \begin{align}\label{Sec_2_1_eq:kernel_equations}
% \Lambda(x) k_x(x,y)+(k(x,y)\Lambda(y))_y+k(x,y)\Sigma(y)=0,
% \end{align}
% \normalsize
% with the boundary condition
% \begin{align}
%     & \Lambda k(x,x)-k(x,x)\Lambda=-\Sigma(x), \\
%      &k(x,0)\Lambda(0)\begin{pmatrix}
%          0 & Q \\ 0 & \text{Id}
%      \end{pmatrix}=0.\label{Sec_2_1_kernel_BC}
% \end{align}
% When $Q\ne 0$, it is possible to show the existence of a unique continuous solution to equations~\eqref{Sec_2_1_eq:kernel_equations}-\eqref{Sec_2_1_kernel_BC}~\cite{Coron2013}. The proof consists in writing~\eqref{Sec_2_1_eq:kernel_equations}-\eqref{Sec_2_1_kernel_BC} as integral equations using the method of characteristics. Then, applying the method of successive approximations, it is possible to show the existence of a unique fixed point and consequently a unique continuous solution to~\eqref{Sec_2_1_eq:kernel_equations}-\eqref{Sec_2_1_kernel_BC}. 
This leads to the following state-feedback controller
\begin{align} \label{Sec_4_1_1_control}
    U(t)=-Ru(t,1)+\int_0^1 \begin{pmatrix}
        0 \\ Id
    \end{pmatrix}^TK(1,y)w(y)dy,
\end{align}
which guarantees the exponential stability of the closed-loop system (and even its finite-time stability). The case $Q=0$ is also solved in~\cite{Coron2013} using a slightly different target system. %However, it requires a slightly different approach (with a different target system), which we do not detail here. 

A collocated state-observer is proposed in~\cite{vazquez2011backstepping} using the measurement $y(t)=v(t,1)$. This observer is defined as a copy of the system plus injection of output error tuned using the backstepping approach to guarantee the exponential stability of the error system (see Section~\ref{sec-obs}). This observer is combined with the state feedback controller~\eqref{Sec_4_1_1_control}  in~\cite{vazquez2011backstepping} to obtain an output-feedback controller that stabilizes the original system~\eqref{Sec_4_1_v_pde_orig_simp} in finite time $2t_F=\frac{2}{\lambda_1}+\frac{2}{\mu_1}$.

\subsubsection{Homodirectional systems of $m$ equations and $n+m$ heterodirectional systems} \label{Sec_4_1_2}

Once the problem of stabilizing two coupled scalar hyperbolic equations had been solved, subsequent studies have concentrated on non-scalar \textit{heterodirectional} hyperbolic systems with $n>0$ rightward propagating states and $m>0$ leftward propagating states. Such systems are referred to as $n+m$ heterodirectional systems. They correspond to equation~\eqref{Sec_4_1_v_pde_orig_simp} without ODEs. %Due to the in-domain and boundary coupling terms, the open-loop system may be unstable.
Note that $0+m$ or $n+0$ heterodirectional systems are qualified as \textit{homodirectional} systems since all transport velocities have the same signs, meaning all PDEs convect in the same direction. Due to the finite length of the spatial domain,  homodirectional systems are inherently stable, but the coupling between states can lead to undesirable transient behaviors, making trajectory planning problems non-trivial.

For $n+m$ heterodirectional systems with \textbf{right-boundary actuation} (i.e. we have $v(t,1)=Ru(t,1)+U(t)$), it has been shown in~\cite{li2010strong} that it is always possible to design a stabilizing feedback controller that stirs the state to zero in a minimum finite-time $t_F$ defined as $t_F=\frac{1}{\lambda_1}+\frac{1}{\mu_1}$ corresponding to the sum of the time needed for the slowest characteristic in each direction to travel the entire spatial domain. Analogously, it is also possible to design a state-observer that reconstructs the system's state in the same minimum finite-time $t_F$.  This minimum time is \emph{critical} in the sense that it is, in general, impossible to satisfy the given zero-terminal conditions if less time is allowed.  %To solve these two minimum-time problems (stabilization and observer design) %Several contributions focused on the design of backstepping controllers (resp. observers) that not only stabilize the system (resp. reconstruct the state) but do it in finite time. 
The backstepping techniques used in~\cite{Coron2013} to solve these two minimum-time stabilization and observer design problems when $n=m=1$ are extended in~\cite{DiMeglio2013a} to design appropriate finite-time output feedback controllers for $n+1$ systems using uncollocated measurements. 
%design appropriate finite-time output feedback controllers
The generalization of this result to an arbitrary number of $m$ leftward propagating equations is first presented in~\cite{Hu2016a}, where a backstepping transformation~\eqref{Sec_4_BS_transf_matrix_form} (without the ODE states and where $\mathbb{V}$ is a second kind Volterra integral operator) is used to map the original system to the following target system
\begin{equation}\label{Sec_4_1_2_target}
\left\{
\begin{array}{l}
\alpha_t+\Lambda^+\alpha_x=\Sigma^{++}\alpha+\Sigma^{+-}\beta\\
\hspace{1cm}+\int_0^xC^+(x,y)\alpha(y)+C^-(x,y)\beta(y)dy,\\
    \beta_t-\Lambda^-\beta_x=G(x)\beta(t,0),\\
    \alpha(0)=Q\beta(0),~\beta(1)=0,
\end{array} 
\right.
\end{equation}
where $C^+$ and $C^-$ are piecewise continuous functions, while the piecewise continuous matrix $G$ is upper triangular.  It can be noticed that contrary to~\cite{Coron2013}, the target system is not composed of transport equations, as, for instance, local cascade coupling terms are involved in the $\beta$-PDEs. Interestingly, for non-scalar $n+m$ hyperbolic systems, it has not been possible so far to derive a backstepping transformation able to get rid of all in-domain coupling terms. There, due to the triangularity of the matrix $G$ in~\eqref{Sec_4_1_2_target}, the control law~\eqref{Sec_4_1_1_control} yields finite-time convergence~\cite{Hu2016a}, but the convergence time $\bar t_F=\frac{1}{\lambda_1}+\sum_{i=1}^m\frac{1}{\mu_i}$ is larger than the theoretical minimum control time $t_F$.

To overcome this limitation, a different target system is proposed in~\cite{Auriol2016} where the non-local coupling terms of the $\beta$-PDE in~\eqref{Sec_4_1_2_target} ($G(x)\beta(t,0$) are replaced by local coupling terms ($\Omega(x)\beta(t,x)$ with $\Omega$ upper-triangular), thus enabling stabilization of the system in the minimum finite time $t_F$.  Both in~\cite{Hu2016a} and \cite{Auriol2016}, finite-time collocated and uncollocated state-observers are designed using a dual backstepping approach. This yields finite-time stabilizing output-feedback controllers.

An alternative target system (and therefore stabilizing control law) is proposed in~\cite{Coron2017}. It is given by
\begin{equation}\label{Sec_4_1_2_target_bis}
\left\{
\begin{array}{l}
\alpha_t+\Lambda^+\alpha_x=G_1(x)\beta(t,0),\\
    \beta_t-\Lambda^-\beta_x=G_2(x)\beta(t,1),\\
    \alpha(0)=Q\beta(0),~\beta(1)=0,
\end{array} 
\right.
\end{equation}
where $G_1$ and $G_2$ are piecewise continuous functions, $G_2$ being upper-triangular. Since $\beta(t,1)=0$, the term $G_2(x)\beta(t,1)$ vanishes, and the target system~\eqref{Sec_4_1_2_target_bis} converges to zero in the desired finite time $t_F$. The corresponding backstepping transformation can still be expressed as~\eqref{Sec_4_BS_transf_matrix_form}, but the operator $\mathbb{V}$ reads as the composition of a Volterra transformation (adjusted from~\cite{Hu2016a}) with an affine triangular integral transformation.
% Although the corresponding controller can still be rewritten as in~\eqref{Sec_4_1_1_control}, the system~\eqref{Sec_4_1_v_pde_orig_simp} is mapped to the target system~\eqref{Sec_4_1_2_target_bis} using two successive backstepping transformations: a backstepping Volterra transformation adjusted from~\cite{Hu2016a} and an affine integral transformation expressed as 
% \begin{align}
%     \beta(t,x)=z(t,x)+\int_0^1 K_z(x,y)z(t,y)dy,
% \end{align}
% where the kernel $K_z$ is upper triangular, which directly implies the invertibility of such a transformation.
Compared to~\cite{Auriol2016}, the target system~\eqref{Sec_4_1_2_target_bis} is simpler and more amenable for analysis since it has fewer coupling terms. Interestingly, if the matrix $Q$ presents specific structural properties, it is possible to modify the control input to stabilize the system in a finite time smaller than $t_F$~\cite{coron2019optimal}.

An important by-product of the analysis carried on in~\cite{Hu2016a} is that the target system~\eqref{Sec_4_1_2_target} (or equivalently the target system~\eqref{Sec_4_1_2_target_bis}) can be used to solve \emph{trajectory tracking} problems: given a known function $\Phi(t)$, find the value of $U(t)$ so that $v(t,0)=\Phi(t)$ for $t\geq t_F$. Indeed, applying the method of the characteristics to the $\beta$-equation in~\eqref{Sec_4_1_2_target}, it is possible to add a tracking component to the control input to solve the desired tracking problem. However, such a control law requires the knowledge of future values of the function $\Phi(t)$. 

%%%%%%%%%%%%%%%%%%%%%%%%%%%%%%%
\subsubsection{The case of isotachic blocks or zero transport speeds} \label{Sec_isotachic}
The case of ``isotachic'' blocks, which are blocks composed of states that share the same transport speed, was first mentioned in a remark in \cite{Hu2016a}. The authors pointed out that in this scenario, the direct application of the backstepping method would lead to singularities in the kernel equations. They suggested a diagonalizing transformation to decouple the isotachic states before applying the backstepping transformation, an idea that was further developed and explored in \cite{chen2023block} for multilayer Timoshenko beams (where multiple layers can potentially share the same physical properties, such as in laminated beams).

Consider a system of $m+n$ hyperbolic PDEs where a subset of $n_i$ states share the same transport speed $\sigma^-_i$, forming an ``isotachic block''. The dynamics of this block can be represented as:
\begin{equation}
    Y_{i,t} = \sigma^-_i Y_{i,x} + \Lambda_{i,i}Y_i + \cdots,
\end{equation}
where $Y_i\in \mathbb{R}^{n_i}$ is the vector of isotachic states and $\Lambda_{i,i}$ represents the coupling between them. The dots represent coupling with other states, which are not relevant. To remove this coupling, a block diagonalizing transformation is applied:
\begin{equation}
    \bar{Y}_i(x,t) = A_i(x) Y_i(t,x),
\end{equation}
where $A_i(x) \in \mathbb{R}^{n_i \times n_i}$ is a diagonal invertible matrix satisfying the following differential equation:
\begin{equation}
    \frac{dA_i(x)}{dx} = \frac{1}{\sigma^-_i} A_i(x) \Lambda_{i,i}, \quad A_i(0) = I_{n_i}.
\end{equation}
This transformation decouples the isotachic states, allowing for the application of the backstepping method without singularities in the kernel equations.

On the other hand, the control of hyperbolic PDEs with zero transport speeds (``atachic'' states) is a challenging problem that has received limited attention in the literature. The presence of atachic states can lead to non-stabilizability of the system, requiring additional assumptions on the system dynamics. In \cite{de2024backstepping}, the authors tackle the stabilization problem for a linear hyperbolic system with two counter-convecting PDEs and multiple atachic states. They provide a verifiable condition for the system to be stabilizable, requiring the atachic subsystem to be asymptotically stable in the absence of coupling with the nonzero-speed states. The control design employs an invertible Volterra transformation for the PDEs with nonzero speeds, leaving the atachic equations unaltered in the target system. The target atachic subsystem is then shown to be input-to-state stable (ISS) with respect to the decoupled and stable counter-convecting nonzero-speed equations. This approach allows for the design of a full-state backstepping controller that exponentially stabilizes the system origin in the $L^2$ sense.
%%%%%%%%%%%%%%%%%%%%%%%
 \subsection{Hyperbolic systems and Integral Delay Equations}  

 %(i.e., $n+m$ heterodirectional systems modeled by~\eqref{Sec_4_1_v_pde_orig_simp}) 
Hyperbolic systems of linear balance laws exhibit notable connections with a specific class of time-delay systems. Initially established by D'Alembert's formula~\cite{DAlembert}, which transforms wave equations into difference equations, this connection extends through the method of characteristics, allowing to represent as difference equations a class of linear conservation laws. %(i.e., equations modeled by~\eqref{Sec_4_1_v_pde_orig_simp} for which the coupling term $\Sigma$ equals zero)
In~\cite{russell1991neutral}, the existence of a mapping between the solutions of first-order hyperbolic PDEs and zero-order neutral systems is proved using spectral methods. This mapping is shown to be unique and is explicitly constructed in~\cite{karafyllis2014relation} for a single hyperbolic equation with a reaction term. Interestingly, the proposed methodology emphasizes the simplicity of stability analysis, as techniques initially developed for time-delay systems~\cite{halebook,michiels2009strong} can then be applied for the analysis of hyperbolic PDEs. We show below how the backstepping approach has been used to generalize the results of~\cite{karafyllis2014relation} to $n+m$ heterodirectional systems. Such connections between hyperbolic systems and time-delay systems can be crucial to designing stabilizing backstepping controllers for complex interconnected systems or to analyzing the robustness properties of the backstepping controllers.

\subsubsection{Time-delay formulation of $n+m$ heterodirectional systems} \label{Sec_time_delay_representation}

The backstepping approach is applied in~\cite{Auriol2019a} to rewrite $n+m$ heterodirectional systems modeled by equation~\eqref{Sec_4_1_v_pde_orig_simp} as \emph{Integral Delay Equations} (IDEs), i.e. difference equations with distributed-delay term. More precisely, backstepping transformations from~\cite{Coron2017} are applied in~\cite{Auriol2019a} without imposing any specific design on the control input so that the last boundary condition in the target system~\eqref{Sec_4_1_2_target_bis} is replaced by $\beta(1)=R\alpha(1)+\int_0^1 N(y)\gamma(t,y)dy+U(t)$ (with $\gamma=(\alpha^\top,\beta^\top)^\top$ and $N$ piecewise continuous). %The control law $U(t)$ is designed in~\cite{Coron2017} to obtain $z(t,1)=0$ and guarantee the finite-time stability of the system. 
Applying the method of characteristics, it is shown in~\cite{Auriol2019a} that the function $Z(t)=\beta(t,1)$ satisfies for all $t\geq t_F$ the following IDE
\small
\begin{align}
    Z(t)=&\sum_{i=1}^{n_0} F_i Z(t-\tau_i)+\int_0^{t_F} H(\nu)Z(t-\nu)d\nu+U, \label{Sec_4_2_IDE}
\end{align}
\normalsize
where $n_0$ is a positive integer (that depends on $n$ and $m$), the delays $0<\tau_1\leq \cdots \leq \tau_{n_0}=t_F$ depend on the transport velocities $\lambda_i$ and $\mu_i$, the constant matrices $F_i$ depend on the boundary coupling term $Q$ and $R$, while the function $H$ is piecewise continuous. The state $Z(t)$ can be seen as a function that belongs to $L^2_{t_F}=L^2([-t_F,0],\mathbb{R}^m)$, i.e.  the Banach space of $L^2$ functions mapping the interval $[-t_F, 0]$ into $\mathbb{R}^m$. The associated norm is defined as $||Z||_{L^2_{t_F}}=(\int_{-t_F}^0 ||Z(t+s)||^2 dt)^{\frac{1}{2}}.$ We refer the reader to~\cite{halebook} for more details on this mathematical framework.
It is shown in~\cite{Auriol2019a} and \cite[Theorem 6.1.3]{auriol2024contributions} that the IDE~\eqref{Sec_4_2_IDE} and the PDE system~\eqref{Sec_4_1_v_pde_orig_simp} have equivalent stability properties. More precisely, there exist $\tau^\star>0$ and two constants $\kappa_0$ and $\kappa_1$, such that for any linear bounded state-feedback law $U(t)$, for any $t>t_F$
    \begin{align}
        \kappa_0 ||Z||_{L^2_{\tau^\star}}^2 \leq ||w||^2_{L^2} \leq \kappa_1 ||Z||_{L^2_{t_F}}^2.
    \end{align}
    Moreover, the exponential stability of~$Z$ in the sense of the~$L^2_{\tau_{n_0}}$-norm is equivalent to the exponential stability of~$w$ in the sense of the spatial~$L^2$ norm.

%The fact that the norms are different on the two sides of the inequality is related to the structure of the difference equation (see, for instance, the design of converse Lyapunov-Krasovskii functions~\cite{pepe2013converse}).
Interestingly, the representation~\eqref{Sec_4_2_IDE} is more amenable for open-loop stability analysis or closed-loop robustness analysis~\cite{damak2015stability,niculescu2001delay}. %As it will be shown in this survey, such a representation was crucial to designing backstepping-based controllers for networks of hyperbolic systems. %For instance, the control law designed in~\cite{Coron2017} corresponds to $U(t)=-\sum_{i=1}^{n_0} F_i Z(t-\tau_i)-\int_0^{\tau_{n_0}} H(\nu)Z(t-\nu)d\nu\nonumber$. 
Taking the Laplace transform of~\eqref{Sec_4_2_IDE} (denoting $s$ the Laplace variable), we obtain the open-loop characteristic equation~\cite{halebook,Auriol2019a}
\begin{align}
    \det\big(\text{Id}-\sum_{i=1}^{n_0} F_i \mathrm{e}^{-\tau_i s}+\int_0^{t_F} H(\nu)\mathrm{e}^{-\nu s}d\nu\big)=0. \label{eq_Laplace_Sec_4}
\end{align}
As shown in~\cite{henry1974linear} the exponential open-loop stability of~\eqref{Sec_4_2_IDE} is equivalent to having a strictly negative spectral abscissa, i.e., all the roots of~\eqref{eq_Laplace_Sec_4} must have a real part smaller than $-\eta<0$. Therefore, the stability analysis of~\eqref{Sec_4_1_v_pde_orig_simp} boils down to the problem of roots location for equation~\eqref{eq_Laplace_Sec_4}. In the case of scalar equations ($n=m=1$) with constant coefficients, it is possible to obtain an explicit expression of the function $H$ and to introduce simple sufficient conditions to guarantee open-loop stability~\cite{saba2019stability}. %Interestingly, the time-delay representation~\eqref{Sec_4_2_IDE} is of specific interest for the robustness analysis of the backstepping controller~\eqref{Sec_4_1_1_control}.

\subsubsection{Delay-robustness of backstepping controllers} \label{Sec_Delay_rob}

%As it has been shown in Section~\ref{Sec_hyp_system}, 
Most contributions in the literature initially focused on designing backstepping controllers stabilizing hyperbolic systems in finite-time, thereby shadowing the robustness properties of the corresponding closed-loop systems~\cite{logemann1996conditions,michiels2009strong}. These robustness constraints may stem from uncertainties in the parameters, disturbances acting on the system, noise on the measurements, neglected dynamics, or delays acting on the actuators. It has been observed (see~\cite{datko1986example,logemann1996conditions}) that for many feedback systems, the introduction of arbitrarily small time delays in the control input may cause instability for any feedback. In particular, for $n=m=1$, it is shown in~\cite{Auriol2018} that the controller~\eqref{Sec_4_1_1_control} designed to guarantee finite-time stability may have zero delay-robustness margins. This result is generalized in~\cite{Auriol2019a} to $n+m$ heterodirectional systems and is directly obtained from the time-delay representation~\eqref{Sec_4_2_IDE}. 

 %It has been observed (see~\cite{datko1986example,logemann1996conditions}) that 
 If the open-loop transfer function of a feedback system has an infinite number of poles in the complex right-half plane, then the introduction of arbitrarily small time delays in the loop may cause instability for any feedback control law~\cite{datko1986example,logemann1996conditions}. Therefore, it is first essential to analyze the open-loop behavior of the system~\eqref{Sec_hyp_system}. Based on the characteristic equation~\eqref{eq_Laplace_Sec_4}, it is required in~\cite{Auriol2019a,Auriol2023} that the open-loop PDE system~\eqref{Sec_4_1_v_pde_orig_simp} with $\Sigma(x)=0$ and in the absence of the ODE states is exponentially stable. %This avoids having an infinite number of unstable poles in open loop.
 %since the integral term vanishes at high frequencies, it is shown in~\cite{Auriol2019a,halebook} that the principal part of the IDE~\eqref{Sec_4_2_IDE} (i.e. the equation $Z(t)=\sum_{i=1}^{n_0}F_iZ(t-\tau_i)$) must be exponentially stable to avoid an infinite number of unstable poles. This condition implies the following requirement for the PDE system~\eqref{Sec_4_1_v_pde_orig_simp}
 % \begin{assum} \cite[Theorem 4]{Auriol2019a} \label{Assum_robust}
 %     The open-loop PDE system~\eqref{Sec_4_1_v_pde_orig_simp} with $\Sigma(x)=0$ and in the absence of the ODE states is exponentially stable. 
 % \end{assum}
This assumption is a \textbf{necessary condition} for delay-robust stabilization. It corresponds to the exponential stability of the principal part of the IDE~\eqref{Sec_4_2_IDE} (i.e., the equation $Z(t)=\sum_{i=1}^{n_0}F_iZ(t-\tau_i)$)~\cite{halebook} and avoids having an infinite number of unstable poles in open loop. %Necessary and sufficient conditions to characterize the stability of the principal part of ~\eqref{Sec_4_2_IDE} can be found in~\cite{halebook}, although they may not be tractable for systems of high dimensions. Therefore, sufficient conditions based on Lyapunov-Krasovskii methods~\cite{pepe2005asymptotic,fridman2002stability} may be more applicable.
For $n=m=1$, this assumption is equivalent to $|\rho q|<1$. For larger values of $n$ and $m$, it can be verified using sufficient conditions based on Lyapunov-Krasovskii methods~\cite{pepe2005asymptotic,fridman2002stability}. 
%\modifJA{Please checked the paragraph I commented below (in the .tex) to confirm its suppression.}
%One should be aware that this assumption does not hold for most practical applications. Models of the form~\eqref{Sec_4_1_v_pde_orig_simp} are simplistic and do not capture phenomena that would be susceptible to making the delay margins non-null. For instance, they may neglect the diffusivity stemming from Kelvin-Voigt damping. Therefore, although the delay-robustness margins for such systems would be poor, it remains theoretically possible to stabilize them robustly.

Even when this assumption is verified, it is shown in \cite{Auriol2019a} that the backstepping controllers may not be robust to input delays as they significantly trade-off delay-robustness for performance to guarantee finite-time stability. More precisely, the controllers with a structure as in~\eqref{Sec_4_1_1_control} are composed of two parts:
\begin{enumerate}
	\item the integral part whose objective is to remove the effect of in-domain couplings;
	\item the term~$-R u(t,1)$ whose objective is to cancel the proximal reflection and to ensure finite-time stability.
\end{enumerate}
It has been observed in~\cite{kyllingstad2009new} that canceling the proximal reflection coefficient can change the dynamics in a way that makes the system unstable. This robustness issue is underlined in~\cite{Auriol2019a}: due to the compensation of the term~$-R u(t,1)$ the backstepping feedback control operator may not be strictly proper~\cite{curtain2012introduction}, thus inducing vanishing robustness margins. 
%These observations are consistent with reports by industrial practitioners on the limitations of the \emph{impedance matching method}. The impedance matching method (see~\cite{egeland2002modeling}) consists of matching the load impedance to the characteristic line impedance. 
%The backstepping finite-time controllers (as~\eqref{Sec_4_1_1_control}) can be expressed as $U(t)=-\sum_{i=1}^{n_0} F_i Z(t-\tau_i)-\int_0^{\tau_{n_0}} H(\nu)Z(t-\nu)d\nu$ \cite{Auriol2019a}. In the presence of an arbitrarily small delay, it can be verified that due to the term $-\sum_{i=1}^{n_0} F_i Z(t-\tau_i)$, the characteristic equation of the closed-loop system may have unstable poles. 
%For instance, in~\cite{kyllingstad2009new}, the authors designed a controller preventing stick-slip oscillations, a class of undesired torsional oscillations characterized by a series of \emph{stick} (a cessation of bit rotation) and \emph{slip} (a sudden release of rotational energy) that often occur in drilling devices. They 
 Two solutions are proposed in the literature to modify the control law~\eqref{Sec_4_1_1_control} and guarantee robustness margins, thereby giving up finite time stabilization
\begin{enumerate}
    \item Cancel only a part of the reflection term $-Ru(t,1)$. This approach is applied in~\cite{Auriol2018c,Auriol2020e} and results in the introduction of tuning parameters in the design, thus guaranteeing potential trade-offs between different specifications (namely delay-robustness and convergence rate). %The gained robustness comes at the price of degraded nominal performance.
    \item Combine the non-strictly proper controller with a well-tuned low-pass filter~\cite{Auriol2023} The resulting control law then becomes strictly proper, which guarantees the existence of robustness margins. This filtering approach simplifies the design of stabilizing controllers by separating the stabilization and robustness concerns. %Specifically, when integrated with the controller, the low-pass filter guarantees delay-robust stability of the closed-loop dynamics. 
    Designing these filters involves leveraging the insight that robustness challenges arise at high frequencies. 
\end{enumerate}

%Both approaches have been used in the literature to tackle these robustness issues when designing robust controllers.

%Dire qu'il n'y a pas d'analyse quantitative pour le moment ?
%%%%%%%%%%%%%%%%%%%%%%%%%%%%%%%%%%%%%%%%%%%
\subsection{Interconnected hyperbolic systems and ODEs} 

The control of systems of hyperbolic PDEs coupled with ODEs is a very active research topic as this class of system naturally emerges when the distributed dynamics is subject to actuator and/or load dynamics. Among possible applications, we can cite the attenuation of mechanical vibrations in drilling devices (where the hyperbolic PDEs represent axial and torsional stress propagation, while the ODEs model the actuation and the Bottom Hole Assembly (BHA) dynamics)~\cite{Saldivar2016}, deepwater construction vessels~\cite{Tang2011}, or the Rijke tube~\cite{de2020backstepping}. 
Most constructive control designs for such systems are based on the backstepping approach. Interestingly, ODE-PDE interconnections also arise when considering dynamical delays with input or sensor delays. Therefore, the backstepping approach has also proved its usefulness in solving advanced problems for time-delay systems.

\subsubsection{Advanced problems for time-delay systems}
%\modifJA{We can remove one column if we do not give the equations.}
%ODE systems subject to input delays can easily be modeled as cascaded PDE-ODE systems, the PDE being a transport equation whose velocity matrix is related to the delays. 
When controlling dynamical systems, actuation or transmission delays are unavoidable. These delays often lead to detrimental effects on performance and stability, which explains why numerous active control approaches have been designed to compensate for them~\cite{smith1959controller,artstein1982linear,Manitius1979}. In this context, predictor-based control entails incorporating predicted future states into the control design to counteract the destabilizing effects of time delays (see~\cite{deng2022predictor} for a survey on predictor-based controllers for time-delay systems). Advancements in predictor-based control techniques have been greatly enhanced by reformulating delays using additional transport equations, as highlighted in \cite{krstic5}. This approach enables compensation for potentially long-time delays using backstepping controllers. More precisely, on the seminal paper~\cite{Krstic2008c}, a re-interpretation of the classical Finite Spectrum Assignment \cite{Manitius1979} is proposed, modeling ODEs with input delays as PDE-ODE interconnections. Consider the following linear ODE with a delayed control input
\begin{align}
    \dot{X}(t)=AX(t)+BU(t-\tau), \label{Sec_4_Delayed_ODE}
\end{align}
where $\tau>0$. Defining $v(t,x)=U(t-\tau(1-x)) \in \mathbb{R},~x \in [0,1]$, this time-delay can be rewritten as the following PDE-ODE system
\begin{equation}\label{Sec_4_ODE_PDE}
\left\{
\begin{array}{l}
\dot{X}(t)=AX(t)+Bv(t,0),\\
    v_t-\frac{1}{\tau}v_x=0,~v(1)=U(t).
\end{array} 
\right.
\end{equation}
Assume that the pair $(A,B)$ is stabilizable (that is there exists a matrix $K_0$ such that $\bar A=A+BK_0$ is Hurwitz).
Using the backstepping approach, the system~\eqref{Sec_4_ODE_PDE} is mapped in~\cite{Krstic2008c} to an exponentially stable target system with the same structure, the matrix $A$ being replaced by $\bar A$. The backstepping transformation reads as~\eqref{Sec_4_BS_transf_matrix_form} (with $\mathbb{I}_1=F_X=0$ and without the $Y$ state). In this simple configuration, it is possible to obtain an explicit expression of the backstepping kernels and express the stabilizing control law as
\begin{align}
    U(t)=K_0(\mathrm{e}^{\tau A}X(t)+\tau\int_0^1\mathrm{e}^{\tau A(1-y)}Bv(t,y)dy),
\end{align}
which is equivalent to the ODE-controller $U(t)=K_0X(t+\tau)$ computed by finite spectrum assignment methods. All in all, this application of the backstepping approach for time-delay systems can be summarized as follows
\begin{enumerate}
    \item Rewrite the time-delay system as a PDE-ODE system;
    \item Map this PDE-ODE system into an appropriate target system using a backstepping transformation;
    \item Compute the corresponding control law.
\end{enumerate}
Such a backstepping predictor-based strategy is successfully applied in~\cite{krstic2010lyapunov} in the case of time-varying input delays. Moreover, an observer reconstructing the state $X(t)$ using delayed measurement is also obtained therein using a dual approach. The case of distributed input delays is solved in~\cite{bekiaris2011lyapunov}, and an extension is proposed in \cite{zhu2020predictor} in the presence of model uncertainties. This backstepping predictor-based strategy is adjusted in~\cite{holloway2019prescribed,espitia2021predictor} (using time-varying kernels) to stabilize~\eqref{Sec_4_1_v_pde_orig_simp} in prescribed time.

%Although initially designed for linear infinite-dimensional systems, 
The backstepping approach could then be extended to non-linear ODEs with input delays. The case of multi-input nonlinear systems with constant delays is solved in~\cite{bekiaris2016predictor}, while the case of time-varying input delays is considered in~\cite{bekiaris2012compensation} for the system
\begin{align}
    \dot{X}(t)=f(X(t),U(t-\tau(t)), \label{Sec_4_Delayed_ODE_NL}
\end{align}
where $f$ is $C^1$ and the time-varying delay $\tau$ is bounded and verifies $|\dot{\tau}|<1$. The additional PDE state can be defined as $v(t,w)=U(\phi(t+x(\phi^{-1}(t)-t)))$, where $\phi(t)=t-\tau(t)$. Assuming the the delay-free ODE is stabilized by $U(t)=\kappa(t,X(t))$  and under an appropriate ISS condition given in~\cite{bekiaris2012compensation}, it is possible to map the system~\eqref{Sec_4_Delayed_ODE_NL} to the asymptotically stable target system
\begin{equation}\label{Sec_4_ODE_PDE_target_NL}
\left\{
\begin{array}{l}
\dot{X}(t)=f(X(t),\kappa(t,X(t))+\beta(t,0)),\\
    \beta_t-\frac{1+x(\frac{d(\phi^{-1}(t))}{dt}-1)}{\phi^{-1}(t)-t}\beta_x=0,~\beta(1)=0.
\end{array} 
\right.
\end{equation} 
The corresponding backstepping transformation reads
\begin{align}
    \beta(t,x)=v(t,x)-\kappa(t+x(\phi^{-1}(t)-t),p(t,x)),
\end{align}
where $p(t,x)=X(t)+(\phi^{-1}(t)-t)\int_0^xf(p(t,y),v(t,y))dy$. In the case where the inverse function $\phi^{-1}$ is not available, approximated prediction horizons are used in~\cite{bresch2018robust} using $\phi^{-1}(t) \approx t+\tau(t)$. Such a design is adjusted in~\cite{zekraoui2023finite} to stabilize a chain of integrators with input delays in finite/fixed time.

%Analogously to what has been presented in Section~\ref{Sec_Delay_rob} for $n+m$ heterodirectional hyperbolic systems, 
The backstepping methodology has proven valuable in analyzing the delay-robustness properties of predictor-based controllers. For linear ODEs, conditions on the constant delay mismatch that guarantee robustness are given in~\cite{krstic2008lyapunov}. Recently, these results have been extended in~\cite{kong2022prediction} to obtain conditions guaranteeing mean-square exponential stability in the presence of a stochastic delay. For nonlinear systems, the case of time-dependent and state-dependent delay uncertainties is tackled in~\cite{bekiaris2013robustness}, and closed-loop robustness is shown if the magnitudes of the delay uncertainty and its partial derivatives are small enough. %\modifJA{Maybe what follows can be moved to the adaptive section?}
 For unknown (but constant) time delays, \textit{adaptive backstepping controllers} have been introduced in~\cite{krstic2009delay} when an upper bound of the delay (denoted $\bar \tau$) is known. 
% \modifJA{
% Consider equation~\eqref{Sec_4_Delayed_ODE} and its PDE-ODE representation \eqref{Sec_4_ODE_PDE} (where $\tau$ is from now on unknown). Denote $\hat \tau(t)$ as the online estimation of the delay. Consider the invertible backstepping transformation (adjusted from~\eqref{Sec_4_transf_BS_ODE})
% \begin{align*}
%  \beta(x)= v(x)-K_0[\int_0^x\hat \tau \mathrm{e}^{\hat \tau A(x-y)}Bv(y)dy+\mathrm{e}^{\hat \tau Ax}X],
% \end{align*}
% We obtain the target system
% \begin{equation}\label{Sec_4_ODE_PDE_target_adaptive}
% \left\{
% \begin{array}{l}
% \dot{X}(t)=(A+BK_0)X(t)+B\beta(t,0),\\
%     \tau(t) \beta_t-\beta_x=-\tilde \tau(t)p(t,x)-\tau \dot{\hat \tau}(t)q(t,x),\\
%     \beta(1)=0,
% \end{array} 
% \right.
% \end{equation}
% where $\tilde \tau(t)=\tau-\hat \tau(t)$, $p(t,x)=K\mathrm{e}^{\hat \tau(t)Ax}[(A+BK_0)X(t)+B\beta(t,0)]$, and $q(t,x)=\int_0^x K_0(\text{Id}+A\hat \tau(t)(x-y))\mathrm{e}^{\hat \tau(t)A(x-y)}Bv(t,y)dy+K\mathrm{e}^{\hat \tau(t)Ax}X(t).$ The adaptation law of $\hat \tau$ is designed in~\cite{bresch2010delay} as
% \begin{align*}
%     \dot{\hat \tau}(t)&=-\eta\text{Proj}_{[\underline \tau,\bar \tau]}\{g(t)\}, ~ \eta>0,  \label{eq_proj_Sec_4}\\
%     g(t)&=\frac{\int_0^1(1+x)\beta(t,x)p(t,x)dx}{1+X(t)^TPX(t)+b\int_0^1(1+x)\beta(t,x)^2dx},
% \end{align*}
% where $\text{Proj}$ refers to the standard projection operator~\cite{krstic2009delay}, and $P$ (resp. $b$) is an appropriate constant matrix (resp. scalar). This adaptation law is chosen to guarantee the negativity of an appropriate Lyapunov-Krasovskii functional.}
Such adaptive control law require knowing the distributed state $u(t,x)$, which is restrictive. An input observer is designed in~\cite{bresch2010delay} to estimate $u(t,x)$. However, global stability is no longer guaranteed in this configuration. This control framework has been extended to systems with distinct delays~\cite{zhu2018pde} or distributed delays~\cite{zhu2020predictor}. Output-feedback controller designs are addressed in~\cite{bresch2012adaptive} using dual approaches while solving input disturbance rejection. Finally, we mention that a neural network-based estimation scheme is proposed in~\cite{chakraborty2017control} to estimate the unknown input delay magnitude.

%%%%%%%%%%%%%%%%%%%%%%%%%%%%%%%%%%%%%%
\subsubsection{Interconnected ODE-PDE systems} \label{Sec_ODE_PDE}

As we have seen above, it has been possible to design backstepping stabilizing controllers for dynamical systems with input and measurement delays. Subsequently, these control techniques have been extended to design observers, controllers, or parameter estimation methods for a general class of interconnected ODE-PDE systems, for which the PDE subsystem is not a simple transport equation. Consider the system~\eqref{Sec_4_1_v_pde_orig_simp} in the absence of the $Y$-ODE state (i.e. $A_1=E_1=C_1=0$). Although different actuator configurations can be considered (as either the ODE state or the PDE state can be actuated), most contributions in the literature focus on the case of an actuated PDE. Therefore, we have $v(1)=Ru(1)+U(t)$. To deal with the joined presence of the ODE and the heterodirectional $n+m$ PDE system, the backstepping transformations found in the literature have the structure~\eqref{Sec_4_BS_transf_matrix_form} (without the $Y$ term and the corresponding terms in the transformation). In~\cite{DiMeglio2018}, such a backstepping transformation is used to map the system~\eqref{Sec_4_1_v_pde_orig_simp} to the exponentially stable target system
\begin{equation}\label{Sec_4_3_ODE-PDE_target}
\left\{
\begin{array}{l}
\dot{X}=(A_0+E_0K_0) X+E_0 \beta(t,0), \\
    \alpha_t+\Lambda^+\alpha_x=\Sigma^{++}\alpha+\Sigma^{+-}\beta+D(x)X(t)\\
\hspace{1cm}+\int_0^xC^+(x,y)\alpha(y)+C^-(x,y)\beta(y)dy,\\
    \beta_t-\Lambda^-\beta_x=G(x)\beta(t,0),\\
    \alpha(0)=Q\beta(0)+C_0X(t),~\beta(1)=0,
\end{array} 
\right.
\end{equation}
where $G$ is strictly upper-triangular. Provided that $(A_0+E_0K_0)$ is Hurwitz, the target system~\eqref{Sec_4_3_ODE-PDE_target} is exponentially stable. This target system can be seen as the natural extension of~\eqref{Sec_4_1_2_target}. Due to the presence of the ODE term in the backstepping transformation, the kernel equations can be rewritten as a coupled PDE-ODE system, whose well-posedness is shown in~\cite[Theorem 3.2]{DiMeglio2018}. Compared to~\eqref{Sec_4_1_1_control}, the stabilizing control input includes an additional term of the form $FX(t)$.
% \begin{align*} \label{Sec_ODE_PDE_control}
%     U(t)=-Ru(t,1)+FX(t)+\int_0^1 \begin{pmatrix}
%         0 \\ Id
%     \end{pmatrix}^TK(1,y)w(y)dy.
% \end{align*}
%However, reinterpreting the ODEs as PDEs with horizontal characteristic lines, the well-posedness proof proposed in~\cite{Hu2016a} can be directly adjusted. Interestingly, the well-posedness of a generic class of hyperbolic PDEs on a triangular domain is given in~\cite{DiMeglio2018}. 

The control law proposed in~\cite{DiMeglio2018} requires canceling all the reflection term $Ru(t,1)$. As it has been emphasized in Section~\ref{Sec_Delay_rob}, this may lead to vanishing robustness margins. Therefore, in the case of a $2 \times 2$ PDE system ($n=m=1$), a new design is provided in~\cite{Auriol2018a} to ensure delay-robust stabilization by preserving some proximal reflection terms in the control law. The control design is based on a rewriting of the system as a time-delay neutral system. This methodology is extended in~\cite{Auriol2019d} in the case of a $n+m$ heterodirectional system. 

% The control input now has two components $U_{PDE}$ (used to stabilize the PDE by killing the in-domain coupling terms) and $U_{ODE}$ (used to stabilize the ODE part). The component $U_{PDE}$ has the structure
% \begin{align*} 
%     U_{PDE}(t)=-\tilde Ru(1)+\tilde FX+\int_0^1 \begin{pmatrix}
%         0 \\ Id
%     \end{pmatrix}^T\tilde K(1,y)w(y)dy,
% \end{align*}
% where $\tilde R$ is such that $|\tilde R Q|+|RQ|<1$. The target system becomes 
% \begin{equation}\label{Sec_4_3_ODE-PDE_target_simple}
% \left\{
% \begin{array}{l}
% \dot{X}=(A_0+E_0K_0) X+E_0 \beta(t,0), \\
%     \alpha_t+\Lambda^+\alpha_x=0,~\beta_t-\Lambda^-\beta_x=0,\\
%     \alpha(0)=Q\beta(0)+C_0X(t),\\
%     \beta(1)=(R-\tilde R)\alpha(1)+U_{ODE}.
% \end{array} 
% \right.
% \end{equation}
% Using the method of characteristics, it is then shown in~\cite{Auriol2018a} that the state $X(t)$ verifies a time-delay neutral equation that can be stabilized using $U_{ODE}$, therefore implying the exponential stability of the interconnected system (due to the underlying cascade structure). Delay-robustness is ensured by preserving some proximal reflection terms in the control law. The methodology developed in~\cite{Auriol2018a} is extended in~\cite{Auriol2019d} in the case of a $n+m$ heterodirectional system. 

An alternative two-step method for the backstepping stabilization of coupled linear heterodirectional hyperbolic PDE-ODE systems is proposed in~\cite{deutscher2019output}
\begin{enumerate}
    \item The first step consists of using a classical Volterra backstepping transformation (that does not depend on the ODE state) to simplify the PDE subsystem;
    \item Then, it is possible to use a second backstepping transformation (that now depends on the ODE state) to decouple the target system into a PDE-ODE cascade.
\end{enumerate} 
The main advantage of this approach is that only the conventional kernel equations have to be solved, whereas the decoupling kernel equations take a simple form. This two-step approach has been adjusted in~\cite{deutscher2019output} to design a state-observer. Adjustments of this multi-step approach are proposed in~\cite{gehring2023control}  for scalar hyperbolic systems ($n=m=1)$ coupled with ODE systems in normal forms.  The proposed backstepping design exploits the strict feedback form of the system to recursively stabilize it. Extensions to the 
case of a non-linear ODE are proposed in~\cite{irscheid2023output}. As the proposed controller requires prediction of the ODE and PDE states, appropriate predictors are designed as the unique solution of general nonlinear Volterra integro-differential equations.

Other ODE-PDE configurations, slightly different from~\eqref{Sec_4_1_2_target}, can be found in the literature. In~\cite{mathiyalagan2022observer}, a backstepping output-feedback controller is designed for a homodirectional system ($n=0$) coupled at both boundaries with an ODE. Moreover, the boundary at $x=1$ is not completely actuated, as we have $v(t,1)=Ru(t,1)+\begin{pmatrix}
         C_1X(t) & U(t)
     \end{pmatrix}^\top.$
This type of under-actuated system can be frequently seen in the control of multiphase flows as it describes the simultaneous flow of materials that differ in their thermodynamic properties.  In~\cite{Hasan2016}, a state-observer is designed for a semilinear PDE-ODE cascade (i.e., $E_0=0$), using the measurement $u(t,1)$. 

 Systems with state-dependent or time-dependent input delays have been addressed, for instance, in \cite{diagne2017time}, \cite{diagne2017compensation}, and \cite{diagne2017control} developing predictor-feedback control designs with applications to 3D printing processes, where material transport in screw extruders creates complex input delay dynamics. 
Finally, we emphasize that specific attention has been recently paid to wave equations coupled with ODEs~\cite{roman2018backstepping,deutscher2021backstepping,Zhou2012}. 

% The backstepping-based control designs presented above have been successfully adapted to deal with systems having a configuration similar to that of~\eqref{Sec_4_3_ODE-PDE}. The multi-step approach derived in~\cite{deutscher2019output} is extended in~\cite{irscheid2023output} to the case of a non-linear ODE. As the proposed controller requires prediction of the ODE and PDE states, appropriate predictors are designed as the unique solution of general nonlinear Volterra integro-differential equations. %The same cascade configuration has been considered in~\cite{ghousein2020adaptive} with a homodirectional linear PDE system but with time-varying $A_0$ and $C_0$ matrices. Moreover, in this contribution, the control input can simultaneously act on the PDE and the ODE. case of time-varying ODE coefficients 

\subsubsection{Interconnected ODE-PDE-ODE systems} \label{Sec_ODE_PDE_ODE}

Many recent publications have recently focused on the control of ODE-PDE-ODE interconnections (where the ODEs correspond to actuator and load dynamics) described by equation~\eqref{Sec_4_1_2_target}. These contributions usually consider that the control input $U(t)$ acts on the ODE state $X(t)$ and that a part of the $Y$-ODE state $y(t)=C_{mes}Y(t)$ is measured. The case of collocated measurements can be solved using symmetry arguments. 

The first application of backstepping to the stabilization of interconnected ODE-PDE-ODE systems is found in~\cite{BouSaba2017ODE-PDE-ODE}, for two scalar heterodirectional transport equations ($n=m=1$) with potentially unstable ODE ``actuator" and ``load" dynamics, where the ``actuator" ODE is required to have an invertible input matrix, and the ``load" ODE is required to be stabilizable. This result obtains a target system where in-domain couplings in the PDEs are removed, as well as the proximal reflection coefficient on the PDEs (i.e., on the side of the actuated ODE) and the coupling from the PDEs to the actuated ODE. Furthermore, in the target system, the two ODEs are stabilized.

Another early contribution can be found in~\cite{anfinsen2018stabilization}, limited to two scalar ODEs, and the PDE system corresponds to a scalar transport equation from the $X$-ODE to the $Y$-ODE. In other words, this contribution shows how to incorporate first-order actuator and sensor dynamics into the controller and observer designs for the scalar 1-D linear hyperbolic PDE system considered in~\cite{Krstic2008c}. The proposed state-feedback controller uses the transformation~\eqref{Sec_4_1_2_target_bis} (without the ODE terms) and does not guarantee the exponential stability of the ODE state $Y$. The case of a scalar heterodirectional PDE system ($n=m=1$) is considered in~\cite{DiMeglio2020} for a scalar actuated ODE (i.e., $p=1$). Under appropriate stabilizability conditions for the ODE systems, the original interconnected system~\eqref{Sec_4_1_v_pde_orig_simp} is mapped to a target system with a similar structure but with $\Sigma=0$ and where $A_0$ and $A_1$ are replaced by Hurwitz matrices. The matrix $E_0$ is also replaced by a matrix $\bar E_0$ chosen to guarantee the closed-loop exponential stability. Rather than canceling the entire proximal wave reflection in the control law, it is assigned the dynamics of a high-pass filter. Therefore, the resulting target system dynamics take the form of a cascade of the distal ODE into a first-order neutral system, which, when stable, is shown to be $w$-stable in the sense~\cite{curtain2012introduction}. As explained in Section~\ref{Sec_Delay_rob}, such a filtering approach robustifies the controller. A state-observer relying on collocated measurements is also designed in~\cite{DiMeglio2020} using a dual approach.

To deal with non-scalar PDE and ODE subsystems, it is possible to extend the two-step approach presented in~\cite{deutscher2019output}. In the first step, the PDE subsystem is mapped into a target system with a simpler structure using a Volterra integral transform. This significantly facilitates the decoupling into a stable ODE-PDE–ODE cascade in the second step, as shown in~\cite{Deutscher2018}. However, obtaining the desired stable ODE–PDE–ODE cascade in the last step requires the introduction of new coordinates to represent the $Y$ ODE in its multivariable Byrnes–Isidori normal form~\cite{isidori2013nonlinear}. For this, a vector relative degree of one is assumed, which is a requirement often met in applications. \\
The limitations induced by this relative degree one condition are partially mitigated in \cite{BouSaba2019_b} for a scalar PDE system, where the $X$-ODE is presumed to be minimum phase for the influencing output of the PDE. This article was also the first to include a strictly proper control for ODE-PDE-ODE systems, by adding an adequately chosen low-pass filter in the control design. As mentioned in previous sections, this can have a large impact on the resulting robustness of the control laws for PDE systems. This methodology was broadened in \cite{Wang2020} to include a state observer, thereby allowing the design of an output-feedback controller. Nonetheless, the PDE subsystem is still assumed to be scalar in this contribution. The extension to the stabilization of an interconnected ODE-$n+m$ heterodirectional PDE-ODE system is proposed in~\cite{Auriol2023}, while a state observer can be found in~\cite{Auriol2022b}. Using the backstepping transformation~\eqref{Sec_4_BS_transf_matrix_form} (with $F_X=0$), the system~\eqref{Sec_4_1_v_pde_orig_simp} is mapped to a target system that can, in turn, be rewritten as an IDE system (coupled with ODEs) with equivalent stability properties \cite[Lemma 6.1.4]{auriol2024contributions}.
% Analogously to what has been already presented in Section~\ref{Sec_time_delay_representation}, the design in both cases is based on a time-delay representation of the ODE-PDE-ODE system~\eqref{Sec_4_3_ODE_PDE_ODE}. The methodology is analogous to the one presented in Section~\ref{Sec_time_delay_representation}: using the backstepping transformation~\eqref{Sec_4_BS_transf_matrix_form} (with $F_X=0$), the system~\eqref{Sec_4_3_ODE_PDE_ODE} is mapped to a target system that can, in turn, be rewritten as an IDE system (coupled with ODEs) with equivalent stability properties \cite[Lemma 6.1.4]{auriol2024contributions}. 
A stabilizing controller can then be designed in the frequency domain~\cite{Auriol2023}. The proposed approach requires the following assumptions.
\begin{enumerate}
    \item The open-loop system without the ODEs and the in-domain coupling terms is exponentially stable.
    \item The pairs $(A_1,E_1)$ and $(A_0,B)$ are stabilizable.
    \item For all complex number $s$ with a positive real part, we have rank $\begin{pmatrix}
        s\text{Id}-A_0 & B \\ C_0 & 0
    \end{pmatrix}=p+n$;
\end{enumerate}
The first condition corresponds to the (delay-) robustness assumption stated in~\cite[Theorem 4]{Auriol2019a} and is necessary to guarantee the delay-robustness of the control input (see Section~\ref{Sec_time_delay_representation}). The second assumption is not very conservative since it is impossible to stabilize the ODE subsystems, independently of the PDE or interconnection structure, without the stabilizability of $(A_1,E_1)$ and $(A_0,[B\ E_0])$ at the very least. The last assumption implies that the matrices $C_0$ and $B$ are not identically zero and that $P_0(s)=C_0(sI-\bar A_0)^{-1}B$ does not have any zeros in $\mathbb{C}^+$ that are common to all its components (and therefore admits a stable right inverse). Although conservative, this last assumption is less restrictive than other conditions existing in the literature and previously mentioned. The robustness of the closed-loop system is guaranteed in~\cite{Auriol2023} using an appropriate filtering procedure (explained in Section~\ref{Sec_time_delay_representation}). 

% that there exist a positive integer $n_0$, positive delays $0<\tau_1\cdots \leq \tau_{n_0}$, constant matrices $G_1$ and $G_2$, and piecewise continuous bounded functions $F_i^X, F_i^Y, F_i^z, H_X, H_Y$ and $H_z$ such that \eqref{Sec_4_3_ODE_PDE_ODE} rewrites 
% \begin{equation}\label{Sec_4_3_ODE-PDE_IDE_form}
% \left\{
% \begin{array}{l}
% z(t)=\sum_{i=1}^{n_0}F_i^zz(t-\tau_i)+G_1Y\\
% \hspace{1cm}+C_0\bar X+\int_0^{\tau_{n_0}}H_z(\nu)z(t-\nu)d\nu,\\
% \dot{\bar X}=A_0 \bar X+\sum_{i=1}^{n_0}F_i^Xz(t-\tau_i)+G_2Y\\\hspace{1cm}+\int_0^{\tau_{n_0}}H_X(\nu)z(t-\nu)d\nu+BU,\\
% \dot{Y}=\bar A_1Y+\sum_{i=1}^{n_0}F_i^Yz(t-\tau_i)\\\hspace{1cm}+\int_0^{\tau_{n_0}}H_Y(\nu)z(t-\nu)d\nu,
% \end{array} 
% \right.
% \end{equation}
% $\bar A_1$ is Hurwitz. It is shown in~\cite{auriol2024contributions} that the results obtained for eq

%\modifJA{I am not satisfied with this part. Federico, you may want to modify it to simplify the exposition.}

Interestingly, the problem of stabilizing an ODE-PDE-ODE system with PDE actuation or measurement has received limited attention in the literature. In \cite{Deutscher2018}, a state observer is developed for such a system, where the measurement corresponds to the PDE state $u(t,1)$. This is a challenging problem as the ODE at the unactuated boundary of the observer is, in this case, subject to a coupling with the PDEs from both boundaries. By making use of the multi-step approach, it is possible to design a state-observer under appropriate observability conditions. Since the observability problem is the dual of the controllability problem, the properties of the observer operator can be readily adapted to those of the controller operator. Recently, an alternative approach has been proposed in~\cite{Auriol2023,auriol2024contributions} using the time-delay representation of the system~\eqref{Sec_4_1_v_pde_orig_simp}.

 %%%%%%%%%%%%%%%%%%%%%%%%%%%%%%%%%%%%%%%%%%%%%%%
        \subsection{Interconnected systems with a chain structure} 
%\modifJA{Again, we can shorten}
The framework proposed by the representation~\eqref{Sec_4_1_v_pde_orig_simp} is comprehensive and can describe various interconnected systems. As shown in~\cite{auriol2024contributions}, using a technique referred to as \textbf{folding} (see~\cite{Andrade2018,Auriol2020d,Auriol2022} for details), it can describe complex networks composed of interconnected hyperbolic/ODEs subsystems. However, such a representation shadows the cascade structure between the different subsystems, which may lead to strong requirements when designing stabilizing controllers. In particular, the contributions presented in Section~\ref{Sec_ODE_PDE_ODE} avoid configurations where the PDE system is under-actuated or under-measured. Consequently, they fail to design stabilizing controllers for interconnected PDE systems when only one of the PDEs is actuated. Consider the following interconnections of $N_0$ PDE of $n_i+m_i$ heterodirectional systems
\begin{equation}\label{Sec_4_1_v_pde_orig_simp_chain}
\left\{
\begin{array}{l}
    (w_i)_t+\Lambda_i(w_i)_x=\Sigma_i(x)w_i,\\
    u_i(0)=Q_iv_i(0)+Q_{i,i+1}u_{i+1}(1),\\
    v_i(1)=R_i u_i(1)+R_{i,i-1}v_{i-1}(0),
\end{array} 
\right.
\end{equation}
where $Q_{N_0,N_0+1}=0$ and where $R_{1,0}v_0=U$ is the control input. The different subsystems are interconnected through their boundaries in a chain structure, where only the extremity of the chain is actuated. In most contributions, the matrices $R_{i,i-1}$ are assumed to be full-row rank~\cite{auriol2024contributions}. Such a chain configuration is considered in~\cite{zhang2023robust} in the case of a scalar $1+1$ system subject to a delayed control input that can be modeled as a transport equation leading to the framework~\eqref{Sec_4_1_v_pde_orig_simp_chain} with $N_0=2$. The system is mapped to an exponentially stable target system using the following backstepping transformation 
\begin{equation*}
\left\{
\begin{array}{l}
    \gamma_2=w_2-\int_0^x K_{22}(x,y)w_2dy, \\
    \gamma_1=w_1-\int_0^x K_{11}(x,y)w_1dy-\int_0^1 K_{12}(x,y)w_2dy.
\end{array} 
\right.
\end{equation*}
Such a transformation is invertible due to its triangular structure. The first line of the transformation moves the in-domain couplings of the second subsystem at the boundary between the two subsystems (where they take the form of integral couplings). Then, the second line of the transformations is used to move these new integral coupling terms and the in-domain coupling terms of the first subsystem to the actuated boundary. %\modifJA{Sentence to make a link with classical backstepping?} 
Such a multi-step backstepping transformation is proposed in~\cite{Auriol2020d} to design a stabilizing output-feedback controller for a chain consisting of an arbitrary number of $1+1$ subsystems. %However, the proposed proof to show the well-posedness of the resulting kernel equations requires additional conditions on the boundary coupling matrices.
Although this approach has enabled significant advances, it only considers scalar subsystems and lacks
adaptability from one chain to another. For instance, adding one additional PDE subsystem in the chain structure is not directly possible. To overcome this limitation, an alternative recursive methodology is introduced in~\cite{Redaud2021}. Assuming several fundamental properties for each subsystem: output trajectory tracking, input-to-state stability, predictability (we can design predictors for the different
states), and observability, it is possible to simply design a recursive stabilizing output-feedback controller using elementary backstepping transformations. Again, the time-delay formulation of the interconnected system is at the core of the design. The proposed approach is modular in that additional subsystems can easily be included. Moreover, this framework can be extended to include ODEs at the end of the chain. The presence of ODE in the chain is also considered in~\cite{Aarsnes2019c} for a PDE-ODE-PDE configuration. The case of two interconnected $n_i+m_i$ heterodirectional systems is tackled in~\cite{Auriol2022} using similar techniques. Finally, this recursive approach is extended in~\cite{auriol2024contributions} to an arbitrary number of heterodirectional subsystems. However, in this latter configuration, the backstepping transformations now contain a time-affine integral term.

There are several situations in which the actuator may be located at an arbitrary node of the chain (and not at its extremity). For instance, when developing traffic control strategies on vast road networks, the actuator (ramp metering) can be located at a crossroad (junction of two roads)~\cite{Yu2023}. It is shown in~\cite{Redaud2022Fredholm} that having an actuator located at one of the intersection nodes of the chain raises challenging controllability questions. For the case of two interconnected $1+1$ subsystems actuated at the junction, it is possible to design a stabilizing backstepping controller. However, the transformation proposed in~\cite{Redaud2022Fredholm} corresponds to a Fredholm integral transformation, which is not always invertible. Interestingly, it is shown that spectral controllability conditions imply the existence and invertibility of such a transformation.

 %%%%%%%%%%%%%%%%%%%%%%%%%%%%%%%%%%%%%%%%%%%%%%%
  \subsection{Output regulation and disturbance rejection}
 % \cite{Anfinsen2015,Deutscher2017c,guo2021robust}%,Deutscher2017,deutscher2021backstepping,,Redaud2022b

%\cite{Aamo2013, Redaud2022b, irscheid2023output}

A natural extension of the analysis of coupled PDE-ODE systems is the consideration of output tracking and disturbance rejection problems with a disturbance or reference model. In this case, when the disturbance and output trajectories can be modeled by a finite-dimensional exosystem, backstepping designs exist for the disturbance rejection and output tracking problems. The main technical difference with respect to the results in the previous section is that the output regulation problem cannot be directly reformulated as a classical PDE-ODE or ODE-PDE-ODE stabilization problem, where the (extended) system is mapped to a stable target system by adequately choosing the backstepping transform. This is due to the fact that the state of the system will not (in general) converge to the origin. There are some results using backstepping techniques that have a first component of motion planning and then use backstepping to track the desired trajectory, for instance, \cite{coc2009}. In this section, we will focus instead on results based on output regulation and disturbance rejection techniques that do not require motion planning, since they are fundamentally different from stabilization problems (even around a reference trajectory).

The first result in trajectory tracking based on a backstepping design can be found in \cite{bresch2009adaptive}, even if the main focus of the paper is the use of adaptive techniques. In this paper, the system under consideration is an ODE with unknown input delay which, as has been previously discussed, is closely related to hyperbolic PDEs. A first result on disturbance rejection specifically for hyperbolic PDEs using the backstepping method can be found a few years later in \cite{Aamo2013}, for $2\times 2$ linear hyperbolic systems including a boundary disturbance generated by an autonomous linear ODE (exosystem). More precisely, the boundary disturbance is considered as acting on the boundary opposite the control action, and the regulated output is a linear function of the PDE states on the same boundary as the disturbance (yet linearly independent). This result was later extended to a regulated output defined as a linear function of the PDE states at an arbitrary point inside the domain in \cite{Anfinsen2015}, and subsequently to general heterodirectional transport equation systems with constant coefficients in \cite{Anfinsen2017}. The corresponding cases for spatially varying coefficients were studied in \cite{Deutscher2017}, \cite{Deutscher2017c} and \cite{Deutscher2017b}. In particular \cite{Deutscher2017c} links the corresponding design to a set of regulator equations depending on the localization of the regulated output. Two related results appearing in the same time period but on slightly different research lines are the disturbance rejection problem for 1D hyperbolic Partial Integro-Differential Equations (PIDEs) in \cite{Xu2017} and the consideration of output regulation for anti-stable wave equations in \cite{Gu2018coupledwave} and coupled linear wave-{ODE} systems in \cite{deutscher2021backstepping}. Other recent results have focused a bit more on the trajectory tracking problem, such as \cite{Redaud2022b} and \cite{redaud2024tracking} with a focus on linear hyperbolic systems and linear ODE exosystems, and \cite{irscheid2023output} focusing on nonlinear ODEs. 
%\modifFBA{Cooperative output regulation for networks + adaptive observers \cite{enderes2024cooperative} }

%\modifFBA{ADD CONFERENCE PAPERS OR NOT ???? Also, general refs not in line with backstepping are probably better not here ? (like \cite{guo2021robust}). Also there is some link with Rafa's (and florent's etc.) paper on hyperbolic systems but it was definitely not the focus and presented very differently...}

In a parallel track, some results have been obtained for analogous problems in the case of parabolic PDEs. In particular, \cite{Meurer2009} proposes a control design for the output-tracking problem of boundary-actuated parabolic PDEs with varying parameters, using both a backstepping and flatness approach. In \cite{Deutscher015outputparabolic}, the output regulation problem for a class of boundary-controlled parabolic PDEs is considered, where the disturbance can act both on the uncontrolled boundary of the PDE and the domain. Related to these results, other problems have also been considered, such as cooperative output-tracking problems that add the difficulty of time-varying structures, as in \cite{deutscher2021backstepping,Deutscher2021cooperativeoutputregulation}. It should be noted that several of the results focusing on output-tracking for parabolic systems via backstepping have followed developments in the hyperbolic cases, which probably explains the fact that papers in this area tend to have fewer and generally more consolidated results.

  %%%%%%%%%%%%%%%%%%%%%%%%%%%%%%%%%%%%%%%%%%%%%%%      

\subsection{Coupled parabolic problems}

Coupled parabolic problems, encompassing systems of multiple parabolic PDEs interconnected through various coupling mechanisms, appear frequently in modeling real-world phenomena. From chemical reactors and thermal processes to biological systems and population dynamics, these interconnected PDE systems pose unique challenges for control design due to the complex interactions between their spatially distributed states.  

\subsubsection{$n \times n$ systems}
Even though the initial results in backstepping dealt with parabolic problems, designs for coupled parabolic systems came later than their hyperbolic counterparts. Early attempts, such as~\cite{Baccoli2015}, focused on simplified cases with constant parameters and explored the use of scalar kernels for the backstepping transformation. However, this approach limited the design flexibility and achievable performance. Subsequent works, notably~\cite{Vazquez2017}, expanded the methodology to include spatially varying reaction terms. They recognized the limitations of scalar kernels and introduced matrix kernels for greater control design versatility. This allowed for the assignment of an arbitrary decay rate for the closed-loop system, addressing a broader range of unstable systems. A crucial aspect of~\cite{Vazquez2017} was the understanding of a profound connection between the kernel equations for coupled parabolic systems and those arising in the control of coupled hyperbolic systems. The authors observed that by defining auxiliary kernel matrices related to the spatial derivatives of the original backstepping kernel, one could transform the coupled second-order hyperbolic PDEs governing the original kernel into a system of first-order hyperbolic PDEs.

Specifically, for a system of $n$ coupled parabolic PDEs with the backstepping transformation defined by the kernel matrix $K(x,\xi) \in \mathbb{R}^{n \times n}$, the authors introduced the auxiliary kernels $L(x,\xi)$ and $R(x,\xi)$ as follows:
\begin{align}
    L(x,\xi) &= \sqrt{\Sigma(x)}K_x(x,\xi) - K_\xi(x,\xi)\sqrt{\Sigma(\xi)},\\
    R(x,\xi) &= \sqrt{\Sigma(x)}K_x(x,\xi) + K_\xi(x,\xi)\sqrt{\Sigma(\xi)},
\end{align}
where $\Sigma(x)$ is the diagonal matrix of diffusion coefficients. This transformation effectively converted the $n \times n$ system of 2nd order hyperbolic kernel equations for $K(x,\xi)$ into a system of $2n \times 2n$ first-order hyperbolic PDEs for $L(x,\xi)$ and $R(x,\xi)$.

Remarkably, the structure of these kernel equations is similar to the ones encountered in the backstepping control of coupled hyperbolic systems\cite{Hu2016a,Coron2013}. This connection unveils a deeper underlying mathematical structure shared by backstepping designs for apparently different classes of PDE systems. Subsequently, \cite{camacho2020boundary} addressed the issue of boundary observer design for coupled reaction-diffusion systems with spatially varying reaction terms, utilizing a similar connection. Expanding on these advancements, \cite{deutscher2018coupled}  proposed a comprehensive backstepping design for coupled linear parabolic PIDEs (Partial Integro-Differential Equations) with distinct diffusion coefficients, spatially varying parameters, and mixed coupled boundary conditions. This work systematically tackled the well-posedness of the complex kernel equations, which constituted a system of coupled second-order hyperbolic PIDEs. By extending the method of successive approximations, the authors presented a constructive solution approach for the kernel equations, culminating in a robust control design methodology for a broader class of coupled parabolic systems, later expanded to time-varying coefficients in \cite{ker2019}. 

While existing results address fully-actuated boundary-coupled parabolic systems, the control of underactuated coupled systems---where the number of controls is less than the number of states---remains an open challenge and an important direction for future research.

% the structure of these first-order hyperbolic PDEs closely resembled the kernel equations encountered in the backstepping control of coupled hyperbolic systems, particularly those investigated in \cite{Hu2016a,Coron2013}. This connection revealed a deeper underlying mathematical structure shared by backstepping designs for seemingly disparate classes of PDE systems. Later, \cite{camacho2020boundary} tackled the problem of boundary observer design for coupled reaction-diffusion systems with spatially varying reaction terms usign a similar connection.

% Furthering these advancements, \cite{Deutscher2018} proposed a comprehensive backstepping design for coupled linear parabolic PIDEs (Partial Integro-Differential Equations) with distinct diffusion coefficients, spatially varying parameters, and mixed boundary conditions. This work systematically addressed the well-posedness of the complex kernel equations, which formed a system of coupled second-order hyperbolic PIDEs. By extending the method of successive approximations, the authors provided a constructive solution procedure for the kernel equations, culminating in a robust control design methodology for a broader class of coupled parabolic systems.

\subsubsection{Couplings of parabolic PDEs with ODEs and chains}
%The interconnection of parabolic PDEs with ODEs introduces new challenges and opportunities for control design, as the dynamics of finite-dimensional systems interact with the spatially distributed behavior of PDEs. 
The interplay between parabolic PDEs and ODEs arises naturally in many applications, where actuator and sensor dynamics, or other lumped parameter subsystems, influence the behavior of systems governed by parabolic PDEs. Initial efforts, such as \cite{krs2009}, considered the control of a single diffusion PDE coupled with a linear ODE.  The authors employed a backstepping transformation incorporating both the PDE and ODE states, leading to a target system that could be stabilized using a boundary controller. The design focused on compensating the effect of the ODE dynamics on the PDE, achieving exponential stabilization of the combined system. A more complex scenario involving a diffusion PDE sandwiched between two ODEs was investigated in \cite{Tang2011}. The system represented a plant with both actuator and sensor dynamics governed by ODEs, coupled with a diffusion process occurring within the plant. The control design involved a two-step backstepping procedure, addressing first the actuator dynamics and then the sensor dynamics.  The resulting controller achieved stabilization of the sandwiched PDE and guaranteed the exponential convergence of the observer error.

Finally, \cite{Deutscher2020} presented a comprehensive backstepping design for coupled parabolic ODE-PDE-ODE systems, where multiple parabolic PDEs with distinct diffusion coefficients and spatially varying parameters were bidirectionally coupled with ODEs at both boundaries. The authors employed a successive backstepping approach, systematically decoupling the system into a stable ODE-PDE-ODE cascade. The design involved solving a sequence of kernel equations, initial value problems, and boundary value problems, enabling the specification of the closed-loop stability margin. This work significantly broadened the scope of backstepping designs for interconnected parabolic PDEs and ODE systems, addressing a wider range of complex coupling structures and boundary conditions.

On the other hand, chains of parabolic equations have not been the object of research until recently. In \cite{Xu2023}, the authors investigate the stabilization of a chain of parabolic PDE-ODE cascades, where the output of each ODE subsystem acts as the control input to the preceding PDE subsystem, and the state of each PDE subsystem enters as control into the corresponding ODE subsystem. This work extends previous results, which primarily focused on chains of scalar ODEs interconnected with pure delay PDEs, by considering more general scenarios. The authors allow for the virtual inputs (the outputs of the ODE subsystems) to be affected by PDE dynamics that are not restricted to pure delays. They include diffusion, counter-convection, and spatially distributed coupling terms in the PDE subsystems, allowing for a broader range of actuator dynamics. Additionally, the ODE subsystems are not limited to scalar ODEs in a strict-feedback configuration; instead, they can be general linear time-invariant (LTI) systems.

\subsubsection{Hyperbolic-Parabolic Problems}

The interconnection of hyperbolic and parabolic PDEs leads to mixed-class systems that capture the interplay between transport and diffusion phenomena. These systems arise in diverse applications, such as fluid-structure interactions, thermoacoustic instabilities, and chemotaxis models, where the finite speed of propagation inherent to hyperbolic dynamics interacts with the smoothing and dissipative effects of parabolic dynamics.  Interestingly, the designs give rise to kernel equations that also exhibit a hyperbolic-parabolic structure. In particular, one of the backstepping kernels dealing with the parabolic state is defined on a square domain and verifies a parabolic equation, a unique feature not normally encountered in backstepping designs.

An early exploration of backstepping control for a hyperbolic-parabolic system can be found in \cite{krstic2009control}. The author considered a reaction-diffusion PDE (parabolic) with a long input delay, which was represented as a transport PDE (hyperbolic). The combined system, although appearing as a cascade, exhibited a challenging coupling structure due to the delayed boundary control action affecting the parabolic PDE. Using backstepping, a stabilizing controller was designed, compensating for the delay and achieving stability in the sense of a combined norm incorporating the states of both PDEs.

Extending these concepts, \cite{Chen2017} investigated the stabilization of an underactuated coupled transport-wave PDE system. The system, motivated by models of acoustic waves interacting with vibrating structures, involved a first-order transport PDE (hyperbolic) coupled with a wave equation (hyperbolic) through a boundary condition. However, the wave equation could be reformulated as a system of two first-order transport equations, allowing for a unified representation as a hyperbolic-parabolic system. The authors employed backstepping to design a boundary controller, achieving exponential stability despite the underactuation.  A similar problem, with bidirectional couplings between a transport equation and a diffusion equation, was also explored in \cite{ghousein2020backstepping}. The authors allowed for spatially distributed coupling terms within the domain of each PDE and designed a backstepping controller that achieved exponential stability for a broad class of interconnected systems. The kernel equations, reflecting the mixed nature of the system, were also of a hyperbolic-parabolic type, requiring sophisticated analytical tools to prove their well-posedness.

A significant contribution towards tackling more general hyperbolic-parabolic PDE-PDE systems was made in \cite{Chen2023b}. This work addressed a system with boundary-dependent and in-domain direct coupling terms within the hyperbolic system domain, making controller design considerably more challenging.  The authors applied a novel backstepping transformation,  leading to a set of coupled kernel equations that reflect the mixed hyperbolic-parabolic nature of the original system.   The well-posedness of these coupled kernel equations was proven using an innovative infinite induction energy series, and the invertibility of the transformations was also established. This approach facilitated the design of a boundary controller that achieved exponential stabilization of the mixed PDE system in the $H^1$ sense. 

Finally, \cite{deutscher2023backstepping} addressed the challenge of input and output delays in coupled linear parabolic PDEs. The authors introduced a transport PDE to model the delays, resulting in a general family of parabolic-hyperbolic systems. They addressed both controller and observer design, overcoming the complexities introduced by the non-local and delayed nature of the system. This work highlighted the applicability of backstepping to increasingly sophisticated interconnected PDE systems.

\subsubsection{Elliptic-parabolic problems}
Elliptic-parabolic problems arise in various applications where the steady-state behavior of one physical process (described by the elliptic equation) influences the transient behavior of another (described by the parabolic equation). While seemingly less common than purely parabolic or hyperbolic systems, these mixed-class systems present unique challenges for control design due to the inherent differences in their mathematical structure and underlying physical phenomena. 

Interestingly, elliptic-parabolic couplings often emerge as a consequence of employing singular perturbation techniques to simplify the analysis of more complex PDE systems. For instance, in \cite{Vazquez2006,convloop}, the authors tackled the stabilization of a thermal convection loop. Using singular perturbation theory, assuming a large Prandtl number, they derive a reduced model in the form of a coupled parabolic-elliptic system which is made strict-feedback with one of the controls. Backstepping was then employed to design a boundary controller for the temperature, leading to an exponentially stable closed-loop system with an elliptic-parabolic structure. The control of flexible structures, particularly the Timoshenko beam model, has also involved elliptic-parabolic couplings. In \cite{Krstic2006}, the authors designed a backstepping boundary controller for a slender Timoshenko beam with Kelvin–Voigt damping. The control design employed a singular perturbation approach, exploiting the smallness of the rotational inertia coefficient, and the resulting closed-loop system exhibited an elliptic-parabolic coupling.

The control of Navier-Stokes and magnetohydrodynamic (MHD) channel flows, studied respectively in \cite{Vazquez2007} and \cite{Xu2008}, also present elliptic-parabolic couplings. In both cases, the pressure verifies an elliptic PDE. In the case of MHD, there is an electrically conducting fluid with an imposed magnetic field which also verifies an elliptic PDE for the electric potential. Backstepping was applied to design a boundary controller that stabilized the velocity field by using the similar principle of solving the elliptic equation and using additional controls to make the solution strict-feedback and thus amenable to backstepping.

A more recent contribution specifically addressing the control of an unstable parabolic-elliptic system can be found in \cite{Alalabi2025}. This work considered a system where the parabolic equation, despite being exponentially stable on its own, became unstable due to coupling with the elliptic equation. The authors employed backstepping to design a boundary controller, mapping the original system into an exponentially stable target system with an elliptic-parabolic structure. This approach resulted in an explicit control law and provided a constructive method for stabilizing this challenging class of coupled PDE systems for some values of the parameters, in the case of a single control.

        %%%%%%%%%%%%%%%%%%%%%%%
        \subsection{Bilateral control of PDEs} 

Across all preceding configurations, we mainly focused on interconnected systems with control or measurement at one end of the domain. Even when dealing with chains of systems, we always assumed that one and only one node of the network describing the chain was fully actuated. In various scenarios, there might be instances where actuators are available at multiple nodes of the chain. For instance, in the context of underbalanced drilling~\cite{aarsnes2016methodology}, sensors are typically located at the surface and, in more modern drilling facilities, at the bottom of the well, which induces a bilateral estimation problem. Having actuators at several nodes induces over-actuated systems and these additional actuators can be used to improve the closed-loop performance of the system (robustness to noise, smaller control effort, better convergence rate). We present here the results obtained for the stabilization of PDE systems with actuation at both ends of the spatial domain. We refer to these situations as \emph{bilateral control}. These bilateral cases are nontrivial extensions as they require specific backstepping transformations.

        \subsubsection{Parabolic case}\label{sec-bil-par}
       The problem of bilateral control for reaction-diffusion equations was formally addressed first in \cite{Vazquez2016}, where the authors proposed a backstepping design for 1-D parabolic PDEs with actuation at both ends. This design considers a different backstepping transformation, where the integral starts at the middle of the domain, as opposed to the classical transformations that start at the uncontrolled boundary. This subtle change, which was inspired by the N-dimensional design in balls for $N=1$ (see Section~\ref{sec-higher}) allows solving the bilateral control problem for a general class of reaction-diffusion equations with spatially varying reaction terms.

A possible way to look at the bilateral backstepping transformation is to divide the domain into two subdomains and define separate backstepping transformations for each subdomain. These transformations are then coupled through compatibility conditions at the middle point of the domain, ensuring the continuity of the state and its spatial derivative. The resulting kernel equations are well-posed and can be solved using the method of successive approximations. The feedback laws are then obtained by evaluating the transformations at the controlled boundaries, resulting in integral operators acting on the state of the system.

Considering a domain of the shape $[-L,L]$, the bilateral backstepping transformation introduced in \cite{Vazquez2016} for a problem similar to the one of Section~\ref{sec-plant} but with controls $U_1$ and $U_2$ at each side of the domain is
\begin{equation}
w(t, x) = u(t, x)-\int_{-x}^{x} K(x, \xi)u(t, \xi)d\xi, \label{eqn-bilateral_transform}
\end{equation}
 The kernel $K(x, \xi)$ satisfies a hyperbolic PDE in an hourglass-shaped domain instead of a triangle; this PDE can be split into two separate PDEs defined on triangular subdomains with coupled boundary conditions. The control laws are obtained from evaluating the transformation \eqref{eqn-bilateral_transform} at the boundaries, resulting in integral operators acting on the state $u$:
\begin{align}
U_1(t) &= \int_{-L}^{L} K(L, \xi)u(t, \xi)d\xi, \label{eqn-bilateral_control1}\\
U_2(t) &= \int_{-L}^{L} K(-L, \xi)u(t, \xi)d\xi. \label{eqn-bilateral_control2}
\end{align}
The paper \cite{Vazquez2016}  also shows that the bilateral control problem can be reformulated as a control problem for a $2 \times 2$ coupled reaction-diffusion system. This is achieved by ``folding'' the original system at the middle point of the domain, which allows using existing backstepping results for coupled parabolic PDEs (see Section~\ref{sec-higher}).

The bilateral backstepping design was also used in \cite{bekiaris2019nonlinear} to address the trajectory tracking problem for a class of viscous Hamilton-Jacobi PDEs, using bilateral actuation and sensing. The design employed a combination of feedback linearizing transformation and backstepping, resulting in explicit observer-based output-feedback laws that achieved local asymptotic stabilization of the desired trajectory. The feasibility conditions imposed by the linearizing transformation (its invertibility is only local) highlighted the regional nature of the stability result (see Section~\ref{sec-nonlinear} for a similar unilateral result on the Burgers equation). A more recent paper \cite{chen2021folding} further explores the folding technique for designing bilateral backstepping output-feedback controllers for an unstable parabolic PDE. The authors introduce the notion of a folding point that can be arbitrarily chosen within the domain, leading to a more general formulation of the bilateral control problem. They also design a state observer based on collocated measurements at an arbitrary interior point, generating exponentially converging state estimates. The combined output feedback scheme, composed of the state-feedback and observer, is shown to stabilize the system. This paper also includes a numerical analysis of how the selection of folding points and measurement points affect the controller and observer responses, illustrating the potential for performance optimization using these design parameters.

        %%%%%%%%%%%%%%%%%%%%%%%%%%%%%%%
        \subsubsection{Hyperbolic case}

Let us now consider system~\eqref{Sec_4_1_v_pde_orig_simp} without the ODE terms and with the boundary conditions
$u(0)=Qv(0)+U_1(t),~v(1)=R u(1)+U_2(t)$ and the measurement $y_1(t)=v(t,0)$ and $y_2(t)=u(t,1)$. For this bilateral configuration, it has been shown in~\cite{li2010strong} that it is always possible to design a stabilizing feedback controller that stirs the state to zero in a minimum finite-time $t^1_F=\max\{\frac{1}{\lambda_1},\frac{1}{\mu_1}\}$ that corresponds to largest time needed for the slowest characteristic in each direction to travel the entire spatial domain. Analogously, it is also possible to design a state-observer that reconstructs the state of the system in the same minimum finite-time $t_F^1$. Explicit finite-time backstepping controllers and observers are proposed in~\cite{Auriol2018b} for this bilateral configuration. The finite-time stable target system is chosen as 
\begin{equation}\label{Sec_4_two_side_target}
\left\{
\begin{array}{l}
    \gamma_t+\Lambda \gamma_x=\begin{pmatrix}
        \Omega(x) & \Gamma(x) \\ \bar \Omega(x) & \bar \Gamma(x)    \end{pmatrix}\gamma, \\
        \alpha(t,0)=0,~\beta(t,1)=0,
\end{array} 
\right.
\end{equation}
where $\Omega$ and $\Bar \Omega$ are upper-triangular continuous matrices, and $\Gamma$ and $\bar \Gamma$ are continuous matrices with a specific cascade structure, which guarantee the finite-time stability in time $t_F^1$ of system~\eqref{Sec_4_two_side_target}. Their presence is necessary to prove the well-posedness of the backstepping transformation. The proposed backstepping transformation in~\cite{Auriol2018b} is given by
\begin{align} \label{Sec_4_bilateral_BS}
    \gamma(t,x)=w(t,x)-\int_0^1Q(x,y)w(t,y).
\end{align}
Although~\eqref{Sec_4_bilateral_BS} looks like a Fredholm transformation, the kernels $Q$ given in~\cite{Auriol2018b}  are actually defined on sub-triangles of the unit square $[0,1]^2$ and the transformation could be rewritten as a convoluted Volterra transformation, thus guaranteeing its invertibility. However, the invertibility was proved in a different way in~\cite{Auriol2018b} using an operator framework inspired from~\cite{Coron2016}. The corresponding kernel equations are solved recursively using analogous arguments as the ones developed~\cite{Auriol2016} and the cascade structure of $\Gamma$ and $\bar \Gamma$. 
% All in all, the control laws read as
% \begin{align*} 
% U_1(t)&=-Qv(t,0)+\begin{pmatrix}
%         Id\\ 0
%     \end{pmatrix}^TQ(0,y)w(y)dy,\\
%     U_2(t)&=-Ru(t,1)+\int_0^1 \begin{pmatrix}
%         0 \\ Id
%     \end{pmatrix}^TQ(1,y)w(y)dy.
% \end{align*}
The state observer can be designed using a dual operator approach~\cite{Auriol2018b} to obtain a stabilizing output-feedback controller.

For the observer design problem, a more straightforward alternative approach is proposed in~\cite{wilhelmsen2021minimum}, where the original $n + m$ system with bilateral sensing is transformed to an appropriate $(n + m) + (n + m)$ system with unilateral sensing via an invertible folding coordinate transformation. This folding transformation is chosen such that the sum of the transport times needed for the slowest characteristic in each direction to travel the entire spatial domain corresponds to $t_F^1$. It is then possible to apply the results of Section~\ref{Sec_4_1_2}.

%%%%%%%%%%%%%%%%%%%%%%%%%%%%%%%%%%%
         \subsection{Timoshenko beams and wave equations} \label{sec-tim} 
Recent designs for control of Timoshenko beams and wave equations using backstepping rely on a preliminary step that involves a Riemann-like transformation to map the original system into a system of first-order hyperbolic PDEs, potentially coupled with ODEs. This transformation greatly facilitates the application of backstepping designs developed for interconnected systems that have been explored throughout this section,  allowing for the design of controllers that achieve rapid stabilization, even in the presence of destabilizing boundary conditions and complex nonlocal terms.

First used formally in \cite{Su2017}, the authors considered the wave equation with velocity recirculation, where the time derivative of the state at one boundary was fed back into the PDE domain. This nonlocal term destroyed the passivity of the wave equation, making traditional control methods inapplicable. The authors used a Riemann transformation to convert the system into a cascade of two transport PDEs and an ODE, and then designed a backstepping controller that achieved exponential stability.

The application of Riemann transformations to Timoshenko beam control was further explored in \cite{Chen2023}. The authors addressed the problem of rapid stabilization for a Timoshenko beam with anti-damping and anti-stiffness at the uncontrolled boundary. They used a Riemann transformation to convert the Timoshenko beam model into a 1-D hyperbolic PIDE-ODE system and employed a backstepping approach to design a controller achieving an arbitrarily fast decay rate. More recently, in \cite{chen2023block}, the authors extended the backstepping design to multilayer Timoshenko beams, addressing the challenge of isotachic states (see Section~\ref{Sec_isotachic}), where multiple states share the same transport speed. They employed a Riemann transformation to convert the multilayer beam model into a 1-D hyperbolic PIDE-ODE system and used a block diagonalization approach to handle the isotachic states. This resulted in a block backstepping controller that achieved rapid stabilization with an arbitrary decay rate.
%%%%%%%%%%%%%%%%%%%%%%%%
\section{Other advanced PDE topics}\label{Sec_advanced_topics}

In this section, we navigate through a collection of advanced and multifaceted problems in the realm of PDEs. These issues, each meriting detailed individual analysis, are unified here for a succinct yet comprehensive exploration, reflecting their complexity and diverse nature.

%%%%%%%%%%%%%%%%%%%%%%%%%%%%%%%%%%%%%%%%%%%%        
        
\subsection{Adaptive control and parameter estimation} \label{sec-adaptive}
As seen in the previous sections, the backstepping technique has provided constructive and sometimes explicit control laws for large classes of (mostly linear) PDE systems. It has also enabled the use of other advanced methods and multi-step design techniques by simplifying the structure of closed-loop PDE systems. One of the first such uses of the backstepping methodology was the extension of adaptive control designs to new classes of PDEs. These designs leveraged the simpler structure of target systems obtained via the backstepping transform and, when available, the explicit nature of the transform, in order to apply different adaptive strategies to the control of both parabolic and hyperbolic PDEs. These strategies also allow for the estimation of boundary parameters in certain classes of systems and are not limited to in-domain parameters. The use of backstepping-based adaptive designs appeared in \cite{Smyshlyaevfirstadaptivecdc}, applied to 1D parabolic PDEs and was rapidly expanded in the series of papers \cite{KrsticadaptivepartI,SmyshlyaevadaptivepartII,SmyshlyaevadaptivepartIII}, as well as extended to higher-dimensional parabolic PDEs in \cite{krsticannualreviewsadaptive}. These and other designs were collected in the book \cite{adaptive}. In the case of delay systems, closely related to hyperbolic PDE systems, the first results for unknown delay were obtained in \cite{krstic2009delay,bresch2010delay} for and ODE with unknown input delay and extended to the problem of trajectory tracking for an output of the ODE in \cite{bresch2009adaptive}. As for hyperbolic PDEs, the first result is \cite{Krstic2010adaptivewave}, followed a few years later by \cite{Bernard2014} for a 1D transport PDE with non-local terms. Closely thereafter, an adaptive output-feedback control law for a wave equation with dynamic boundary appeared in \cite{BreschPietri2014a}. After these initial results, many papers on adaptive schemes based on the backstepping technique have appeared in the literature, for instance, for general linear hyperbolic PDEs in \cite{Anfinsen2016a}. A good overview of these designs can be found in \cite{anfinsen2019adaptive}.

\subsection{Extremum seeking for PDEs} \label{sec-ES}

In the book \cite{oliveira2022extremum}, Oliveira and Krstic assemble a multitude of algorithms for seeking extrema and Nash equilibria (in non-cooperative multi-agent scenarios) in the presence of parabolic and hyperbolic PDEs, including delays. While the presence of stable PDEs can be ignored by resorting to sufficiently slow seeking, when the PDEs are compensated with PDE backstepping the limit on the convergence speed is removed and the convergence rate is assignable with Newton algorithms.

\subsection{Extension to higher-dimensional domains}\label{sec-higher} 
The original backstepping design for PDEs, as explained in Section~\ref{Sec_primer}, focused on a one-dimensional reaction-diffusion equation. This limitation to 1-D systems stemmed from the challenge of deriving and solving the kernel equations for higher-dimensional domains, where the complexity of the integral transformations and the associated PDEs governing the kernels increases substantially. Despite the inherent difficulties,  there have been considerable efforts in extending the backstepping method to higher-dimensional domains.  These extensions have primarily focused on specific geometries and boundary conditions, leveraging the symmetry or special properties of the domains to simplify the kernel equations and enable the design of explicit controllers. 

For many models in fluid dynamics and other areas, the design of
boundary controllers and observers can be greatly simplified by
exploiting the concept of \emph{spatial invariance}~\cite{bamieh2002distributed}.  This property, also
known as translational invariance, arises when the system dynamics
and geometry are invariant under translations in one or more
spatial coordinates. This allows for a reduction in the system's
dimensionality by transforming the spatially invariant coordinates
into a parameter or set of parameters. %Essentially, spatial invariance allows substituting spatial derivatives and spatial dependence for algebraic multiplication and dependence on a new parameter (or set of parameters) that replaces the spatially invariant coordinates. This transforms the
%original PDE system into a family (or ensemble) of parameterized lower dimensional systems, which are easier to analyze and control.
The controllers and observers can then be designed for these
simpler systems and subsequently reconstructed in the original
physical space. This procedure requires that the system possess specific properties, namely that the spatially invariant coordinates lie within a
group (usually the real line $\mathbb{R}$ or a periodic domain), actuation
and sensing are fully distributed over those coordinates, and the
dynamics and geometry are invariant with respect to translations
in them. By applying tools like the Fourier transform (for infinite
domains) or Fourier series/spherical harmonics (for periodic domains), one can convert the spatially invariant coordinates into wave numbers or harmonics, transforming
spatial derivatives into algebraic multiplications. This simplifies
the control design and analysis, allowing for the application of
techniques like backstepping to each wave number separately. The earliest designs based on these ideas were for flow control applications, and are explored in Section~\ref{sec-flowcontrol}. Later, \cite{Vazquez2016disk} and \cite{vazquez2019sphere} introduced backstepping designs for reaction-diffusion equations on, respectively, a 2-D disk and a 3-D sphere. The authors exploited the radial symmetry of these domains to simplify the kernel equations and obtained explicit control laws that achieved exponential stabilization.  These results were generalized in \cite{vazquez2016explicit} to address the control of reaction-diffusion equations on a n-dimensional ball, employing spherical harmonics and  Bessel functions to derive the backstepping kernel and feedback gains. Interestingly, the 1-dimensional case can be exploited to obtain the bilateral design described in Section~\ref{sec-bil-par}. More recently, \cite{R.2022ball}  presented a backstepping design for a reaction-diffusion system in a ball of arbitrary dimension with radially varying reaction using the power series method to overcome the singularities in the kernel equations. Other works exploring higher-dimensional designs include~\cite{liu2020boundary}, which exploits symmetry to apply the classic backstepping designs to simple geometries, and the book \cite{meurer2}, which provides a comprehensive overview of backstepping designs for higher-dimensional PDEs in several domains, including those for plates, shells, and other flexible structures. 

\subsection{Nonlinear PDEs}\label{sec-nonlinear}
While the backstepping methodology has achieved remarkable success in the realm of linear PDEs, its extension to nonlinear PDEs poses significant challenges. The nonlinear nature of these systems often precludes the direct application of linear techniques, necessitating the development of innovative approaches to address the inherent complexities.  One common strategy for controlling unstable nonlinear PDEs is to linearize the system around the desired equilibrium point and then design a feedback controller for the resulting linear model. This approach, employed in works like \cite{Coron2013,vazquez-nonlinear3}, guarantees local exponential stability in a neighborhood of the equilibrium. 

However, it inherits fundamental limitations of linearization. Stability properties are, at best, valid only for small perturbations around the linearization point, and in some instances, the linearized system can be exponentially stable while the nonlinear PDE system lacks even local stability.

%However, it inherits the limitations of linearization, as the stability properties are valid only for small perturbations around the linearization point. 
In certain cases, a linearizing transformation may exist that allows for the conversion of the nonlinear PDE into a linear one. This approach, showcased in \cite{Krstic2009,krstic2008nonlinearB} for the Burgers equation and in \cite{bekiaris2019nonlinear} for a Hamilton-Jacobi nonlinear PDE, relies on finding a suitable linearizing change of variables, which in this case is the Cole-Hopf transformation \cite[p. 206]{Evans2010}. The resulting linear system can then be stabilized using standard backstepping methods, yielding a controller that is nonlinear in the original state. The success of this approach hinges on the existence of such a transformation, which is not guaranteed for all nonlinear PDEs and is typically not global. For systems involving a nonlinear ODE coupled with a PDE, backstepping can sometimes be employed to compensate for the PDE dynamics and achieve stabilization.  Several examples are presented in the works \cite{BekiarisLiberis2014,BekiarisLiberis2013,BekiarisLiberis2016}, where the authors address the problem of compensating wave or delay dynamics for nonlinear ODE systems. The design involves a predictor-based feedback law, compensating for the dynamics introduced by the different PDEs. Interestingly, in some cases, a linear feedback law may suffice to stabilize a nonlinear PDE.  In \cite{karafyllis2019global}, the authors present a global stabilization result for a class of nonlinear reaction-diffusion PDEs using boundary feedback. The key insight is that for certain nonlinearities, a simple linear feedback law can induce a ``damping'' effect when the norm of the state becomes large, leading to global exponential stability in the $L^2$ norm. However, the most general and ambitious approach for nonlinear PDE control is the one based on \emph{spatial} Volterra series\footnote{A Volterra series is a functional series representation of a nonlinear system, generalizing the concept of a single Volterra integral. While a single Volterra integral captures the output of a system as a weighted sum of input values, a Volterra series extends this notion by incorporating multiple integrals of increasing order, capturing higher-order nonlinear interactions between input values.}.  The works \cite{vaz2008_a,vaz2008_b} introduce a framework that employs spatial Volterra series nonlinear operators both in the state transformation and the feedback law, generalizing the concept of feedback linearization to infinite dimensions. This approach allows for the stabilization of a broad class of nonlinear parabolic PDEs, but its complexity and the inherent challenges in proving the global invertibility of the transformation (in general, only local invertibility can be proved) have limited its application to specific cases. 

\subsection{Time-Varying PDEs and Prescribed/Fixed-Time Stabilization}

The extension of backstepping to time-varying PDEs, i.e. where the coefficients in the PDE depend on time, introduces significant challenges compared to the time-invariant case. Beyond a few cases where explicit solutions are available~\cite{smy2005b}, the primary obstacle stems from the complex interplay between the time-varying coefficients and the spatial dynamics of the PDE. For parabolic PDEs,  \cite{vazquez2008control,ker2019} demonstrate that the backstepping method can be applied to systems with time-varying coefficients, but the stability analysis requires imposing restrictions on the regularity of these coefficients. Specifically, the time-varying parameters are typically assumed to belong to certain Gevrey classes, which guarantees the existence of solutions for the kernel equations with specific smoothness properties. This reliance on Gevrey regularity stems from the inherent smoothing properties of parabolic PDEs, where the solution at any time instant is influenced by the entire past history of the system. An implementation to solve the resulting kernel equations is given in~\cite{jadachowski2012efficient}. In contrast, for hyperbolic PDEs, the situation is less restrictive, as highlighted in \cite{deu2016,anfinsen2020stabilization,Coron2021,coron2021boundary}.  The finite speed of propagation inherent to hyperbolic dynamics allows for a more localized influence of the time-varying coefficients. As a result, backstepping designs for time-varying hyperbolic PDEs can be developed without relying on Gevrey regularity, as the solution at a given point in spacetime is affected only by a limited portion of the time-varying parameter history.  Recently, alternative gain-scheduling approaches have been proposed in~\cite{karafyllis2021event,auriol2024event}. Instead of handling time-varying backstepping kernels that account for both temporal and spatial variations in the coefficients, these approaches consider only the spatial variations of the coupling coefficients. This method remains valid as long as the coupling coefficients are sampled over time, thereby reducing the kernel PDEs to involve only space-dependent coefficients between successive sampling or triggering instants. These methods do not require solving time-varying backstepping kernel equations in real-time resulting in reduced computational burden and broader applicability.

The goal of achieving stabilization within a prescribed time, independent of initial conditions, further complicates the design.  This is where the concepts of \emph{fixed-time} and \emph{prescribed-time} stabilization become relevant. \emph{Finite-time stability} refers to systems that reach equilibrium within a finite time, but this time depends on the initial conditions. Fixed-time stability, a stronger notion, guarantees convergence to the equilibrium within a time that is uniformly bounded with respect to the initial conditions.  In prescribed-time stability, the settling time can be explicitly prescribed a priori, independently of both initial conditions and system parameters. For parabolic PDEs, achieving finite-time or fixed-time stability typically requires the use of time-varying feedback gains.  In \cite{coron2017null}, this is achieved for space-varying reaction-diffusion systems by constructing a sequence of backstepping kernels mapping the  plant into a sequence of target systems with successively increased damping coefficients, which allows to designing periodic time-varying
feedback laws. Alternatively, in \cite{espitia2019boundary},  a single backstepping transformation with a time-varying kernel and a single target system with a time-varying coefficient is considered, ensuring convergence within a predetermined time, but the approach only covers constant-coefficient reaction-diffusion systems. Similarly,  \cite{espitia2022sensor,espitia2021predictor} focuses on prescribed-time stabilization and estimation of LTI systems with input delay, employing a time-varying predictor feedback based on a transport PDE representation of the delay and time-varying backstepping kernels. These approaches rely on selecting a time-varying target system that exhibits the desired finite-time convergence property.  For parabolic PDEs, this is essential because these systems do not naturally exhibit finite-time convergence to the origin. In contrast, hyperbolic equations can achieve finite-time convergence to zero, but with a minimum time determined by the propagation speed of the characteristics. This distinction highlights the fundamental difference between the temporal behavior of these two classes of PDEs.  In a related context,  \cite{steeves2020prescribed} presents prescribed-time observers and controllers for the linearized Schrödinger equation using time-varying gains, also based on backstepping. The observer gain diverges as time approaches the prescribed terminal time, but this growth is carefully managed to guarantee prescribed–time convergence of the estimate. The authors also demonstrate prescribed–time output regulation by combining the estimation and control results.

\subsection{Moving Boundary Domains}\label{sec-movingboundary}

The control and estimation of PDEs with time-varying spatial domains, known as moving boundary problems, present considerable challenges beyond those encountered in fixed-domain scenarios. The moving boundaries introduce time-dependent spatial domains and boundary conditions, leading to a more complex mathematical framework and necessitating specialized control and observer designs. Backstepping has shown remarkable flexibility in adapting to these challenges, with several contributions showcasing its efficacy in tackling moving boundary problems. One of the earliest works in this domain is \cite{izadi2015backstepping}, where an output feedback control design for moving boundary parabolic PDEs is proposed using backstepping.  The authors introduce a new coordinate transformation that incorporates both the state of the system and the moving boundary, leading to a fixed-domain target system but a time-varying PDE. A backstepping observer is designed to estimate the state using boundary measurements, and the estimated state is then used in a feedback law that compensates for the moving boundary effects. The problem of boundary stabilization of one-dimensional cross-diffusion systems in a moving domain is studied in \cite{CAUVINVILA2023251}. In~\cite{kog2018} the authors address the stabilization of the one-phase Stefan problem, a classic moving boundary problem describing melting or solidification processes. They employ a backstepping transformation that explicitly accounts for the moving boundary, mapping the Stefan problem into a stable target system with a fixed boundary. This approach leads to a boundary controller that effectively regulates the interface position and ensures exponential stability in a suitable norm. Building upon this work, \cite{kog2020_f} extends the backstepping design to the broader class of two-phase Stefan problems. Recently, \cite{szc2022} explores adaptive control in moving boundary problems for the wave equation, by rewriting the problem into an LTI form with time-varying input delay and unknown parameters to which adaptive backstepping designs can be applied. See also Section~\ref{sec-appmovingboundary} for some applications of these designs.

\subsection{Fredholm operators and non-strict feedback plants}\label{sec-fredholm}
As mentioned in the primer, in Section \ref{Sec_primer}, while the original backstepping method was developed using Volterra integral transformations, see \eqref{eqn-backstepping}, in some cases this does not give enough degrees of freedom to find an adequate target system and well-posed kernel equation for the control design. The first specific example where a different type of transform was looked for is likely \cite{SmyshlyaevEulerBernoulli2009}, to assign an arbitrary decay rate for an Euler-Bernoulli beam equation. This early example builds a control law based on a previous design for the Schrödinger equation, which uses a standard Volterra transformation, and then adapts it to obtain stabilization to zero of the Euler-Bernoulli beam equation instead of to a constant. Unlike later backstepping designs, this design depends on the system dynamics when defining the properties of the transforms, these are not seen as operators mapping from generic function spaces (e.g. $L^2$ to $L^2$), but relating solutions to two specific differential equations. This is reflected in the transformation structure, where the time derivative of the PDE state appears alongside the value of the state itself. In this particular case, the Fredholm transformation linking the two sets of solutions is found after imposing a particular form of the control action in order to obtain explicit expressions for the closed-loop system. The next examples in the literature are \cite{BEKIARISLIBERIS2010DistributedWave}, with a wave PDE in the actuation path of an ODE, and \cite{bekiaris2011lyapunov}, in the case of distributed input delays in an ODE. In these examples, the direct backstepping transformation is written in a more general form.  Nevertheless, this transformation, due to the particular structure of the system under consideration, has good invertibility properties. In particular, the inverse transform only contains Volterra integrals, which seems to indicate that the transformation is not a generic Fredholm transform. Another particular structure, this time concerning parabolic equations with integral terms, was presented in \cite{Guo2014Fredholm} for a spatially non-causal reaction-diffusion equation. In this case, the particular form of the integral kernel in the non-local term is exploited to obtain an explicit solution for the Fredholm transformation and its invertibility. After these first examples, two different early approaches emerged to tackle more general Fredholm transformations, a contraction approach in  \cite{BribiescaACC2014} and \cite{BribiescaArgomedo2015} that required little structure from the coefficients yet only worked for ``small'' coefficients; and a more technically involved approach tailored to self-adjoint and skew-adjoint operators in the case of a Korteweg-de Vries equation in \cite{CORON2014Fredholm} and a Kuramoto-Sivashinsky equation in \cite{CORON2015Fredholm}, respectively. In these last two papers, the existence and invertibility of the transform were related to controllability conditions.
Many recent results follow more closely a later approach, presented in \cite{Coron2016}, which linked the existence and invertibility of the Fredholm transform to a controllability condition in the case of first-order differential equations, yet required less specific conditions on the operators studied.

Following these initial results, Fredholm transformations have been successfully extended in the hyperbolic case to the finite-time stabilization of hyperbolic balance laws \cite{Coron2017}, to the output-feedback control of underactuated hyperbolic systems in \cite{Redaud2022Fredholm}, as well as internal stabilization of transport equations in \cite{Zhang2019InternalFredholm} and \cite{Zhang2022FredholmInternal}, and the stabilization of a linearized water tank system in \cite{coron2022stabilization}. In the related time-delay case, the stabilization of Integral-Delay equations has been studied in \cite{Auriol2024Fredholm}. Another related result using a Fredholm transformation concerns the linearized bilinear Schrödinger equation in \cite{coron2018rapid}. In the parabolic case, several results have been presented in the literature, including partially separable kernels in \cite{Tsubakino2014Semiseparable}, rapid stabilization of degenerate parabolic equation in \cite{gagnon2021fredholm} and a heat equation in \cite{Gagnon2022FredholmHeat}.  Finally, recent contributions~\cite{gagnon2022fredholm,hayat2024fredholm} have developed a general framework for Fredholm backstepping, extending beyond self-adjoint, skew-adjoint, or 1D systems. This framework relies on the existence of a Riesz basis of eigenvectors for the differential operator.  
 An interesting development in the approximation of both the standard backstepping and the Fredholm-transformation based backstepping using convex optimization can be found in \cite{ascencio2018backstepping}, including some technical notions on invertibility of the class of transformations obtained using these approximations.

\subsection{Sampled-data and Event-triggered systems}\label{sec-sampled}
As discussed in section \ref{sec-adaptive}, one of the strengths of the backstepping method has been its capacity to simplify convoluted structures in PDE systems and allow for an equivalent stability analysis based on a simpler target system. This advantage has been exploited not only in adaptive and parameter estimation methods but also in sampled-data and Event-triggered stabilization and estimation schemes, which generally exploit the easier structure of the target system in order to construct Lyapunov functions for the closed-loop system (or to estimate the rate of decay of some norm of the original system). 

The availability of such Lyapunov functions has allowed for sampled-data control and estimation schemes for transport PDEs (with potential non-local terms) in \cite{KARAFYLLIS2017SampledTransport}, as well as for $2 \times 2$ hyperbolic PDEs in \cite{davo2019SampledData}. Another related result, concerning quantization-related restrictions in the measurements and control inputs can be found in \cite{BekiarisLiberis2020}. In the parabolic case, \cite{KARAFYLLIS2018SampledParabolic} presents a sampled-data boundary-feedback scheme for 1-D PDEs. 

Event-based control results include $2\times 2$ linear hyperbolic systems in \cite{Espitia2018EventBased}, as well as an observer-based version in \cite{Espitia2020}. Another related system can be found in \cite{Espitia2022}, where traffic flow on cascaded roads is considered under event-triggered output-feedback. As for parabolic PDEs, existing results include \cite{ESPITIA2021EventTriggeredReactionDiffusion} for reaction-diffusion PDEs with constant coefficients, as well as \cite{Rathnayake2022ObserverEventTriggered} for an observer-based control and \cite{Rathnayake2022SampledDataEventTriggered} considering also sampled-data. For interconnected PDE-ODE systems, some results exist, such as \cite{Li2021} in the hyperbolic case, and \cite{Wang2023EventTriggeredParabolicPDEODE} in the parabolic case.

        %\begin{enumerate} \setcounter{enumi}{5}
      %  \item Time-varying PDEs \cite{Coron2021,deu2016,ker2019}%,anfinsen2020stabilization} 
      %  \& prescribed-time designs  \cite{espitia2022sensor,steeves2020prescribed,espitia2019boundary} RAFA
        %\item Fredholm transformations and non-strict feedback plants \cite{BribiescaArgomedo2015,Coron2016,redaud2022fredholm}%,coron2022stabilization,coron2018rapid,gagnon2021fredholm} 
        %FEDERICO 
      %  \item Moving boundary domains \cite{kog2020_f,kog2018}%,kog2021_b,kog2020_e,szc2022}
       % RAFA
        %\item Sampled data and Event-triggered systems %\cite{BekiarisLiberis2020,davo2019SampledData,Espitia2020,%Espitia2022,Espitia2018EventBased,Espitia2016,
        %Li2021,Wang2020}
        %FEDERICO
    %\end{enumerate}
    \section{Applications of backstepping}\label{sec-apps}

    %FEDERICO coordinates
    
    Backstepping for PDEs has shown remarkable utility in specific applications described by infinite-dimensional models, demonstrating the method's versatility and effectiveness in addressing complex engineering challenges. 
    
 %   \begin{enumerate}
 %       \item Traffic \cite{Burkhardt2021,Yu2019,Yu2021,Yu2022,Yu2023} RAFA
 %       \item Moving boundary problems: 3-D printing \cite{and2019,kog2020_a,kog2020_c}, polar ice estimation \cite{kog2020_d} RAFA
 %       \item Flow Control problems: Convection loop \cite{convloop,Vazquez2006}, Navier-Stokes \cite{coc2009,Vazquez2007}, Magnetohydrodynamics \cite{vaz2008_d,Xu2008}, Multi-Layer flows (Saint-Venant–Exner) \cite{diagne2017backstepping}, Pipe flows \cite{Aamo2016,Anfinsen2022} RAFA + ALL
 %       \item Rijke tube (thermoacustic instabilities) \cite{Andrade2018,de2020backstepping} RAFA
  %      \item Batteries \cite{kog2021,kog2022} FEDERICO + RAFA
  %      \item Drilling problems \cite{Aarsnes2019,Aarsnes2019b,Auriol2021,BreschPietri2014a,Sagert2013,Wang2020} FEDERICO + JEAN
   %     \item Others \cite{kog2020_b,wan2018,wang2022pde} FEDERICO 
   %     \item Multi-agent systems \cite{meurer2011finite,qi2019control,frihauf2010leader,qi2017wave,qi2015multi,freudenthaler2020experiments} RAFA
   % \end{enumerate}

\subsection{Traffic Control}
Freeway traffic modeling and control have been intensively
investigated over the past decades. Backstepping has proven to be an effective method for designing boundary controllers that address stop-and-go oscillations in congested traffic, often modeled by the nonlinear Aw-Rascle-Zhang (ARZ) model~\cite{Yu2019, Yu2022}. This model, consisting of second-order hyperbolic PDEs, accurately captures the complex dynamics of traffic density and velocity. Early backstepping designs focused on stabilizing traffic flow on a single freeway segment, considering both upstream and downstream control problems using ramp metering~\cite{Yu2019}. They employed backstepping transformations to map the original system to a target system with desirable stability properties, yielding full-state feedback control laws and collocated/anti-collocated observers. These designs achieved exponential stability and finite-time convergence to equilibrium, validated through simulations. Later works expanded the methodology to address the more challenging problem of simultaneous downstream and upstream stabilization on cascaded freeway segments connected by a junction, considering varying road conditions~\cite{Yu2022}. This scenario requires a more sophisticated control approach due to the interconnected nature of the system and distinct traffic scenarios on different segments. Backstepping-based output-feedback controllers were designed to suppress stop-and-go oscillations, leveraging boundary measurements of flow rate and velocity at the junction or the outlet of the segment. These designs ensured robust stabilization of the under-actuated network of two hyperbolic PDE systems, effectively mitigating congestion on cascaded roads. In \cite{guan2021,guan2023}, these works further extended the application of backstepping to the control of multi-class traffic flow, incorporating the influence of disturbances and bottlenecks, they designed an optimal observer-based output feedback controller for a heterogeneous traffic flow model with disturbances at the inlet of a considered road segment. Leveraging the backstepping approach, the proposed controller guaranteed the integral input-to-state stability of the closed-loop system and effectively addressed traffic breakdown caused by a combination of high traffic demand, bottlenecks, and driver disturbances. The controller design involved constructing an observer that estimates the states of the linearized ARZ traffic flow model using only boundary measurements. In \cite{guan2023disturb},  the backstepping method was employed to design an observer for a heterogeneous quasilinear traffic flow system, explicitly addressing the impact of disturbances at the inlet of the road section. This work demonstrated the applicability of backstepping for state observation in the presence of perturbations, ensuring accurate estimation of the traffic state despite uncertainties and nonlinearities. Backstepping has also been applied to address the control of traffic flow with multiple vehicle classes, each characterized by its own size and driver behavior~\cite{Burkhardt2021}. This requires a multi-class traffic model, typically represented by a set of coupled hyperbolic PDEs. By designing an output-feedback controller based on backstepping and an anti-collocated observer, researchers achieved finite-time damping of stop-and-go waves by measuring the velocities and densities of different vehicle classes at the inlet of the considered track section. In another work, a backstepping observer was developed and validated with real freeway data to estimate congested freeway traffic states based on the ARZ model~\cite{Yu2021}. This observer design only requires boundary measurements of flow and velocity, making it suitable for practical applications with limited sensor coverage. Finally, recent contributions used the backstepping approach to address the problem of mean-square stabilization of mixed-autonomy traffic PDE systems~\cite{zhang2024mean}. 

\subsection{Applications of Moving Boundary Domain Designs}\label{sec-appmovingboundary}
The moving boundary domain problems described in Section~\ref{sec-movingboundary} arise in various applications, including 3-D printing, estimation of polar ice melting and freezing, and crystal growth processes.  In these applications, the interface between different phases (e.g., solid and liquid, ice and water) is dynamic, and its position plays a crucial role in determining the system's overall behavior. In \cite{kog2020_c}, the authors address the problem of regulating temperature in 3-D printers based on selective laser sintering. The dynamics of the melt pool is described by a moving boundary problem, where the interface between the solid and liquid phases evolves according to the heat input from the laser. The work \cite{kog2020_a} proposes a laser sintering control design for metal additive manufacturing by PDE backstepping. The objective is to stabilize the depth of the melt pool at a desired position. The authors utilize a 1-D approximation of the Stefan problem in the vertical direction, incorporating the in-domain effect of laser power; these results are robust even when applied to a more complex two-phase Stefan model and in the presence of measurement bias. On the other hand, \cite{kog2020_d} deals with the challenging problem of estimating the thickness of polar ice sheets, which involves a moving boundary due to the melting and freezing processes. The authors employ a backstepping observer that utilizes measurements of the ice surface elevation to estimate the ice thickness, enabling the monitoring of changes in polar ice volume and its impact on sea level rise.  Another example includes \cite{cai2016nonlinear}, which addresses nonlinear stabilization through wave PDE dynamics with a moving uncontrolled boundary, with applications to oil drilling where the domain length changes with time and depends on the drill bit speed.  Finally, in \cite{and2019}, the authors present an alternative design and implementation of a backstepping controller for regulating temperature in 3-D printers based on selective laser sintering which aims to achieve rapid and accurate temperature regulation, ensuring proper layer deposition and achieving desired print quality. The authors demonstrate the effectiveness of their approach through experimental validation on a real 3-D printer.

\subsection{Flow Control}\label{sec-flowcontrol}

Flow control, encompassing the manipulation of fluid flows to achieve desired behaviors, represents a fertile ground for applications of PDE backstepping. Backstepping has enabled the design of boundary controllers and observers for a variety of flow models, including the Navier-Stokes equations, magnetohydrodynamic (MHD) systems, multi-layer flows, and pipe flows, offering new avenues for enhancing flow stability, reducing drag, and mitigating undesired flow phenomena. One of the earliest successes of PDE backstepping in flow control was the stabilization of the linearized 2-D Navier-Stokes Poiseuille flow using a closed-form feedback controller~\cite{Vazquez2007}. This work demonstrated the feasibility of designing explicit boundary control laws for a simplified flow model, paving the way for more complex scenarios such as transitions between profiles~\cite{vazquez2008control}. Subsequent efforts, such as~\cite{coc2009}, focused on motion planning and trajectory tracking for three-dimensional Poiseuille flow (which result in the well-known Orr-Sommerfeld-Squire equations). By combining backstepping with a motion planning algorithm, the authors designed a controller that drives the velocity field of the flow to a desired trajectory, illustrating the potential for achieving precise flow manipulation. The control of MHD channel flows involving electrically conducting fluids subject to magnetic fields has also benefited from backstepping designs. In~\cite{Xu2008}, the authors addressed the stabilization of the linearized 2-D MHD channel flow using backstepping boundary control. This work highlighted the applicability of backstepping to more complex flow models, incorporating the interaction between fluid dynamics and electromagnetism.  Further explorations of MHD state estimation with boundary sensors were presented in~\cite{vaz2008_d}. A common theme in these flow control applications is the utilization of spatial invariance to simplify the control design (see Section~\ref{sec-higher}; techniques like the Fourier transform can be employed to convert the PDE system into an ensemble of parameterized lower-dimensional systems. This allows for the application of backstepping to each wave number separately, simplifying the control design and analysis.

For multi-layer flows, often encountered in geophysical and environmental applications, backstepping has proven valuable in stabilizing the linearized Saint-Venant–Exner model~\cite{diagne2017backstepping}. This model, describing the flow of a fluid over a movable bed, exhibits challenging dynamics due to the coupling between the fluid and sediment transport. The control of pipe flows, crucial for a wide range of industrial processes, has also been tackled with backstepping. Notably,~\cite{Aamo2016} presents a method for leak detection, size estimation, and localization in pipe flows, while~\cite{Anfinsen2022} extends this approach to branched pipe flows.

    \subsection{Control of thermal and thermo-acoustic instabilities}
    Backstepping has proven effective in controlling instabilities arising in fluid systems subject to thermal gradients. These instabilities, often characterized by complex, nonlinear dynamics, pose unique challenges for control design due to the interplay between fluid flow and heat transfer. One such example is the Rayleigh-Bénard convective instability, a classic problem in fluid mechanics that arises when a fluid is heated from below and cooled from above. In~\cite{Vazquez2006}, a boundary feedback control law was designed to stabilize the flow in a two-dimensional thermal convection loop model. The design leveraged a combination of singular perturbation theory and backstepping for multidimensional systems. The inverse Prandtl number acted as a singular perturbation parameter, allowing for the decoupling of the system into a slow subsystem (temperature) and a fast subsystem (velocity).  Backstepping was then employed to design a state-feedback controller for the slow subsystem. An observer was designed using singular perturbations and dual backstepping methods, leading to an output feedback controller~\cite{convloop}. Similarly, thermoacoustic instabilities, often encountered in combustion systems like gas turbines and jet engines, present a challenging control problem due to the complex interplay between combustion dynamics and acoustic pressure fluctuations. The Rijke tube, an academic experiment that exhibits thermoacoustic oscillations without combustion, provides an accesible platform for exploring these instabilities. In~\cite{Andrade2018},  the problem of boundary stabilization of thermoacoustic oscillations in a linearized Rijke tube model was addressed using backstepping. Remarkably, the design yielded an exact piecewise-differentiable expression for the kernels of this transformation, allowing for the derivation of an explicit feedback control law.  In~\cite{de2020backstepping}, a Luenberger-type state observer was designed for the linearized Rijke tube model, leveraging a single boundary acoustic pressure sensor.  Experimental validation on a prototype Rijke tube demonstrated the effectiveness of this observer design in real-world scenarios.

     \subsection{Control of drilling devices}

    When extracting resources deep in the earth's subsurface, such as oil, gas, minerals, and thermal energy, it is necessary to drill long slender boreholes from the surface to the subsurface target. The drilling process involves transferring rotational and axial energy, which can result in complex drill string vibrations. These vibrations need to be controlled to improve the rate of penetration.  While earlier models used simplified mass-spring representations, more effective distributed models (namely wave equations) have been developed to capture high-frequency behaviors~\cite{Aarsnes2019b}. Such field-validated PDE models are coupled with ODEs to account for the actuator dynamics and the interaction between the drilling device and the rock. As theoretical developments emerged, it became possible to use the backstepping approach to create stabilizing controllers and state-observers for these interconnected ODE-PDE-ODE systems~\cite{BreschPietri2014a,Sagert2013}. Although earlier results overlooked the actuator dynamics, the most recent developments involve event-triggered controllers~\cite{Wang2023EventTriggeredParabolicPDEODE} and adaptive observer designs with non-linear terms to address Coulomb friction~\cite{Aarsnes2019,Auriol2021}.

    \subsection{Multi-Agent Systems}
Multi-agent systems, involving the coordination and control of multiple interacting agents, have witnessed a surge in research interest due to their potential in various applications, including robotics, autonomous vehicles, and sensor networks. While traditionally modeled using ODEs, a shift towards continuum models based on PDEs has emerged, enabling scalable control design and analysis, particularly when dealing with a large number of agents that can be approached by a continuum. In this context, backstepping has demonstrated promising capabilities for achieving coordinated motion and formation control in multi-agent systems.  The first work in this area is~\cite{frihauf2010leader}, where the authors introduced a PDE-based approach for leader-enabled deployment of multiple agents onto planar curves. By modeling the agents as a continuum, they employed backstepping to design a boundary controller that guides the agents to converge to the desired curve while maintaining a prescribed spacing. Building upon this concept,~\cite{meurer2011finite}  extended the PDE-based backstepping design to achieve finite-time multi-agent deployment. The authors considered a nonlinear PDE motion planning approach, utilizing a transport equation to model the agent dynamics; the approach ensures finite-time convergence of the agents to the desired formation, thus demonstrating the capability for rapid deployment in time-critical scenarios. Later works further expanded the application of backstepping in multi-agent control. In~\cite{qi2015multi}, the authors addressed multi-agent deployment in three dimensions using multidimensional backstepping (see Section~\ref{sec-higher}) based on a reaction-diffusion model; this approach has been later expanded in~\cite{zhang2024multi} considering more complex (radially-varying) models to extend the number of achievable configurations. Expanding on the repertoire of PDE models for multi-agent systems,~\cite{qi2017wave}  explored the use of the wave equation for time-varying formation control, leveraging the inherent wave propagation characteristics to design a controller that achieves dynamic formation reconfiguration.  More recently,~\cite{qi2019control} presented a PDE-based method for control of multi-agent systems with input delay. This work addressed the practical challenge of communication delays in multi-agent coordination, utilizing a transport equation to model the delay dynamics. By employing backstepping, a controller was designed to compensate for the delay and achieve stable formation control, ensuring robustness to communication latency. Finally,~\cite{freudenthaler2020experiments} provided experimental validation of PDE-based multi-agent formation control, combining flatness and backstepping. The authors implemented their control design on a group of mobile robots, demonstrating the real-world feasibility and effectiveness of the PDE-based backstepping approach for achieving coordinated motion in multi-agent systems.
    
\section{Issues in Practical Implementation}\label{sec-implementation}
% JEAN coordinates
% Topics covered:
%     \begin{enumerate}
%         \item Model reduction \cite{Auriol2019c,Woittennek2017} JEAN
%         \item Numerical solution of the kernel equations, including explicit solutions \cite{Ascencio2017,Vazquez2014a,wilhelmsen2022explicit} RAFA + ALL
%         %Results from Aamo's team
%         \item Neural Networks for efficient computation of backstepping \cite{Bhan2023,Shi2022} RAFA + JEAN
%     \end{enumerate}
% %\modifJA{Maybe the NN part can be merged in the two other parts. I have done it for 7.2} 

Although the discourse on the practical implementation of controllers and observers derived from PDE backstepping has been relatively sparse, the significance of implementation considerations is undeniable. Recently, the incorporation of neural networks into this area has sparked significant interest, offering a promising way to advance the application of backstepping methodologies.

\subsection{Solving the kernel equations}
%Explicit solutions + Solvers
As explained in Section~\ref{Sec_primer}, the implementation of backstepping controllers and state observers requires solving an appropriate set of kernel equations, e.g., equations~(\ref{eqn-kerneleq})--(\ref{eqn-kerneleqbc2}) for the reaction-diffusion PDE~\eqref{eqn-reactdifint}. Generally, these equations do not have a known closed-form solution and must be approximated numerically.  In addition, recent backstepping designs for coupled systems (both parabolic~\cite{Vazquez2017} and hyperbolic~\cite{Hu2016a}) exhibit discontinuous behavior which must be accurately captured with the numerical approximation. However, until very recently, there have been no publications specifically devoted to this topic. Only a few sections and appendices spread through the backstepping literature gave some clues about numerical algorithms, which include finite difference approximations of the kernel equations~\cite{krs2008,ker2019, deutscher2021backstepping,anfinsen2019adaptive}, the use of symbolical successive approximation series~\cite{convloop}, or the numerical solution of the integral version of the kernel equations~\cite{jadachowski2012efficient,Auriol2022}. A more recent approach uses power series to obtain solutions for the backstepping kernel equations when coefficients are analytic. This idea was first seen in~\cite{ascencio2018backstepping} (without much analysis of the convergence of the method itself, but rather using ideas from convex optimization to approach the kernels as best as possible while obtaining stability) and the convergence of the series was first rigorously analyzed in~\cite{R.2022ball} for a multidimensional case that presents singularities at the origin and thus is not amenable to other methods. Recent works~\cite{vazquez2023power,lin2024towards} leverage complex analysis to show that the method can be applied to most kernel equations found in practice, providing respectively Mathematica and Matlab code to the reader that can be adapted to most particular cases.

In addition, during the development of the backstepping method, it was discovered that in certain simple cases, the kernel solutions could be expressed in closed form~\cite{Smyshlyaev2004,smy2005b}. Having an explicit expression of the kernels offers several advantages. For example, it provides a better understanding of the control law's structure and its dependence on the different parameters of the system. Furthermore, having an analytical expression makes implementation simpler and more precise. Most importantly, when the parameters of the PDE plant are unknown, it is crucial to compute the controller gains in a timely manner. This becomes possible when these gains are available as explicit functions of the plant parameters.
\newline
For the parabolic reaction-diffusion equation~\eqref{eqn-reactdifint} with constant reaction term, it has been demonstrated in~\cite{krs2008} that the solutions to the kernel equations~(\ref{eqn-kerneleq})--(\ref{eqn-kerneleqbc2}) can be explicitly expressed using the first-order modified Bessel function of the first kind. Interestingly, for $2\times 2$ hyperbolic systems with constant coefficients (i.e. system~\eqref{Sec_4_1_v_pde_orig_simp} with $n=m=1$ and without ODEs), the kernels can also be expressed using modified Bessel functions of the first kind~\cite{Vazquez2014a}. An alternative expression based on the generalized Marcum Q-function of first order is also available~\cite{marcum1950table}. For such $2\times 2$ hyperbolic systems,  it is possible to use similar techniques to obtain an explicit expression of the inverse kernels~\cite{saba2019stability} (see Section~\ref{sec-inv} for the definition of the inverse kernels). This explicit expression is useful for deriving simple and easy-to-compute stability conditions for these hyperbolic systems, as discussed in the same paper.  Finally, extensions to $n+m$ heterodirectional hyperbolic systems with exponential in-domain coupling terms have been proposed in~\cite{wilhelmsen2022explicit}.

% For these systems, the class of constant-coefficient parameters gives a very wide range of plants, with 7 distinct parameters that can have any value, with only the speeds of propagation being restricted to be positive.

    \subsection{Integration and model reduction}

    The success of classical control design algorithms, such as PID controllers, in industrial engineering applications has been attributed to their simplicity and low computational burden. Implementing backstepping output-feedback controllers requires expertise and computing power. Real-time state estimation of distributed states using rapidly converging observers, particularly those based on backstepping, can be computationally expensive and prohibitive in some cases. This burden is further emphasized when dealing with networks of PDEs systems. To handle the computational demands of backstepping output-feedback controllers, it is essential to approximate them, usually through finite-dimensional systems. However, the approximation method should ensure satisfactory closed-loop properties, particularly closed-loop stability. Various late-lumping approximation methods aiming to reduce the computational burden of the proposed backstepping controllers while maintaining comparable closed-loop performance have been suggested in~\cite{Woittennek2017, ecklebe2017approximation} or~\cite{Auriol2019c}. 
    
    In~\cite{Woittennek2017}, a method has been proposed for computing the bounded part of the backstepping control operator. The proposed approach relies on a finite-dimensional representation of the state, which helps in efficiently computing the feedback law. By choosing suitable finite-dimensional approximation spaces, it is possible to simplify the computation of backstepping transformations. However, there are no guarantees of convergence. Additionally, such an approach does not cover the implementation of observer schemes. In \cite{Auriol2019c}, general assumptions are provided regarding approximation schemes for backstepping state-feedback controllers. %The paper demonstrates stability for three simple test cases by combining a Lyapunov analysis with the backstepping method as an analysis tool.
    In this paper, the approximated control law $U_{ap}$ is rewritten as $U_{ap}=U+(U_{ap}-U)$ where $U$ is the nominal non-approximated backstepping controller. The difference between the approximated control law and the non-approximated one $(U_{ap}-U)$ is considered as a disturbance term, whose effect on the closed-loop system can be estimated using Input-to-State Stability estimates~\cite{Auriol2019c}. Closed-loop stability is then guaranteed if this disturbance term is small enough. The Lyapunov analysis initiated in~\cite{Auriol2019c} on simple examples could be extended to guarantee the convergence of approximation schemes under generic conditions as long as an ISS-Lyapunov functional is available.

\subsection{Machine Learning and PDE backstepping}

In a recent study~\cite{yu2021reinforcement}, backstepping controllers were compared to reinforcement learning controllers for addressing stop-and-go traffic congestion on a freeway segment. The study found that while reinforcement learning controllers can achieve similar performance to backstepping controllers with less computation time, they typically require one thousand episodes of iterative training on a simulation model. 
Conducting this training in a collision-free manner within real traffic is impractical for PDE-based models, although some successful implementations exist for simpler ODE-based traffic models \cite{lee2025traffic}. However, even in these cases, convergence during training cannot be guaranteed.  As a result, reinforcement learning controllers are not a straightforward or entirely secure replacement for model-based control in real traffic systems. However, it is possible to successfully combine backstepping with machine learning. A recent development~\cite{Bhan2023} introduced a promising approach, involving a direct deep-neural network approximation of the backstepping controller for a transport equation with re-circulation. This model-based machine-learning backstepping controller offers established convergence guarantees using the DeepONet general approximation theorem, along with a Lyapunov analysis. 
%This approach is currently being extended to a wide range of systems, requiring the existence of ISS-Lyapunov functionals for the target system. 
Current progress in this area demonstrate the enormous potential of neural operators to effectively learn backstepping kernels (or alternatively, gains, as explained in~\cite{vazquez2024gain}) while maintaining stability guarantees. This is reflected in the latest works, which include designing neural operator-based observers and controllers for parabolic PDEs \cite{krstic2023neural,redaud2024physics}, extending the approach to more complex hyperbolic PDE systems \cite{wang2023backstepping}, developing adaptive neural operator-based control designs \cite{lamarque2024adaptive}, and even implementing gain scheduling techniques using neural operators \cite{lamarque2024gain}. 

\section{Conclusions: open problems and future lines of work}\label{sec-conclusions}
%FEDERICO coordinates

The paper culminates with a reflective discourse on unresolved challenges and potential directions for future research in the realm of PDE backstepping, offering insights into the frontier areas of this evolving field. Some of these topics might not have many available references yet, but they still represent interesting areas of active research.

%\begin{enumerate}
%\item Underactuated systems~\cite{auriol2020robust} FEDERICO + JEAN 
%\item Backstepping controllers with in-domain actuation~\cite{Zhang2019InternalFredholm} FEDERICO
%\item Backstepping and stochastic systems~\cite{Auriol2023b} JEAN
%\item Towards parametrizable target systems~\cite{ramirez2017backstepping} JEAN
 
%\item Use of neural operators in gain scheduling and adaptive control~\cite{lamarque2024gain} RAFA -> I THINK THIS IS ALREADY COVERED BUT PERHAPS WE WANT TO MOVE IT HERE
%\item More general multidimensional domains RAFA
%\item Ensembles of PDEs FEDERICO -> the multidimensional domains can also be tied to ensembles 
%\end{enumerate}

\subsection{Underactuated systems}

As emphasized throughout this survey, most backstepping controllers are designed assuming that the state at one of the domain boundaries is fully actuated. This is motivated by underlying controllability results that guarantee the possibility of stabilizing controllers to be designed without any constraints on the system parameters~\cite{li2010strong}. On the other hand, underactuated systems are those where the dimension of the control input is smaller than the dimension of the boundary state. For example, this configuration occurs when considering system~\eqref{Sec_4_1_v_pde_orig_simp} without ODEs and with the boundary condition $v(t,1)=Ru(t,1)+BU(t)$ (where the rank of $B$ is smaller than the dimension of $v$). Currently, there are limited results in the literature for designing backstepping controllers for underactuated systems, and such systems are often analyzed as interconnected systems. Therefore, most existing designs require structural conditions such as cascade structures~\cite{auriol2024contributions,Chen2023b,deutscher2023backstepping}. One attempt has been made for the backstepping stabilization of a $1+2$ hyperbolic system with one control input in~\cite{Auriol2020c}. This current lack of literature presents a compelling opportunity for promising investigations in the near future.

\subsection{In-domain actuation}
Despite the many developments concerning PDE backstepping, an area that has resisted generic results has been the extension to in-domain actuation (or sensing). As a result of the original Volterra operator structure chosen, backstepping has been developed as a boundary control and observation method.
As a result, most of the references given in this section are conference papers. 

Some early results concerning potential extensions to in-domain actuation are \cite{Tsubakino2012InDomain} and \cite{WANG2013indomain} in the parabolic case. These two results require very specific conditions on the shape of the actuation in the former and specific domain lengths in the latter (similar to rational dependence). An extension to observer design with measurements of weighted averages can also be found in \cite{TSUBAKINO2015observerweighted}. For hyperbolic systems, early results can be found in \cite{Zhu2022distributed} and \cite{REDAUD2023indomain}. In particular, the last result links with the Fredholm transformation results in Section \ref{sec-fredholm}. Seeing the recent advances in Fredholm transformation development, it is possible that more generic results in this area could be obtained soon.

\subsection{Backstepping and stochastic systems}

Disturbances and uncertainties are prevalent in numerous systems, leading to performance issues. As demonstrated throughout this survey, backstepping-based control approaches have been adjusted to address these robustness challenges for PDE systems. However, many of these techniques focus on worst-case scenarios, leading to overly conservative controller designs. External disturbances often occur randomly, motivating the need for robust control of switching or stochastic PDEs. Despite this, the backstepping approach remains a powerful tool for control design and robustness analysis.

% Disturbances and uncertainties are prevalent in numerous systems, often resulting in adverse effects on the performance in closed-loop. As demonstrated throughout this survey, various backstepping-based control approaches have been devised to address challenges such as disturbance attenuation, disturbance rejection, and robust control for PDE systems. Most of the presented contributions focus on the worst-case scenario, potentially leading to overly conservative controller designs. Moreover, when considering applications the external disturbances most often appear in a random way. This motivates the analysis and robust control of switching or stochastic PDEs. In this context, the backstepping approach still appears as a powerful tool for control design or robustness analysis.

In~\cite{kong2022prediction}, the authors considered a time-delay system in which the delay $\tau(t)$ was treated as a stochastic process following a Markov model with a finite number of possible values. They developed a prediction-based controller specifically tailored for a constant reference delay value. A backstepping analysis showed that mean-square closed-loop stability could be ensured as long as the stochastic delays were sufficiently close to the nominal delay on average. An extension of this approach was proposed in~\cite{zhang2023adaptive} to design a robust adaptive controller in the case of uncertain dynamics. The same methodology was applied in~\cite{Auriol2023b} to prove the mean-square exponential stability of two coupled hyperbolic equations (as in ~\eqref{Sec_4_1_v_pde_orig_simp} with $n=m=1$ and without the ODE) with random parameters. Recent work in \cite{guan2024robustness} addresses the robustness of reaction-diffusion PDE predictor-feedback to stochastic delay perturbations, analyzing how random phenomena affect actuator compensation and establishing mean-square exponential stability under stochastic disturbances.

% The problem of prediction-based control for dynamical systems subject to stochastic input delays was considered for the first time in~\cite{kong2022prediction}. The objective is to stabilize the system 
% \begin{align}
%     \dot{X}(t)=AX(t)+BU(t-\tau(t)), \label{Sec_4_Delayed_ODE_stocha}
% \end{align}
% where, compared to equation~\eqref{Sec_4_Delayed_ODE}, $\tau(t)$ is now stochastic and is defined as a Markov process with a finite number of possible values (random switching system). The solution proposed in~\cite{kong2022prediction} consists of designing a prediction-based controller defined for a reference constant delay value. Mean-square closed-loop stability can be guaranteed if the stochastic delays are close enough to the nominal one in expectancy. The associated robustness analysis relies on a PDE formulation of the time-delay system and on backstepping transformations adjusted from~\cite{Krstic2008c}. An extension of this approach is proposed in~\cite{zhang2023adaptive} to design a robust adaptive controller in the case of uncertain dynamics. Interestingly, the same methodology is applied in~\cite{Auriol2023b} to prove the mean-square exponential stability of two coupled hyperbolic equations (as in equation~\eqref{Sec_4_1_v_pde_orig_simp} with $n=m=1$ and without the ODE) with stochastic parameters.

Mean-square exponential stabilization of a stochastic reaction-diffusion equation has been considered in~\cite{wu2022disturbance}. The backstepping approach was used to design a nominal stabilizing controller that guarantees the mean-square exponential stability. This result suggests most designs presented in this survey to design stabilizing controllers for interconnected ODE-PDE systems could be extended to mean-square stabilizing controllers for the stochastic counterparts. However, stabilizing in expectation should not be the only objective, and one must make sure the variance remains bounded to ensure reliability and closed-loop robustness. In this context, the methodology developed in~\cite{velho2024stabilization} for a hyperbolic PDE interconnected with a stochastic differential equation seems promising. Designing backstepping-based risk-averse methods that scale with many different risk measures is a challenging task for future research.

\subsection{Towards parametrizable target systems}

One of the main challenges of the backstepping method is finding an appropriate target system. This target system should be simple enough to design a stabilizing control law while also guaranteeing the existence of a transformation mapping the original system to this target system. The choice of target system directly affects the closed-loop system's performance, and the broader issue of identifying reachable target systems remains unresolved. For heterodirectional hyperbolic equations, most contributions typically choose finite-time stable target systems~\cite{Coron2021}, thereby shadowing the robustness properties of the corresponding closed-loop systems~\cite{logemann1996conditions,michiels2009strong}. However, there might be value in exploring more complex target systems. For instance, in~\cite{auriol2020robust}, tuning parameters are introduced, guaranteeing potential trade-offs between delay-robustness and convergence rate. However, these tuning parameters are limited in their action as they only affect the boundary conditions of the system. For the simple parabolic system presented in Section~\ref{Sec_primer} the coefficient $c$ introduced in the target system~\eqref{eqn-target} can also be seen as a tunable parameter that modifies the closed-loop decay rate.

For the $n+m$ hyperbolic systems presented in Section~\ref{Sec_hyp_systems}, general target systems could be obtained by preserving dissipative in-domain couplings. For instance, the target systems associated with the system~\eqref{Sec_4_1_v_pde_orig_simp} could be parametrized as
%This would require precise knowledge of their influence in terms of stability. For instance, for $n+m$ heterodirectional systems (as described by equation~\eqref{Sec_4_1_v_pde_orig_simp}), possible target systems could read
\begin{equation}\label{Syst_param}
\left\{
\begin{array}{l}
    \gamma_t+\Lambda(x)\gamma_x=\bar \Sigma(x)\gamma,\\
    \alpha(0)=Q\beta(0),~\beta(1)=\bar R \alpha(1),
\end{array} 
\right.
\end{equation}
with an appropriate choice of $\bar \Sigma$ and $\bar R$ to guarantee closed-loop stability. It was shown in~\cite{redaud2024domain} that the target system~\eqref{Syst_param} was always reachable from the original system~\eqref{Sec_4_1_v_pde_orig_simp}.
On the other hand, the Port-Hamiltonian approach (PHS)~\cite{le2005dirac,ramirez2017backstepping} can provide a structured methodology to emphasize and leverage the intrinsic physical properties of the original system (such as passivity, dissipativity, and reversibility) therefore guiding appropriate selections of $\bar \Sigma$ and $\bar R$ in the target system~\eqref{Syst_param}. %These properties can be used to define well-posed, exponentially stable target system candidates. In~\cite{ramirez2017backstepping}, the backstepping approach is, for instance, used to design boundary controllers for a certain class of linear PHS. This is achieved by expressing the target system as a PHS with a tunable in-domain dissipation term. 
All in all, incorporating degrees of freedom into the design can result in modular controllers that allow for various trade-offs between performance specifications, such as robustness margins, sensitivity, and convergence rate. Although we may have an idea of the qualitative impact of tuning parameters, quantifying the performance of the resulting controllers in terms of performance specifications remains an area of ongoing research.
%In this context, the Port-Hamiltonian approach (PHS)~\cite{le2005dirac} could be a path to follow as it corresponds to a multi-physical and modular energy-based representation that considers the system's natural physical properties (e.g., passivity, dissipativity, reversibility). It offers a structured methodology to highlight and capitalize on the intrinsic physical properties of the system that could then be advantageously used to define well-posed, exponentially stable target system candidates. In~\cite{ramirez2017backstepping}, the backstepping approach is used to develop boundary controllers for a class of linear Port-Hamiltonian Systems with constant parameters. The corresponding target system is a PHS with a tunable in-domain dissipation term. Recently, the PHS framework was successfully applied in \cite{redaud2024domain} to introduce degrees of freedom in the control design to obtain a class of easily parametrizable, exponentially stable closed-loop systems with a specified energy decay. The associated controllers are obtained using the backstepping methodology, combining classical Volterra transformations with innovative time-affine transforms. Introducing degrees of freedom in the design leads to modular controllers that can guarantee potential trade-offs between performance specifications (such as robustness margins, sensitivity, or convergence rate). If the qualitative effect of the tuning parameters can be understood,  quantifying the performance of the resulting controllers with respect to performance specifications is still an open-field of research.

\subsection{Multi-dimensional Systems}

While Section \ref{sec-higher} discusses existing results for higher-dimensional domains, most backstepping applications remain fundamentally one-dimensional. Current multi-dimensional results are primarily limited to special geometries---such as balls, disks, and boxes---where symmetry or separability allows reduction to one or more 1D systems. Even in these cases, the resulting kernel equations become significantly more complex, often requiring sophisticated tools from complex analysis and special functions theory.
A general theory for multi-dimensional backstepping that does not rely on dimensional reduction or special geometric properties would represent a significant breakthrough. Such a theory would need to address several challenges: the formulation of appropriate target systems in multiple dimensions, the well-posedness of the resulting multi-dimensional kernel equations, and the development of practical methods for computing the resulting kernel functions. Success in this direction could dramatically expand the applicability of backstepping to important physical systems that are inherently multi-dimensional, such as fluid flows, electromagnetic fields, and complex thermal systems. Addressing some of these challenges, reference~\cite{vazquez2025backstepping} explores these issues and investigates a potential direction using domain extension techniques, particularly for handling irregular domains.

 \bibliographystyle{unsrt}
 \bibliography{Cleaned_refs}

\begin{thebibliography}{100}

\bibitem{vaz2008}
R.~Vazquez and M.~Krstic.
\newblock {\em {Control of Turbulent and Magnetohydrodynamic Channel Flows}}.
\newblock Systems \& Control: Foundations \& Applications. Birkhäuser, 2008.

\bibitem{kog2020_e}
S.~Koga and M.~Krstic.
\newblock {\em Materials Phase Change {PDE} Control \& Estimation}.
\newblock Systems \& Control: Foundations \& Applications. Birkhäuser, 2020.

\bibitem{BekiarisLiberis2013}
N.~Bekiaris-Liberis and M.~Krstic.
\newblock {\em Nonlinear Control Under Nonconstant Delays}.
\newblock {S}{I}{A}{M}, 2013.

\bibitem{krstic5}
M.~Krstic.
\newblock {\em Delay Compensation for {N}onlinear, {A}daptive, and {P}{D}{E}
  Systems}.
\newblock Systems \& Control: Foundations \& Applications. Birkhäuser, 2009.

\bibitem{meurer2}
T.~Meurer.
\newblock {\em Control of Higher-Dimensional {PDEs}: Flatness and Backstepping
  Designs}.
\newblock Communications and Control Engineering. Springer, 2012.

\bibitem{adaptive}
A.~Smyshlyaev and M.~Krstic.
\newblock {\em Adaptive Control of Parabolic {PDE}s}.
\newblock Princeton University Press, 2010.

\bibitem{wang2022pde}
J.~Wang and M.~Krstic.
\newblock {\em {PDE} Control of String-Actuated Motion}, volume~73.
\newblock Princeton University Press, 2022.

\bibitem{Yu2023}
H.~Yu and M.~Krsti{\'c}.
\newblock {\em Traffic Congestion Control by {PDE} Backstepping}.
\newblock Systems \& Control: Foundations \& Applications. Birkhäuser, 2022.

\bibitem{anfinsen2019adaptive}
H.~Anfinsen and O.~M. Aamo.
\newblock {\em Adaptive control of hyperbolic {PDE}s}.
\newblock Communications and Control Engineering. Springer, 2019.

\bibitem{krs2008}
M.~Krstic and A.~Smyshlyaev.
\newblock {\em Boundary control of {P}{D}{E}s: {A} course on backstepping
  designs}, volume~16.
\newblock {S}{I}{A}{M}, 2008.

\bibitem{colton}
D.~Colton.
\newblock The solution of initial-boundary value problems for parabolic
  equations by the method of integral operators.
\newblock {\em Journal of Differential Equations}, 26:181--190, 1977.

\bibitem{seid1984}
T.~I. Seidman.
\newblock Two results on exact boundary control of parabolic equations.
\newblock {\em Applied Mathematics and Optimization}, 11:145--152, 1984.

\bibitem{tsinias1989sufficient}
J.~Tsinias.
\newblock Sufficient {L}yapunov-like conditions for stabilization.
\newblock {\em Mathematics of control, Signals and Systems}, 2:343--357, 1989.

\bibitem{koditschek1987adaptive}
D.~E. Koditschek.
\newblock Adaptive techniques for mechanical systems.
\newblock In {\em Fifth Yale Workshop on Applications of Adaptive Systems
  Theory}, pages 259--265, 1987.

\bibitem{byrnes1989new}
C.~I. Byrnes and A.~Isidori.
\newblock New results and examples in nonlinear feedback stabilization.
\newblock {\em Systems \& Control Letters}, 12(5):437--442, 1989.

\bibitem{sontag1989further}
E.~D. Sontag and H.~J. Sussmann.
\newblock Further comments on the stabilizability of the angular velocity of a
  rigid body.
\newblock {\em Systems \& Control Letters}, 12(3):213--217, 1989.

\bibitem{KKK}
M.~Krstic, I.~Kanellakopoulos, and P.~V. Kokotovic.
\newblock {\em Nonlinear and Adaptive Control Design}.
\newblock Wiley, New York, 1995.

\bibitem{coron1998stabilization}
J.-M. Coron and B.~d'Andrea Novel.
\newblock Stabilization of a rotating body beam without damping.
\newblock {\em IEEE Transactions on Automatic Control}, 43(5):608--618, 1998.

\bibitem{liu2000backstepping}
W.~Liu and M.~Krstic.
\newblock Backstepping boundary control of {B}urgers’ equation with actuator
  dynamics.
\newblock {\em Systems \& Control Letters}, 41(4):291--303, 2000.

\bibitem{Boskovic2001}
D.~M. Bo{\v{s}}kovi{\'c}, M.~Krstic, and W.J. Liu.
\newblock Boundary control of an unstable heat equation via measurement of
  domain-averaged temperature.
\newblock {\em IEEE Transactions on Automatic Control}, 46(12):2022--2028,
  2001.

\bibitem{Boskovic2003}
D.~M. Bo{\v{s}}kovi{\'c}, A.~Balogh, and M.~Krsti{\'c}.
\newblock Backstepping in infinite dimension for a class of parabolic
  distributed parameter systems.
\newblock {\em Mathematics of Control, Signals, and Systems}, 16(1):44--75,
  2003.

\bibitem{Balogh2002}
A.~Balogh and M.~Krstic.
\newblock Infinite dimensional backstepping-style feedback transformations for
  a heat equation with an arbitrary level of instability.
\newblock {\em European Journal of Control}, 8(2):165--175, 2002.

\bibitem{Smyshlyaev2004}
A.~Smyshlyaev and M.~Krstic.
\newblock Closed-form boundary state feedbacks for a class of 1-{D} partial
  integro-differential equations.
\newblock {\em IEEE Transactions on Automatic Control}, 49(12):2185--2202,
  2004.

\bibitem{liu2003boundary}
W.~Liu.
\newblock Boundary feedback stabilization of an unstable heat equation.
\newblock {\em SIAM Journal on Control and Optimization}, 42(3):1033--1043,
  2003.

\bibitem{morris2010control}
K.~Morris and W~Levine.
\newblock Control of systems governed by partial differential equations.
\newblock {\em The control theory handbook}, 2010.

\bibitem{Smyshlyaev2005}
A.~Smyshlyaev and M.~Krstic.
\newblock Backstepping observers for a class of parabolic {P}{D}{E}s.
\newblock {\em Systems \& Control Letters}, 54(7):613--625, 2005.

\bibitem{Lasiecka2000}
I.~Lasiecka and R.~Triggiani.
\newblock {\em Control Theory for Partial Differential Equations: Continuous
  and Approximation Theories}.
\newblock Cambridge University Press, 2000.

\bibitem{Curtain2012}
R.~F. Curtain and H.~Zwart.
\newblock {\em Introduction to Infinite-Dimensional Systems Theory: A
  State-Space Approach}, volume~71 of {\em Texts in Applied Mathematics}.
\newblock Springer, 2020.

\bibitem{krener2021optimal}
A.~Krener.
\newblock Optimal boundary control of a nonlinear reaction diffusion equation
  via completing the square and {A}l’brekht’s method.
\newblock {\em IEEE Transactions on Automatic Control}, 67(12):6698--6709,
  2021.

\bibitem{Bensoussan2007}
A.~Bensoussan, G.~Da Prato, M.~C. Delfour, and S.~K. Mitter.
\newblock {\em Representation and Control of Infinite Dimensional Systems}.
\newblock Systems \& Control: Foundations \& Applications. Birkhäuser, 2007.

\bibitem{afshar2025}
S.~Afshar, F.~Germ, and K.~Morris.
\newblock Extended {K}alman filter-based observer design for semilinear
  infinite-dimensional systems.
\newblock {\em IEEE Transactions on Automatic Control}, 69(6):3631--3646, 2024.

\bibitem{lagnese2012domain}
J.~E. Lagnese and G.~Leugering.
\newblock {\em Domain decomposition methods in optimal control of partial
  differential equations}, volume 148.
\newblock Birkh{\"a}user, 2012.

\bibitem{baker2000finite}
J.~Baker and P.~D. Christofides.
\newblock Finite-dimensional approximation and control of non-linear parabolic
  {PDE} systems.
\newblock {\em International Journal of Control}, 73(5):439--456, 2000.

\bibitem{lhachemi2022finite}
H.~Lhachemi and C.~Prieur.
\newblock Finite-dimensional observer-based boundary stabilization of
  reaction--diffusion equations with either a {D}irichlet or {N}eumann boundary
  measurement.
\newblock {\em Automatica}, 135:109955, 2022.

\bibitem{mazenc2011strict}
F.~Mazenc and C.~Prieur.
\newblock Strict {L}yapunov functions for semilinear parabolic partial
  differential equations.
\newblock {\em Mathematical Control and Related Fields}, 1(2):231--250, 2011.

\bibitem{Coron2007}
J.-M. Coron.
\newblock {\em Control and Nonlinearity}.
\newblock American Mathematical Society, 2007.

\bibitem{Rudolph2003}
J.~Rudolph.
\newblock {\em Flatness Based Control of Distributed Parameter Systems}.
\newblock Shaker Verlag, 2003.

\bibitem{bohm1998model}
M.~B{\"o}hm, M.~A. Demetriou, S.~Reich, and I.~G. Rosen.
\newblock Model reference adaptive control of distributed parameter systems.
\newblock {\em SIAM Journal on Control and Optimization}, 36(1):33--81, 1998.

\bibitem{Jacob-2012-book}
B.~Jacob and H.~J. Zwart.
\newblock {\em Linear Port-{H}amiltonian Systems on Infinite-dimensional
  Spaces}.
\newblock Operator Theory: Advances and Applications. Springer, 2012.

\bibitem{vanderSchaft_2002_DPS}
A.~J. van~der Schaft and B.~M. Maschke.
\newblock {H}amiltonian formulation of distributed-parameter systems with
  boundary energy flow.
\newblock {\em Journal of Geometry and Physics}, 42(1):166--194, 2002.

\bibitem{duindam2009modeling}
V.~Duindam, A.~Macchelli, S.~Stramigioli, and H.~Bruyninckx.
\newblock {\em Modeling and Control of Complex Physical Systems: The
  Port-{H}amiltonian Approach}.
\newblock Springer, 2009.

\bibitem{wu2017port}
Z.-G. Wu and S.~M Swei.
\newblock Port-{H}amiltonian modeling and control of beam vibration.
\newblock {\em Mechanical Systems and Signal Processing}, 90:322--336, 2017.

\bibitem{Macchelli-2017-TAC_synthPHS}
A.~Macchelli, Y.~Le~Gorrec, H.~Ramirez, and H.~Zwart.
\newblock On the synthesis of boundary control laws for distributed
  port-{H}amiltonian systems.
\newblock {\em IEEE Transactions on Automatic Control}, 62(4):1700--1713, 2017.

\bibitem{Evans2010}
L.~C. Evans.
\newblock {\em Partial Differential Equations}, volume~19 of {\em Graduate
  Studies in Mathematics}.
\newblock American Mathematical Society, Providence, Rhode Island, 2 edition,
  2010.

\bibitem{yoshida1960lectures}
K.~Yoshida.
\newblock {\em Lectures on differential and integral equations}, volume~10.
\newblock Interscience Publishers, 1960.

\bibitem{moura2013batteries}
S.~J. Moura, N.~A. Chaturvedi, and M.~Krstic.
\newblock Adaptive partial differential equation observer for battery
  state-of-charge/state-of-health estimation via an electrochemical model.
\newblock {\em Journal of Dynamic Systems, Measurement, and Control},
  136(1):011015, 2013.

\bibitem{moura2017batteries}
S.~J. Moura, F.~Bribiesca-Argomedo, R.~Klein, A.~Mirtabatabaei, and M.~Krstic.
\newblock Battery state estimation for a single particle model with electrolyte
  dynamics.
\newblock {\em IEEE Transactions on Control Systems Technology},
  25(2):453--468, 2017.

\bibitem{qi2015multi}
J.~Qi, R.~Vazquez, and M.~Krstic.
\newblock Multi-agent deployment in 3-{D} via {PDE} control.
\newblock {\em IEEE Transactions on Automatic Control}, 60(4):891--906, 2015.

\bibitem{Krstic2008c}
M.~Krstic and A.~Smyshlyaev.
\newblock Backstepping boundary control for first-order hyperbolic {PDE}s and
  application to systems with actuator and sensor delays.
\newblock {\em Systems \& Control Letters}, 57(9):750 -- 758, 2008.

\bibitem{Smyshlyaev2009}
A.~Smyshlyaev and M.~Krstic.
\newblock Boundary control of an anti-stable wave equation with anti-damping on
  the uncontrolled boundary.
\newblock {\em Systems \& Control Letters}, 58(8):617--623, 2009.

\bibitem{Krstic2008a}
M.~Krstic, B.-Z. Guo, A.~Balogh, and A.~Smyshlyaev.
\newblock Output-feedback stabilization of an unstable wave equation.
\newblock {\em Automatica}, 44(1):63--74, 2008.

\bibitem{Smyshlyaev2010}
A.~Smyshlyaev, E.~Cerpa, and M.~Krstic.
\newblock Boundary stabilization of a 1-{D} wave equation with in-domain
  antidamping.
\newblock {\em {SIAM} Journal on Control and Optimization}, 48(6):4014--4031,
  2010.

\bibitem{Zhou2012}
Z.~Zhou and S.~Tang.
\newblock Boundary stabilization of a coupled wave-{ODE} system with internal
  anti-damping.
\newblock {\em International Journal of Control}, 85(11):1683--1693, 2012.

\bibitem{wu2020output}
X.-H. Wu, H.~Feng, and B.-Z. Guo.
\newblock Output feedback stabilization for 1-{D} wave equation with variable
  coefficients and non-collocated observation.
\newblock {\em Systems \& Control Letters}, 145:104780, 2020.

\bibitem{Krstic2006}
M.~Krstic, A.~A. Siranosian, and A.~Smyshlyaev.
\newblock Backstepping boundary controllers and observers for the slender
  {T}imoshenko beam: Part {I}-{D}esign.
\newblock In {\em 2006 American Control Conference}, pages 2412--2417. IEEE,
  2006.

\bibitem{krstic2008beam}
M.~Krstic, B.-Z. Guo, A.~Balogh, and A.~Smyshlyaev.
\newblock Control of a tip-force destabilized shear beam by observer-based
  boundary feedback.
\newblock {\em SIAM Journal on Control and Optimization}, 47(2):553--574, 2008.

\bibitem{Aamo2005}
O.~M. Aamo, A.~Smyshlyaev, and M.~Krstic.
\newblock Boundary control of the linearized {G}inzburg-{L}andau model of
  vortex shedding.
\newblock {\em SIAM Journal on Control and Optimization}, 43(6):1953--1971,
  2005.

\bibitem{krstic2011boundary}
M.~Krstic, B.-Z. Guo, and A.~Smyshlyaev.
\newblock Boundary controllers and observers for the linearized
  {S}chr{\"o}dinger equation.
\newblock {\em {SIAM} Journal on Control and Optimization}, 49(4):1479--1497,
  2011.

\bibitem{xiang2019null}
S.~Xiang.
\newblock Null controllability of a linearized {K}orteweg--de {V}ries equation
  by backstepping approach.
\newblock {\em SIAM Journal on Control and Optimization}, 57(2):1493--1515,
  2019.

\bibitem{bastin2016stability}
G.~Bastin and J.-M. Coron.
\newblock {\em Stability and boundary stabilization of 1-{D} hyperbolic
  systems}.
\newblock Progress in Nonlinear Differential Equations and Their Applications.
  Birkhäuser, 2016.

\bibitem{Coron2013}
J.-M. Coron, R.~Vazquez, M.~Krstic, and G.~Bastin.
\newblock Local exponential ${H}^2$ stabilization of a 2$\times$ 2 quasilinear
  hyperbolic system using backstepping.
\newblock {\em {SIAM} Journal on Control and Optimization}, 51(3):2005--2035,
  2013.

\bibitem{vazquez2011backstepping}
R.~Vazquez, M.~Krstic, and J.-M. Coron.
\newblock Backstepping boundary stabilization and state estimation of a
  2$\times$ 2 linear hyperbolic system.
\newblock In {\em 2011 50th IEEE conference on decision and control and
  european control conference}, pages 4937--4942. IEEE, 2011.

\bibitem{coron2008dissipative}
J.-M. Coron, G.~Bastin, and B.~d'Andr{\'e}a Novel.
\newblock Dissipative boundary conditions for one-dimensional nonlinear
  hyperbolic systems.
\newblock {\em {SIAM} Journal on Control and Optimization}, 47(3):1460--1498,
  2008.

\bibitem{li2010strong}
T.~Li and B.~Rao.
\newblock Strong (weak) exact controllability and strong (weak) exact
  observability for quasilinear hyperbolic systems.
\newblock {\em Chinese Annals of Mathematics, Series B}, 31(5):723--742, 2010.

\bibitem{DiMeglio2013a}
F.~Di~Meglio, R.~Vazquez, and M.~Krstic.
\newblock Stabilization of a system of $ n+ 1$ coupled first-order hyperbolic
  linear {PDE}s with a single boundary input.
\newblock {\em IEEE Transactions on Automatic Control}, 58(12):3097--3111,
  2013.

\bibitem{Hu2016a}
L.~Hu, F.~Di~Meglio, R.~Vazquez, and M.~Krstic.
\newblock Control of homodirectional and general heterodirectional linear
  coupled hyperbolic {P}{D}{E}s.
\newblock {\em IEEE Transactions on Automatic Control}, 61(11):3301--3314,
  2016.

\bibitem{Auriol2016}
J.~Auriol and F.~Di~Meglio.
\newblock Minimum time control of heterodirectional linear coupled hyperbolic
  {PDE}s.
\newblock {\em Automatica}, 71:300--307, 2016.

\bibitem{Coron2017}
J.-M. Coron, L.~Hu, and G.~Olive.
\newblock Finite-time boundary stabilization of general linear hyperbolic
  balance laws via {F}redholm backstepping transformation.
\newblock {\em Automatica}, 84:95--100, 2017.

\bibitem{coron2019optimal}
J.-M. Coron and H.-M. Nguyen.
\newblock Optimal time for the controllability of linear hyperbolic systems in
  one-dimensional space.
\newblock {\em {SIAM} Journal on Control and Optimization}, 57(2):1127--1156,
  2019.

\bibitem{chen2023block}
G.~Chen, R.~Vazquez, J.~Qiao, and M.~Krstic.
\newblock Block backstepping for isotachic hyperbolic {PDEs} and multilayer
  {T}imoshenko beams.
\newblock {\em arXiv preprint arXiv:2310.11416}, 2023.

\bibitem{de2024backstepping}
G.~A. de~Andrade, R.~Vazquez, I.~Karafyllis, and M.~Krstic.
\newblock Backstepping control of a hyperbolic {PDE} system with zero
  characteristic speed states.
\newblock {\em IEEE Transactions on Automatic Control}, 2024.

\bibitem{DAlembert}
J.~D'Alembert.
\newblock Suite des recherches sur la courbe que forme une corde tendue, mise
  en vibration.
\newblock {\em Histoire de l'Académie Royale des Sciences et des Belles
  Lettres de Berlin}, pages 220--249, 1749.

\bibitem{russell1991neutral}
D.~L. Russell.
\newblock Neutral {PDE} canonical representations of hyperbolic systems.
\newblock {\em The Journal of Integral Equations and Applications}, pages
  129--166, 1991.

\bibitem{karafyllis2014relation}
I.~Karafyllis and M.~Krstic.
\newblock On the relation of delay equations to first-order hyperbolic partial
  differential equations.
\newblock {\em ESAIM: Control, Optimisation and Calculus of Variations},
  20(3):894--923, 2014.

\bibitem{halebook}
J.~K. Hale and S.~M. Verduyn~Lunel.
\newblock {\em Introduction to functional differential equations}.
\newblock Applied Mathematical Sciences. Springer, 1993.

\bibitem{michiels2009strong}
W.~Michiels, T.~Vyhl{\'\i}dal, P.~Z{\'\i}tek, H.~Nijmeijer, and D.~Henrion.
\newblock Strong stability of neutral equations with an arbitrary delay
  dependency structure.
\newblock {\em {SIAM} Journal on Control and Optimization}, 48(2):763--786,
  2009.

\bibitem{Auriol2019a}
J.~Auriol and F.~Di~Meglio.
\newblock An explicit mapping from linear first order hyperbolic {P}{D}{E}s to
  difference systems.
\newblock {\em Systems \& Control Letters}, 123:144--150, 2019.

\bibitem{auriol2024contributions}
J.~Auriol.
\newblock {\em Contributions to the robust stabilization of networks of
  hyperbolic systems}.
\newblock HDR, Universit{\'e} Paris Saclay, 2024.

\bibitem{damak2015stability}
S.~Damak, M.~Di~Loreto, and S.~Mondi{\'e}.
\newblock Stability of linear continuous-time difference equations with
  distributed delay: Constructive exponential estimates.
\newblock {\em International Journal of Robust and Nonlinear Control},
  25(17):3195--3209, 2015.

\bibitem{niculescu2001delay}
S.-I. Niculescu.
\newblock {\em Delay effects on stability: a robust control approach}, volume
  269 of {\em Lecture Notes in Control and Information Sciences}.
\newblock Springer, 2001.

\bibitem{henry1974linear}
D.~Henry.
\newblock Linear autonomous neutral functional differential equations.
\newblock {\em Journal of Differential Equations}, 15(1):106--128, 1974.

\bibitem{saba2019stability}
D.~Bou~Saba, F.~Bribiesca-Argomedo, J.~Auriol, M.~Di~Loreto, and F.~Di~Meglio.
\newblock Stability analysis for a class of linear $2 \times 2$ hyperbolic
  {PDE}s using a backstepping transform.
\newblock {\em IEEE Transactions on Automatic Control}, 65(7):2941--2956, 2019.

\bibitem{logemann1996conditions}
H.~Logemann, R.~Rebarber, and G.~Weiss.
\newblock Conditions for robustness and nonrobustness of the stability of
  feedback systems with respect to small delays in the feedback loop.
\newblock {\em {SIAM} Journal on Control and Optimization}, 34(2):572--600,
  1996.

\bibitem{datko1986example}
R.~Datko, J.~Lagnese, and M.~P. Polis.
\newblock An example on the effect of time delays in boundary feedback
  stabilization of wave equations.
\newblock {\em {SIAM} Journal on Control and Optimization}, 24(1):152--156,
  1986.

\bibitem{Auriol2018}
J.~Auriol, U.~J.~F. Aarsnes, P.~Martin, and F.~Di~Meglio.
\newblock Delay-robust control design for heterodirectional linear coupled
  hyperbolic {P}{D}{E}s.
\newblock {\em IEEE Transactions on Automatic Control}, 63(10):3551--3557,
  2018.

\bibitem{Auriol2023}
J.~Auriol, F.~Bribiesca-Argomedo, and F.~Di~Meglio.
\newblock Robustification of stabilizing controllers for {ODE}-{PDE}-{ODE}
  systems: a filtering approach.
\newblock {\em Automatica}, 147:110724, 2023.

\bibitem{pepe2005asymptotic}
P.~Pepe.
\newblock On the asymptotic stability of coupled delay differential and
  continuous time difference equations.
\newblock {\em Automatica}, 41(1):107--112, 2005.

\bibitem{fridman2002stability}
E.~Fridman.
\newblock Stability of linear descriptor systems with delay: a {L}yapunov-based
  approach.
\newblock {\em Journal of Mathematical Analysis and Applications},
  273(1):24--44, 2002.

\bibitem{kyllingstad2009new}
A.~Kyllingstad and P.~J. Nessj{\o}en.
\newblock A new stick-slip prevention system.
\newblock In {\em SPE/IADC Drilling Conference and Exhibition}. Society of
  Petroleum Engineers, 2009.

\bibitem{curtain2012introduction}
R.~F. Curtain and H.~Zwart.
\newblock {\em An introduction to infinite-dimensional linear systems theory},
  volume~21 of {\em Texts in Applied Mathematics}.
\newblock Springer, 2012.

\bibitem{Auriol2018c}
J.~Auriol.
\newblock {\em Robust design of backstepping controllers for systems of linear
  hyperbolic {PDE}s}.
\newblock PhD thesis, PSL Research University, 2018.

\bibitem{Auriol2020e}
J.~Auriol and F.~Di~Meglio.
\newblock Robust output feedback stabilization for two heterodirectional linear
  coupled hyperbolic {PDE}s.
\newblock {\em Automatica}, 115:108896, 2020.

\bibitem{Saldivar2016}
B.~Saldivar, S.~Mondi{\'e}, S.-I. Niculescu, H.~Mounier, and I.~Boussaada.
\newblock A control oriented guided tour in oilwell drilling vibration
  modeling.
\newblock {\em Annual Reviews in Control}, 42:100 -- 113, 2016.

\bibitem{Tang2011}
S.~Tang and C.~Xie.
\newblock State and output feedback boundary control for a coupled {PDE}--{ODE}
  system.
\newblock {\em Systems \& Control Letters}, 60(8):540--545, 2011.

\bibitem{de2020backstepping}
G.~A. de~Andrade, R.~Vazquez, and D.~J. Pagano.
\newblock Backstepping-based estimation of thermoacoustic oscillations in a
  {Rijke} tube with experimental validation.
\newblock {\em IEEE Transactions on Automatic Control}, 65(12):5336--5343,
  2020.

\bibitem{smith1959controller}
O.~Smith.
\newblock A controller to overcome dead time.
\newblock {\em ISA Journal}, 6(2):28--33, 1959.

\bibitem{artstein1982linear}
Z.~Artstein.
\newblock {Linear systems with delayed controls: a reduction}.
\newblock {\em IEEE Transactions on Automatic Control}, 27(4):869--879, 1982.

\bibitem{Manitius1979}
A.~Manitius and A.~Olbrot.
\newblock Finite spectrum assignment problem for systems with delays.
\newblock {\em IEEE Transactions on Automatic Control}, 24(4):541--552, 1979.

\bibitem{deng2022predictor}
Y.~Deng, V.~L{\'e}chapp{\'e}, E.~Moulay, Z.~Chen, B.~Liang, F.~Plestan, and
  Q.-L. Han.
\newblock Predictor-based control of time-delay systems: a survey.
\newblock {\em International Journal of Systems Science}, 53(12):2496--2534,
  2022.

\bibitem{krstic2010lyapunov}
M.~Krstic.
\newblock Lyapunov stability of linear predictor feedback for time-varying
  input delay.
\newblock {\em IEEE Transactions on Automatic Control}, 55(2):554--559, 2010.

\bibitem{bekiaris2011lyapunov}
N.~Bekiaris-Liberis and M.~Krstic.
\newblock Lyapunov stability of linear predictor feedback for distributed input
  delays.
\newblock {\em IEEE Transactions on Automatic Control}, 56(3):655--660, 2011.

\bibitem{zhu2020predictor}
Y.~Zhu, M.~Krstic, and H.~Su.
\newblock Predictor feedback for uncertain linear systems with distributed
  input delays.
\newblock {\em IEEE Transactions on Automatic Control}, 65(12):5344--5351,
  2020.

\bibitem{holloway2019prescribed}
J.~Holloway and M.~Krstic.
\newblock Prescribed-time output feedback for linear systems in controllable
  canonical form.
\newblock {\em Automatica}, 107:77--85, 2019.

\bibitem{espitia2021predictor}
N.~Espitia and W.~Perruquetti.
\newblock Predictor-feedback prescribed-time stabilization of {LTI} systems
  with input delay.
\newblock {\em IEEE Transactions on Automatic Control}, 67(6):2784--2799, 2021.

\bibitem{bekiaris2016predictor}
N.~Bekiaris-Liberis and M.~Krstic.
\newblock Predictor-feedback stabilization of multi-input nonlinear systems.
\newblock {\em IEEE Transactions on Automatic Control}, 62(2):516--531, 2016.

\bibitem{bekiaris2012compensation}
N.~Bekiaris-Liberis and M.~Krstic.
\newblock Compensation of time-varying input and state delays for nonlinear
  systems.
\newblock {\em Journal of Dynamic Systems, Measurement, and Control},
  134(1):011009, 2011.

\bibitem{bresch2018robust}
D.~Bresch-Pietri, F.~Mazenc, and N.~Petit.
\newblock Robust compensation of a chattering time-varying input delay with
  jumps.
\newblock {\em Automatica}, 92:225--234, 2018.

\bibitem{zekraoui2023finite}
S.~Zekraoui, N.~Espitia, and W.~Perruquetti.
\newblock Finite/fixed-time stabilization of a chain of integrators with input
  delay via {PDE}-based nonlinear backstepping approach.
\newblock {\em Automatica}, 155:111095, 2023.

\bibitem{krstic2008lyapunov}
M.~Krstic.
\newblock Lyapunov tools for predictor feedbacks for delay systems: Inverse
  optimality and robustness to delay mismatch.
\newblock {\em Automatica}, 44(11):2930--2935, 2008.

\bibitem{kong2022prediction}
S.~Kong and D.~Bresch-Pietri.
\newblock Prediction-based controller for linear systems with stochastic input
  delay.
\newblock {\em Automatica}, 138:110149, 2022.

\bibitem{bekiaris2013robustness}
N.~Bekiaris-Liberis and M.~Krstic.
\newblock Robustness of nonlinear predictor feedback laws to time-and
  state-dependent delay perturbations.
\newblock {\em Automatica}, 49(6):1576--1590, 2013.

\bibitem{krstic2009delay}
M.~Krstic and D.~Bresch-Pietri.
\newblock Delay-adaptive full-state predictor feedback for systems with unknown
  long actuator delay.
\newblock In {\em 2009 American control conference}, pages 4500--4505. IEEE,
  2009.

\bibitem{bresch2010delay}
D.~Bresch-Pietri and M.~Krstic.
\newblock Delay-adaptive predictor feedback for systems with unknown long
  actuator delay.
\newblock {\em IEEE Transactions on Automatic Control}, 55(9):2106--2112, 2010.

\bibitem{zhu2018pde}
Y.~Zhu, M.~Krstic, and H.~Su.
\newblock {PDE} boundary control of multi-input {LTI} systems with distinct and
  uncertain input delays.
\newblock {\em IEEE Transactions on Automatic Control}, 63(12):4270--4277,
  2018.

\bibitem{bresch2012adaptive}
D.~Bresch-Pietri, J.~Chauvin, and N.~Petit.
\newblock Adaptive control scheme for uncertain time-delay systems.
\newblock {\em Automatica}, 48(8):1536--1552, 2012.

\bibitem{chakraborty2017control}
I.~Chakraborty, S.~Mehta, E.~Doucette, and W.~Dixon.
\newblock Control of an input delayed uncertain nonlinear system with adaptive
  delay estimation.
\newblock In {\em 2017 American Control Conference (ACC)}, pages 1779--1784.
  IEEE, 2017.

\bibitem{DiMeglio2018}
F.~Di~Meglio, F.~Bribiesca-Argomedo, L.~Hu, and M.~Krstic.
\newblock Stabilization of coupled linear heterodirectional hyperbolic
  {P}{D}{E}--{O}{D}{E} systems.
\newblock {\em Automatica}, 87:281--289, 2018.

\bibitem{Auriol2018a}
J.~Auriol, F.~Bribiesca-Argomedo, D.~Bou~Saba, M.~Di~Loreto, and F.~Di~Meglio.
\newblock Delay-robust stabilization of a hyperbolic {P}{D}{E}--{O}{D}{E}
  system.
\newblock {\em Automatica}, 95:494--502, 2018.

\bibitem{Auriol2019d}
J.~Auriol and F.~Bribiesca-Argomedo.
\newblock Delay-robust stabilization of an $n + m$ hyperbolic {PDE}-{ODE}
  system.
\newblock In {\em Conference on Decision and Control}, pages 4964--4970. IEEE,
  2019.

\bibitem{deutscher2019output}
J.~Deutscher, N.~Gehring, and R.~Kern.
\newblock Output feedback control of general linear heterodirectional
  hyperbolic {PDE}-{ODE} systems with spatially-varying coefficients.
\newblock {\em International Journal of Control}, 92(10):2274--2290, 2019.

\bibitem{gehring2023control}
N.~Gehring, A.~Irscheid, J.~Deutscher, F.~Woittennek, and J.~Rudolph.
\newblock Control of distributed-parameter systems using normal forms: an
  introduction.
\newblock {\em at - Automatisierungstechnik}, 71(8):624--646, 2023.

\bibitem{irscheid2023output}
A.~Irscheid, J.~Deutscher, N.~Gehring, and J.~Rudolph.
\newblock Output regulation for general heterodirectional linear hyperbolic
  {PDE}s coupled with nonlinear {ODE}s.
\newblock {\em Automatica}, 148:110748, 2023.

\bibitem{mathiyalagan2022observer}
K.~Mathiyalagan, A.~Nidhi, H.~Su, and T.~Renugadevi.
\newblock Observer and boundary output feedback control for coupled
  {ODE}-transport {PDE}.
\newblock {\em Applied Mathematics and Computation}, 426:127096, 2022.

\bibitem{Hasan2016}
A.~Hasan, O.~M. Aamo, and M.~Krstic.
\newblock Boundary observer design for hyperbolic {PDE}--{ODE} cascade systems.
\newblock {\em Automatica}, 68:75--86, 2016.

\bibitem{diagne2017time}
M.~Diagne, N.~Bekiaris-Liberis, and M.~Krstic.
\newblock Time-and state-dependent input delay-compensated bang-bang control of
  a screw extruder for {3D} printing.
\newblock {\em International Journal of Robust and Nonlinear Control},
  27(17):3727--3757, 2017.

\bibitem{diagne2017compensation}
M.~Diagne, N.~Bekiaris-Liberis, and M.~Krstic.
\newblock Compensation of input delay that depends on delayed input.
\newblock {\em Automatica}, 85:362--373, 2017.

\bibitem{diagne2017control}
M.~Diagne, N.~Bekiaris-Liberis, A.~Otto, and M.~Krstic.
\newblock Control of transport {PDE}/nonlinear {ODE} cascades with
  state-dependent propagation speed.
\newblock {\em IEEE Transactions on Automatic Control}, 62(12):6278--6293,
  2017.

\bibitem{roman2018backstepping}
C.~Roman, D.~Bresch-Pietri, E.~Cerpa, C.~Prieur, and O.~Sename.
\newblock Backstepping control of a wave {PDE} with unstable source terms and
  dynamic boundary.
\newblock {\em IEEE Control Systems Letters}, 2(3):459--464, 2018.

\bibitem{deutscher2021backstepping}
J.~Deutscher and J.~Gabriel.
\newblock A backstepping approach to output regulation for coupled linear
  wave--{ODE} systems.
\newblock {\em Automatica}, 123:109338, 2021.

\bibitem{BouSaba2017ODE-PDE-ODE}
D.~Bou~Saba, F.~Bribiesca-Argomedo, M.~Di~Loreto, and D.~Eberard.
\newblock Backstepping stabilization of 2×2 linear hyperbolic {PDEs} coupled
  with potentially unstable actuator and load dynamics.
\newblock In {\em 2017 IEEE 56th Annual Conference on Decision and Control
  (CDC)}, pages 2498--2503, 2017.

\bibitem{anfinsen2018stabilization}
H.~Anfinsen and O.~M. Aamo.
\newblock Stabilization of a linear hyperbolic {PDE} with actuator and sensor
  dynamics.
\newblock {\em Automatica}, 95:104--111, 2018.

\bibitem{DiMeglio2020}
F.~Di~Meglio, P.-O. Lamare, and U.~J.~F. Aarsnes.
\newblock Robust output feedback stabilization of an {ODE}--{PDE}--{ODE}
  interconnection.
\newblock {\em Automatica}, 119:109059, 2020.

\bibitem{Deutscher2018}
J.~Deutscher, N.~Gehring, and R.~Kern.
\newblock Output feedback control of general linear heterodirectional
  hyperbolic {ODE}--{PDE}--{ODE} systems.
\newblock {\em Automatica}, 95:472--480, 2018.

\bibitem{isidori2013nonlinear}
A.~Isidori.
\newblock {\em Nonlinear control systems {II}}.
\newblock Communications and Control Engineering. Springer, 1999.

\bibitem{BouSaba2019_b}
D.~Bou~Saba, F.~Bribiesca-Argomedo, M.~Di~Loreto, and D.~Eberard.
\newblock Strictly proper control design for the stabilization of $2\times2$
  linear hyperbolic {ODE-PDE-ODE} systems.
\newblock In {\em 2019 IEEE 58th Conference on Decision and Control (CDC)},
  pages 4996--5001. IEEE, 2019.

\bibitem{Wang2020}
J.~Wang and M.~Krstic.
\newblock Delay-compensated control of sandwiched {ODE}--{PDE}--{ODE}
  hyperbolic systems for oil drilling and disaster relief.
\newblock {\em Automatica}, 120:109131, 2020.

\bibitem{Auriol2022b}
J.~Auriol and F.~Bribiesca-Argomedo.
\newblock Observer design for $n+ m$ linear hyperbolic {ODE}-{PDE}-{ODE}
  systems.
\newblock {\em IEEE Control Systems Letters}, 7:283--288, 2022.

\bibitem{Andrade2018}
G.~A. de~Andrade, R.~Vazquez, and D.~J. Pagano.
\newblock Backstepping stabilization of a linearized {O}{D}{E}--{P}{D}{E}
  {R}ijke tube model.
\newblock {\em Automatica}, 96:98--109, 2018.

\bibitem{Auriol2020d}
J.~Auriol.
\newblock Output feedback stabilization of an underactuated cascade network of
  interconnected linear {PDE} systems using a backstepping approach.
\newblock {\em Automatica}, 117:108964, 2020.

\bibitem{Auriol2022}
J.~Auriol and D.~Bresch-Pietri.
\newblock Robust state-feedback stabilization of an underactuated network of
  interconnected $n+ m$ hyperbolic {PDE} systems.
\newblock {\em Automatica}, 136:110040, 2022.

\bibitem{zhang2023robust}
J.~Zhang and J.~Qi.
\newblock Robust stabilization of 2$\times$ 2 first-order hyperbolic {PDE}s
  with uncertain input delay.
\newblock {\em Automatica}, 157:111235, 2023.

\bibitem{Redaud2021}
J.~Redaud, J.~Auriol, and S.-I. Niculescu.
\newblock Output-feedback control of an underactuated network of interconnected
  hyperbolic {PDE}-{ODE} systems.
\newblock {\em Systems \& Control Letters}, 154:104984, 2021.

\bibitem{Aarsnes2019c}
U.~J.~F. Aarsnes, R.~Vazquez, F.~Di~Meglio, and M.~Krstic.
\newblock Delay robust control design of under-actuated {PDE}-{ODE}-{PDE}
  systems.
\newblock In {\em 2019 American Control Conference (ACC)}, pages 3200--3205.
  IEEE, 2019.

\bibitem{Redaud2022Fredholm}
J.~Redaud, J.~Auriol, and S.-I. Niculescu.
\newblock Stabilizing output-feedback control law for hyperbolic systems using
  a {Fredholm} transformation.
\newblock {\em IEEE Transactions on Automatic Control}, 67(12):6651--6666,
  2022.

\bibitem{coc2009}
J.~Cochran and M.~Krstic.
\newblock Motion planning and trajectory tracking for three-dimensional
  {Poiseuille} flow.
\newblock {\em Journal of Fluid Mechanics}, 626:307--332, 2009.

\bibitem{bresch2009adaptive}
D.~Bresch-Pietri and M.~Krstic.
\newblock Adaptive trajectory tracking despite unknown input delay and plant
  parameters.
\newblock {\em Automatica}, 45(9):2074--2081, 2009.

\bibitem{Aamo2013}
O.~M. Aamo.
\newblock Disturbance rejection in 2x2 linear hyperbolic systems.
\newblock {\em IEEE Transactions on Automatic Control}, 58(5):1095 -- 1106, May
  2013.

\bibitem{Anfinsen2015}
H.~Anfinsen and O.~M. Aamo.
\newblock Disturbance rejection in the interior domain of linear $2 \times 2$
  hyperbolic systems.
\newblock {\em IEEE Transactions on Automatic Control}, 60(1):186--191, 2015.

\bibitem{Anfinsen2017}
H.~Anfinsen and O.~M. Aamo.
\newblock Model reference adaptive control of $n+ 1$ coupled linear hyperbolic
  {PDE}s.
\newblock {\em Systems \& Control Letters}, 109:1--11, 2017.

\bibitem{Deutscher2017}
J.~Deutscher.
\newblock Finite-time output regulation for linear $2\times 2$ hyperbolic
  systems using backstepping.
\newblock {\em Automatica}, 75:54--62, 2017.

\bibitem{Deutscher2017c}
J.~Deutscher.
\newblock Backstepping design of robust state feedback regulators for linear $2
  \times 2$ hyperbolic systems.
\newblock {\em IEEE Transactions on Automatic Control}, 62(10):5240--5247,
  2017.

\bibitem{Deutscher2017b}
J.~Deutscher.
\newblock Output regulation for general linear heterodirectional hyperbolic
  systems with spatially-varying coefficients.
\newblock {\em Automatica}, 85:34--42, 2017.

\bibitem{Xu2017}
X.~Xu and S.~Dubljievic.
\newblock Output regulation for a class of linear boundary controlled
  first-order hyperbolic {PIDE} systems.
\newblock {\em Automatica}, 85:43--52, 2017.

\bibitem{Gu2018coupledwave}
J.-J. Gu, J.-M. Wang, and Y.-P. Guo.
\newblock Output regulation of anti-stable coupled wave equations via the
  backstepping technique.
\newblock {\em IET Control Theory \& Applications}, 12(4):431--445, 2018.

\bibitem{Redaud2022b}
J.~Redaud, F.~Bribiesca-Argomedo, and J.~Auriol.
\newblock Practical output regulation and tracking for linear {ODE}-hyperbolic
  {PDE}-{ODE} systems.
\newblock In {\em Advances in distributed parameter systems}, Advances in
  Delays and Dynamics, pages 143--169. Springer, 2022.

\bibitem{redaud2024tracking}
J.~Redaud, F.~Bribiesca-Argomedo, and J.~Auriol.
\newblock Output regulation and tracking for linear {ODE}-hyperbolic
  {PDE}–{ODE} systems.
\newblock {\em Automatica}, 162:111503, 2024.

\bibitem{Meurer2009}
T.~Meurer and A.~Kugi.
\newblock Tracking control for boundary controlled parabolic {PDE}s with
  varying parameters: Combining backstepping and differential flatness.
\newblock {\em Automatica}, 45(5):1182--1194, 2009.

\bibitem{Deutscher015outputparabolic}
J.~Deutscher.
\newblock A backstepping approach to the output regulation of boundary
  controlled parabolic {PDE}s.
\newblock {\em Automatica}, 57:56--64, 2015.

\bibitem{Deutscher2021cooperativeoutputregulation}
J.~Deutscher.
\newblock Cooperative output regulation for a network of parabolic systems with
  varying parameters.
\newblock {\em Automatica}, 125:109446, 2021.

\bibitem{Baccoli2015}
A.~Baccoli, A.~Pisano, and Y.~Orlov.
\newblock Boundary control of coupled reaction--diffusion processes with
  constant parameters.
\newblock {\em Automatica}, 54:80--90, 2015.

\bibitem{Vazquez2017}
R.~Vazquez and M.~Krstic.
\newblock Boundary control of coupled reaction-advection-diffusion systems with
  spatially-varying coefficients.
\newblock {\em IEEE Transactions on Automatic Control}, 62:2026--2033, 2017.

\bibitem{camacho2020boundary}
L.~Camacho-Solorio, R.~Vazquez, and M.~Krstic.
\newblock Boundary observers for coupled diffusion--reaction systems with
  prescribed convergence rate.
\newblock {\em Systems \& Control Letters}, 135:104586, 2020.

\bibitem{deutscher2018coupled}
J.~Deutscher and S.~Kerschbaum.
\newblock Backstepping control of coupled linear parabolic {PIDE}s with
  spatially varying coefficients.
\newblock {\em IEEE Transactions on Automatic Control}, 63(12):4218--4233,
  2018.

\bibitem{ker2019}
S.~Kerschbaum and J.~Deutscher.
\newblock Backstepping control of coupled linear parabolic {PDE}s with space
  and time dependent coefficients.
\newblock {\em IEEE Transactions on Automatic Control}, 65(7):3060--3067, 2019.

\bibitem{krs2009}
M.~Krstic.
\newblock Compensating actuator and sensor dynamics governed by diffusion
  {PDE}s.
\newblock {\em Systems \& Control Letters}, 58(5):372--377, 2009.

\bibitem{Deutscher2020}
J.~Deutscher and N.~Gehring.
\newblock Output feedback control of coupled linear parabolic
  {ODE}--{PDE}--{ODE} systems.
\newblock {\em IEEE Transactions on Automatic Control}, 66(10):4668--4683,
  2020.

\bibitem{Xu2023}
X.~Xu, L.~Liu, M.~Krstic, and G.~Feng.
\newblock Stabilization of chains of linear parabolic {PDE}--{ODE} cascades.
\newblock {\em Automatica}, 148:110763, 2023.

\bibitem{krstic2009control}
M.~Krstic.
\newblock Control of an unstable reaction--diffusion {PDE} with long input
  delay.
\newblock {\em Systems \& Control Letters}, 58(10-11):773--782, 2009.

\bibitem{Chen2017}
S.~Chen, R.~Vazquez, and M.~Krstic.
\newblock Stabilization of an underactuated coupled transport-wave {PDE}
  system.
\newblock In {\em 2017 American Control Conference (ACC)}, pages 2504--2509.
  IEEE, 2017.

\bibitem{ghousein2020backstepping}
M.~Ghousein and E.~Witrant.
\newblock Backstepping control for a class of coupled hyperbolic-parabolic
  {PDE} systems.
\newblock In {\em 2020 American Control Conference (ACC)}, pages 1600--1605.
  IEEE, 2020.

\bibitem{Chen2023b}
G.~Chen, R.~Vazquez, Z.~Liu, and H.~Su.
\newblock Backstepping control of an underactuated hyperbolic--parabolic
  coupled {PDE} system.
\newblock {\em IEEE Transactions on Automatic Control}, 69(2):1218--1225, 2023.

\bibitem{deutscher2023backstepping}
J.~Deutscher, N.~Gehring, and N.~Jung.
\newblock Backstepping control of coupled general hyperbolic-parabolic
  {PDE}-{PDE} systems.
\newblock {\em IEEE Transactions on Automatic Control}, 69(5):3465--3472, 2024.

\bibitem{Vazquez2006}
R.~Vazquez and M.~Krstic.
\newblock Explicit integral operator feedback for local stabilization of
  nonlinear thermal convection loop {PDE}s.
\newblock {\em Systems \& Control Letters}, 55(8):624--632, 2006.

\bibitem{convloop}
R.~Vazquez and M.~Krstic.
\newblock Boundary observer for output-feedback stabilization of thermal-fluid
  convection loop.
\newblock {\em IEEE Transactions on Control Systems Technology},
  18(4):789--797, 2010.

\bibitem{Vazquez2007}
R.~Vazquez and M.~Krstic.
\newblock A closed-form feedback controller for stabilization of the linearized
  2-{D} {Navier--Stokes Poiseuille} system.
\newblock {\em IEEE Transactions on Automatic Control}, 52(12):2298--2312,
  2007.

\bibitem{Xu2008}
C.~Xu, E.~Schuster, R.~Vazquez, and M.~Krstic.
\newblock Stabilization of linearized 2{D} magnetohydrodynamic channel flow by
  backstepping boundary control.
\newblock {\em Systems \& Control Letters}, 57(10):805--812, 2008.

\bibitem{Alalabi2025}
A.~Alalabi and K.~Morris.
\newblock Boundary control and observer design via backstepping for a coupled
  parabolic–elliptic system.
\newblock {\em Automatica}, 174:112154, 2025.

\bibitem{aarsnes2016methodology}
U.~J.~F. Aarsnes, F.~Di~Meglio, R.~Graham, and O.~M. Aamo.
\newblock A methodology for classifying operating regimes in
  underbalanced-drilling operations.
\newblock {\em SPE Journal}, 21(02):423--433, 2016.

\bibitem{Vazquez2016}
R.~Vazquez and M.~Krstic.
\newblock Bilateral boundary control of one-dimensional first-and second-order
  {P}{D}{E}s using infinite-dimensional backstepping.
\newblock In {\em Conference on Decision and Control (CDC)}, pages 537--542.
  IEEE, 2016.

\bibitem{bekiaris2019nonlinear}
N.~Bekiaris-Liberis and R.~Vazquez.
\newblock Nonlinear bilateral output-feedback control for a class of viscous
  {H}amilton--{J}acobi {PDE}s.
\newblock {\em Automatica}, 101:223--231, 2019.

\bibitem{chen2021folding}
S.~Chen, R.~Vazquez, and M.~Krstic.
\newblock Folding bilateral backstepping output-feedback control design for an
  unstable parabolic {PDE}.
\newblock {\em IEEE Transactions on Automatic Control}, 67(5):2389--2404, 2021.

\bibitem{Auriol2018b}
J.~Auriol and F.~Di~Meglio.
\newblock Two-sided boundary stabilization of heterodirectional linear coupled
  hyperbolic {PDE}s.
\newblock {\em IEEE Transactions on Automatic Control}, 63(8):2421--2436, 2018.

\bibitem{Coron2016}
J.-M. Coron, L.~Hu, and G.~Olive.
\newblock Stabilization and controllability of first-order integro-differential
  hyperbolic equations.
\newblock {\em Journal of Functional Analysis}, 271(12):3554--3587, 2016.

\bibitem{wilhelmsen2021minimum}
N.~Wilhelmsen, H.~Anfinsen, and O.~M. Aamo.
\newblock Minimum time observer designs for $n+ m$ linear hyperbolic systems
  with unilateral, bilateral or pointwise in-domain sensing.
\newblock {\em European Journal of Control}, 61:50--67, 2021.

\bibitem{Su2017}
L.~Su, W.~Guo, J.-M. Wang, and M.~Krstic.
\newblock Boundary stabilization of wave equation with velocity recirculation.
\newblock {\em IEEE Transactions on Automatic Control}, 62(9):4760--4767, 2017.

\bibitem{Chen2023}
G.~Chen, R.~Vazquez, and M.~Krstic.
\newblock Rapid stabilization of {Timoshenko} beam by {PDE} backstepping.
\newblock {\em IEEE Transactions on Automatic Control}, 69(2):1141--1148, 2023.

\bibitem{Smyshlyaevfirstadaptivecdc}
A.~Smyshlyaev and M.~Krstic.
\newblock Output-feedback adaptive control for parabolic {PDE}s with spatially
  varying coefficients.
\newblock In {\em Proceedings of the 45th IEEE Conference on Decision and
  Control}, pages 3099--3104, 2006.

\bibitem{KrsticadaptivepartI}
M.~Krstic and A.~Smyshlyaev.
\newblock Adaptive boundary control for unstable parabolic {PDE}s—{P}art {I}:
  {L}yapunov design.
\newblock {\em IEEE Transactions on Automatic Control}, 53(7):1575--1591, 2008.

\bibitem{SmyshlyaevadaptivepartII}
A.~Smyshlyaev and M.~Krstic.
\newblock Adaptive boundary control for unstable parabolic {PDE}s—{P}art
  {II}: Estimation-based designs.
\newblock {\em Automatica}, 43(9):1543--1556, 2007.

\bibitem{SmyshlyaevadaptivepartIII}
A.~Smyshlyaev and M.~Krstic.
\newblock Adaptive boundary control for unstable parabolic {PDE}s—{P}art
  {III}: Output feedback examples with swapping identifiers.
\newblock {\em Automatica}, 43(9):1557--1564, 2007.

\bibitem{krsticannualreviewsadaptive}
M.~Krstic and A.~Smyshlyaev.
\newblock Adaptive control of {PDE}s.
\newblock {\em Annual Reviews in Control}, 32(2):149--160, 2008.

\bibitem{Krstic2010adaptivewave}
M.~Krstic.
\newblock Adaptive control of an anti-stable wave {PDE}.
\newblock {\em Dynamics of Continuous, Discrete and Impulsive Systems, Series
  A: Mathematical Analysis}, 17(6):853--882, 2010.

\bibitem{Bernard2014}
P.~Bernard and M.~Krstic.
\newblock Adaptive output-feedback stabilization of non-local hyperbolic
  {PDE}s.
\newblock {\em Automatica}, 50(10):2692--2699, 2014.

\bibitem{BreschPietri2014a}
D.~Bresch-Pietri and M.~Krstic.
\newblock Adaptive output feedback for oil drilling stick-slip instability
  modeled by wave {P}{D}{E} with anti-damped dynamic boundary.
\newblock In {\em American Control Conference (ACC), 2014}, pages 386--391.
  IEEE, 2014.

\bibitem{Anfinsen2016a}
H.~Anfinsen, M.~Diagne, O.~M. Aamo, and M.~Krstic.
\newblock Boundary parameter and state estimation in general linear hyperbolic
  {PDE}s.
\newblock {\em IFAC-PapersOnLine}, 49(8):104--110, 2016.

\bibitem{oliveira2022extremum}
T.~R. Oliveira and M.~Krstic.
\newblock {\em Extremum seeking through delays and PDEs}.
\newblock SIAM, 2022.

\bibitem{bamieh2002distributed}
B.~Bamieh, F.~Paganini, and M.~Dahleh.
\newblock Distributed control of spatially invariant systems.
\newblock {\em IEEE Transactions on Automatic Control}, 47(7):1091--1107, 2002.

\bibitem{Vazquez2016disk}
R.~Vazquez and M.~Krstic.
\newblock Boundary control of a singular reaction-diffusion equation on a disk.
\newblock {\em IFAC-PapersOnLine}, 49(8):74--79, 2016.

\bibitem{vazquez2019sphere}
R.~Vazquez and M.~Krstic.
\newblock Boundary control and estimation of reaction--diffusion equations on
  the sphere under revolution symmetry conditions.
\newblock {\em International Journal of Control}, 92(1):2--11, 2019.

\bibitem{vazquez2016explicit}
R.~Vazquez and M.~Krstic.
\newblock Explicit output-feedback boundary control of reaction-diffusion
  {PDE}s on arbitrary-dimensional balls.
\newblock {\em ESAIM: Control, Optimisation and Calculus of Variations},
  22(4):1078--1096, 2016.

\bibitem{R.2022ball}
R.~Vazquez, J.~Zhang, J.~Qi, and M.~Krstic.
\newblock Kernel well-posedness and computation by power series in backstepping
  output feedback for radially-dependent reaction--diffusion {PDE}s on
  multidimensional balls.
\newblock {\em Systems $\&$ Control Letters}, 177:105538, 2023.

\bibitem{liu2020boundary}
X.~Liu and C.~Xie.
\newblock Boundary control of reaction--diffusion equations on
  higher-dimensional symmetric domains.
\newblock {\em Automatica}, 114:108832, 2020.

\bibitem{vazquez-nonlinear3}
R.~Vazquez, J.-M. Coron, M.~Krstic, and G.~Bastin.
\newblock Collocated output-feedback stabilization of a $2\times 2$ quasilinear
  hyperbolic system using backstepping.
\newblock In {\em 2012 American Control Conference (ACC)}, pages 2202--2207,
  2012.

\bibitem{Krstic2009}
M.~Krstic, L.~Magnis, and R.~Vazquez.
\newblock Nonlinear control of the viscous {B}urgers equation: Trajectory
  generation, tracking, and observer design.
\newblock {\em Journal of Dynamic Systems, Measurement, and Control},
  131(2):021012, 2009.

\bibitem{krstic2008nonlinearB}
M.~Krstic, L.~Magnis, and R.~Vazquez.
\newblock Nonlinear stabilization of shock-like unstable equilibria in the
  viscous {Burgers PDE}.
\newblock {\em IEEE Transactions on Automatic Control}, 53(7):1678--1683, 2008.

\bibitem{BekiarisLiberis2014}
N.~Bekiaris-Liberis and M.~Krstic.
\newblock Compensation of wave actuator dynamics for nonlinear systems.
\newblock {\em IEEE Transactions on Automatic Control}, 59(6):1555--1570, June
  2014.

\bibitem{BekiarisLiberis2016}
N.~Bekiaris-Liberis and M.~Krstic.
\newblock Stability of predictor-based feedback for nonlinear systems with
  distributed input delay.
\newblock {\em Automatica}, 70:195--203, 2016.

\bibitem{karafyllis2019global}
I.~Karafyllis and M.~Krstic.
\newblock Global stabilization of a class of nonlinear reaction-diffusion
  partial differential equations by boundary feedback.
\newblock {\em {SIAM} Journal on Control and Optimization}, 57(6):3723--3748,
  2019.

\bibitem{vaz2008_a}
R.~Vazquez and M.~Krstic.
\newblock Control of 1-{D} parabolic {PDE}s with {V}olterra
  nonlinearities---{P}art {I}: {D}esign.
\newblock {\em Automatica}, 44:2778--2790, 2008.

\bibitem{vaz2008_b}
R.~Vazquez and M.~Krstic.
\newblock Control of 1-{D} parabolic {PDE}s with {V}olterra
  nonlinearities---{P}art {II}: {A}nalysis.
\newblock {\em Automatica}, 44:2791--2803, 2008.

\bibitem{smy2005b}
A.~Smyshlyaev and M.~Krstic.
\newblock On control design for {PDE}s with space-dependent diffusivity or
  time-dependent reactivity.
\newblock {\em Automatica}, 41(9):1601--1608, 2005.

\bibitem{vazquez2008control}
R.~Vazquez, E.~Tr{\'e}lat, and J.-M. Coron.
\newblock Control for fast and stable laminar-to-high-{Reynolds}-numbers
  transfer in a {2D} {Navier-Stokes} channel flow.
\newblock {\em Discrete and Continuous Dynamical Systems - B}, 10(4):925--956,
  2008.

\bibitem{jadachowski2012efficient}
L.~Jadachowski, T.~Meurer, and A.~Kugi.
\newblock An efficient implementation of backstepping observers for
  time-varying parabolic {PDE}s.
\newblock {\em IFAC Proceedings Volumes}, 45(2):798--803, 2012.

\bibitem{deu2016}
A.~Deutschmann, L.~Jadachowski, and A.~Kugi.
\newblock Backstepping-based boundary observer for a class of time-varying
  linear hyperbolic {PIDE}s.
\newblock {\em Automatica}, 68:369--377, 2016.

\bibitem{anfinsen2020stabilization}
H.~Anfinsen and O.~M. Aamo.
\newblock Stabilization and tracking control of a time-variant linear
  hyperbolic {PIDE} using backstepping.
\newblock {\em Automatica}, 116:108929, 2020.

\bibitem{Coron2021}
J.-M. Coron and H.-M. Nguyen.
\newblock On the optimal controllability time for linear hyperbolic systems
  with time-dependent coefficients.
\newblock {\em arXiv preprint arXiv:2103.02653}, 2021.

\bibitem{coron2021boundary}
J.-M. Coron, L.~Hu, G.~Olive, and P.~Shang.
\newblock Boundary stabilization in finite time of one-dimensional linear
  hyperbolic balance laws with coefficients depending on time and space.
\newblock {\em Journal of Differential Equations}, 271:1109--1170, 2021.

\bibitem{karafyllis2021event}
I.~Karafyllis, N.~Espitia, and M.~Krstic.
\newblock Event-triggered gain scheduling of reaction-diffusion {PDEs}.
\newblock {\em SIAM Journal on Control and Optimization}, 59(3):2047--2067,
  2021.

\bibitem{auriol2024event}
J.~Auriol and N.~Espitia.
\newblock Event-triggered gain scheduling of 2$\times$ 2 hyperbolic {PDEs} with
  time and space varying coupling coefficients.
\newblock In {\em IEEE Conference on Decision and Control}, 2024.

\bibitem{coron2017null}
J.-M. Coron and H.-M. Nguyen.
\newblock Null controllability and finite time stabilization for the heat
  equations with variable coefficients in space in one dimension via
  backstepping approach.
\newblock {\em Archive for Rational Mechanics and Analysis}, 225:993--1023,
  2017.

\bibitem{espitia2019boundary}
N.~Espitia, A.~Polyakov, D.~Efimov, and W.~Perruquetti.
\newblock Boundary time-varying feedbacks for fixed-time stabilization of
  constant-parameter reaction--diffusion systems.
\newblock {\em Automatica}, 103:398--407, 2019.

\bibitem{espitia2022sensor}
N.~Espitia, D.~Steeves, W.~Perruquetti, and M.~Krstic.
\newblock Sensor delay-compensated prescribed-time observer for {LTI} systems.
\newblock {\em Automatica}, 135:110005, 2022.

\bibitem{steeves2020prescribed}
D.~Steeves, M.~Krstic, and R.~Vazquez.
\newblock Prescribed--time estimation and output regulation of the linearized
  {S}chr{\"o}dinger equation by backstepping.
\newblock {\em European Journal of Control}, 55:3--13, 2020.

\bibitem{izadi2015backstepping}
M.~Izadi and S.~Dubljevic.
\newblock Backstepping output-feedback control of moving boundary parabolic
  {PDE}s.
\newblock {\em European Journal of Control}, 21:27--35, 2015.

\bibitem{CAUVINVILA2023251}
J.~Cauvin-Vila, V.~Ehrlacher, and A.~Hayat.
\newblock Boundary stabilization of one-dimensional cross-diffusion systems in
  a moving domain: Linearized system.
\newblock {\em Journal of Differential Equations}, 350:251--307, 2023.

\bibitem{kog2018}
S.~Koga, M.~Diagne, and M.~Krstic.
\newblock Control and state estimation of the one-phase {S}tefan problem via
  backstepping design.
\newblock {\em IEEE Transactions on Automatic Control}, 64(2):510--525, 2018.

\bibitem{kog2020_f}
S.~Koga and M.~Krstic.
\newblock Single-boundary control of the two-phase {S}tefan system.
\newblock {\em Systems \& Control Letters}, 135:104573, 2020.

\bibitem{szc2022}
M.~Szczesiak and H.~I. Basturk.
\newblock Adaptive boundary control for wave {PDE} on a domain with moving
  boundary and with unknown system parameters in the boundary dynamics.
\newblock {\em Automatica}, 145:110526, 2022.

\bibitem{SmyshlyaevEulerBernoulli2009}
A.~Smyshlyaev, B.-Z. Guo, and M.~Krstic.
\newblock Arbitrary decay rate for {Euler-Bernoulli} beam by backstepping
  boundary feedback.
\newblock {\em IEEE Transactions on Automatic Control}, 54(5):1134--1140, 2009.

\bibitem{BEKIARISLIBERIS2010DistributedWave}
N.~Bekiaris-Liberis and M.~Krstic.
\newblock Compensating the distributed effect of a wave {PDE} in the actuation
  or sensing path of {MIMO LTI} systems.
\newblock {\em Systems \& Control Letters}, 59(11):713--719, 2010.

\bibitem{Guo2014Fredholm}
C.~Guo, C.~Xie, and C.~Zhou.
\newblock Stabilization of a spatially non-causal reaction–diffusion equation
  by boundary control.
\newblock {\em International Journal of Robust and Nonlinear Control},
  24(1):1--17, 2014.

\bibitem{BribiescaACC2014}
F.~Bribiesca-Argomedo and M.~Krstic.
\newblock Backstepping-forwarding boundary control design for first-order
  hyperbolic systems with {Fredholm} integrals.
\newblock In {\em 2014 American Control Conference}, pages 5428--5433, 2014.

\bibitem{BribiescaArgomedo2015}
F.~Bribiesca-Argomedo and M.~Krstic.
\newblock Backstepping-forwarding control and observation for hyperbolic {PDE}s
  with {F}redholm integrals.
\newblock {\em IEEE Transactions on Automatic Control}, 60(8):2145--2160, 2015.

\bibitem{CORON2014Fredholm}
J.-M. Coron and Q.~Lü.
\newblock Local rapid stabilization for a {Korteweg–de Vries} equation with a
  {Neumann} boundary control on the right.
\newblock {\em Journal de Mathématiques Pures et Appliquées},
  102(6):1080--1120, 2014.

\bibitem{CORON2015Fredholm}
J.-M. Coron and Q.~Lü.
\newblock Fredholm transform and local rapid stabilization for a
  {Kuramoto–Sivashinsky} equation.
\newblock {\em Journal of Differential Equations}, 259(8):3683--3729, 2015.

\bibitem{Zhang2019InternalFredholm}
C.~Zhang.
\newblock Finite-time internal stabilization of a linear {1-D} transport
  equation.
\newblock {\em Systems \& Control Letters}, 133:104529, 2019.

\bibitem{Zhang2022FredholmInternal}
C.~Zhang.
\newblock Internal rapid stabilization of a 1-{D} linear transport equation
  with a scalar feedback.
\newblock {\em Mathematical Control and Related Fields}, 12(1):169--200, 2022.

\bibitem{coron2022stabilization}
J.-M. Coron, A.~Hayat, S.~Xiang, and C.~Zhang.
\newblock Stabilization of the linearized water tank system.
\newblock {\em Archive for Rational Mechanics and Analysis}, 244(3):1019--1097,
  2022.

\bibitem{Auriol2024Fredholm}
J.~Auriol.
\newblock Stabilization of integral delay equations by solving {Fredholm}
  equations.
\newblock {\em IEEE Control Systems Letters}, 8:676--681, 2024.

\bibitem{coron2018rapid}
J.-M. Coron, L.~Gagnon, and M.~Morancey.
\newblock Rapid stabilization of a linearized bilinear 1-{D} {S}chr{\"o}dinger
  equation.
\newblock {\em Journal de Math{\'e}matiques Pures et Appliqu{\'e}es},
  115:24--73, 2018.

\bibitem{Tsubakino2014Semiseparable}
D.~Tsubakino, F.~Bribiesca-Argomedo, and M.~Krstic.
\newblock Backstepping-forwarding control of parabolic {PDEs} with partially
  separable kernels.
\newblock In {\em 53rd IEEE Conference on Decision and Control}, pages
  5236--5241, 2014.

\bibitem{gagnon2021fredholm}
L.~Gagnon, P.~Lissy, and S.~Marx.
\newblock A {Fredholm} transformation for the rapid stabilization of a
  degenerate parabolic equation.
\newblock {\em {SIAM} Journal on Control and Optimization}, 59(5):3828--3859,
  2021.

\bibitem{Gagnon2022FredholmHeat}
L.~Gagnon, A.~Hayat, S.~Xiang, and C.~Zhang.
\newblock Fredholm transformation on {Laplacian} and rapid stabilization for
  the heat equation.
\newblock {\em Journal of Functional Analysis}, 283(12):109664, 2022.

\bibitem{gagnon2022fredholm}
L.~Gagnon, A.~Hayat, S.~Xiang, and C.~Zhang.
\newblock Fredholm backstepping for critical operators and application to rapid
  stabilization for the linearized water waves.
\newblock {\em Annales de l'Institut Fourier}, 2025.

\bibitem{hayat2024fredholm}
A.~Hayat and E.~Loko.
\newblock Rapid stabilization of general linear systems with {$F$}-equivalence.
\newblock {\em hal-04597173v2}, 2024.

\bibitem{ascencio2018backstepping}
P.~Ascencio, A.~Astolfi, and T.~Parisini.
\newblock Backstepping {PDE} design: A convex optimization approach.
\newblock {\em IEEE Transactions on Automatic Control}, 63(7):1943--1958, 2018.

\bibitem{KARAFYLLIS2017SampledTransport}
I.~Karafyllis and M.~Krstic.
\newblock Sampled-data boundary feedback control of {1-D} linear transport
  {PDEs} with non-local terms.
\newblock {\em Systems \& Control Letters}, 107:68--75, 2017.

\bibitem{davo2019SampledData}
M.~A. Davó, D.~Bresch-Pietri, C.~Prieur, and F.~Di~Meglio.
\newblock Stability analysis of a ${\text{2}}\times {\text{2}}$ linear
  hyperbolic system with a sampled-data controller via backstepping method and
  looped-functionals.
\newblock {\em IEEE Transactions on Automatic Control}, 64(4):1718--1725, 2019.

\bibitem{BekiarisLiberis2020}
N.~Bekiaris-Liberis.
\newblock Hybrid boundary stabilization of linear first-order hyperbolic {PDEs}
  despite almost quantized measurements and control input.
\newblock {\em Systems \& Control Letters}, 146:104809, 2020.

\bibitem{KARAFYLLIS2018SampledParabolic}
I.~Karafyllis and M.~Krstic.
\newblock Sampled-data boundary feedback control of {1-D} parabolic {PDEs}.
\newblock {\em Automatica}, 87:226--237, 2018.

\bibitem{Espitia2018EventBased}
N.~Espitia, A.~Girard, N.~Marchand, and C.~Prieur.
\newblock Event-based boundary control of a linear $2\times 2$ hyperbolic
  system via backstepping approach.
\newblock {\em IEEE Transactions on Automatic Control}, 63(8):2686--2693, 2018.

\bibitem{Espitia2020}
N.~Espitia.
\newblock Observer-based event-triggered boundary control of a linear
  2$\times$2 hyperbolic systems.
\newblock {\em Systems \& Control Letters}, 138:104668, 2020.

\bibitem{Espitia2022}
N.~Espitia, J.~Auriol, H.~Yu, and M.~Krstic.
\newblock Traffic flow control on cascaded roads by event-triggered output
  feedback.
\newblock {\em International Journal of Robust and Nonlinear Control},
  32(10):5919--5949, 2022.

\bibitem{ESPITIA2021EventTriggeredReactionDiffusion}
N.~Espitia, I.~Karafyllis, and M.~Krstic.
\newblock Event-triggered boundary control of constant-parameter
  reaction–diffusion {PDEs}: A small-gain approach.
\newblock {\em Automatica}, 128:109562, 2021.

\bibitem{Rathnayake2022ObserverEventTriggered}
B.~Rathnayake, M.~Diagne, N.~Espitia, and I.~Karafyllis.
\newblock Observer-based event-triggered boundary control of a class of
  reaction–diffusion {PDEs}.
\newblock {\em IEEE Transactions on Automatic Control}, 67(6):2905--2917, 2022.

\bibitem{Rathnayake2022SampledDataEventTriggered}
B.~Rathnayake, M.~Diagne, and I.~Karafyllis.
\newblock Sampled-data and event-triggered boundary control of a class of
  reaction–diffusion {PDEs} with collocated sensing and actuation.
\newblock {\em Automatica}, 137:110026, 2022.

\bibitem{Li2021}
X.~Li, Y.~Liu, J.~Li, and F.~Li.
\newblock Adaptive event-triggered control for a class of uncertain hyperbolic
  {PDE}-{ODE} cascade systems.
\newblock {\em International Journal of Robust and Nonlinear Control}, 2021.

\bibitem{Wang2023EventTriggeredParabolicPDEODE}
J.~Wang and M.~Krstic.
\newblock Event-triggered adaptive control of a parabolic {PDE–ODE} cascade
  with piecewise-constant inputs and identification.
\newblock {\em IEEE Transactions on Automatic Control}, 68(9):5493--5508, 2023.

\bibitem{Yu2019}
H.~Yu and M.~Krstic.
\newblock Traffic congestion control for {Aw-Rascle-Zhang} model.
\newblock {\em Automatica}, 100:38--51, 2019.

\bibitem{Yu2022}
H.~Yu, J.~Auriol, and M.~Krstic.
\newblock Simultaneous downstream and upstream output-feedback stabilization of
  cascaded freeway traffic.
\newblock {\em Automatica}, 136:110044, 2022.

\bibitem{guan2021}
L.~Guan, L.~Zhang, and C.~Prieur.
\newblock Optimal observer-based output feedback controller for traffic
  congestion with bottleneck.
\newblock {\em International Journal of Robust and Nonlinear Control},
  31(15):7087--7106, 2021.

\bibitem{guan2023}
L.~Guan, L.~Zhang, and C.~Prieur.
\newblock Controller design for heterogeneous traffic with bottleneck and
  disturbances.
\newblock {\em Automatica}, 148:110790, 2023.

\bibitem{guan2023disturb}
L.~Guan, C.~Prieur, L.~Zhang, and R.~Vazquez.
\newblock State observation for heterogeneous quasilinear traffic flow system
  with disturbances.
\newblock {\em Mathematics of Control, Signals, and Systems}, 36:601--628,
  2023.

\bibitem{Burkhardt2021}
M.~Burkhardt, H.~Yu, and M.~Krstic.
\newblock Stop-and-go suppression in two-class congested traffic.
\newblock {\em Automatica}, 125:109381, 2021.

\bibitem{Yu2021}
H.~Yu, Q.~Gan, A.~Bayen, and M.~Krstic.
\newblock {PDE} traffic observer validated on freeway data.
\newblock {\em IEEE Transactions on Control Systems Technology}, 29(3), 2021.

\bibitem{zhang2024mean}
Y.~Zhang, H.~Yu, J.~Auriol, and M.~Pereira.
\newblock Mean-square exponential stabilization of mixed-autonomy traffic {PDE}
  system.
\newblock {\em Automatica}, 170:111859, 2024.

\bibitem{kog2020_c}
S.~Koga, D.~Straub, M.~Diagne, and M.~Krstic.
\newblock Stabilization of filament production rate for screw extrusion-based
  polymer three-dimensional-printing.
\newblock {\em Journal of Dynamic Systems, Measurement, and Control},
  142(3):031005, 2020.

\bibitem{kog2020_a}
S.~Koga, M.~Krstic, and J.~Beaman.
\newblock Laser sintering control for metal additive manufacturing by {PDE}
  backstepping.
\newblock {\em IEEE Transactions on Control Systems Technology},
  28(5):1928--1939, 2020.

\bibitem{kog2020_d}
S.~Koga and M.~Krstic.
\newblock Arctic sea ice state estimation from thermodynamic {PDE} model.
\newblock {\em Automatica}, 112:108713, 2020.

\bibitem{cai2016nonlinear}
X.~Cai and M.~Krstic.
\newblock Nonlinear stabilization through wave {PDE} dynamics with a moving
  uncontrolled boundary.
\newblock {\em Automatica}, 68:27--38, 2016.

\bibitem{and2019}
G.~A. de~Andrade, R.~Vazquez, D.~J. Pagano, and J.~Mascheroni.
\newblock Design and implementation of a backstepping controller for regulating
  temperature in 3{D} printers based on selective laser sintering.
\newblock In {\em 2019 IEEE 58th Conference on Decision and Control (CDC)},
  pages 1183--1188. IEEE, 2019.

\bibitem{vaz2008_d}
R.~Vazquez, E.~Schuster, and M.~Krstic.
\newblock Magnetohydrodynamic state estimation with boundary sensors.
\newblock {\em Automatica}, 44(10):2517--2527, 2008.

\bibitem{diagne2017backstepping}
A.~Diagne, M.~Diagne, S.~Tang, and M.~Krstic.
\newblock Backstepping stabilization of the linearized
  {S}aint-{V}enant--{E}xner model.
\newblock {\em Automatica}, 76:345--354, 2017.

\bibitem{Aamo2016}
O.~M. Aamo.
\newblock Leak detection, size estimation and localization in pipe flows.
\newblock {\em {IEEE} Transactions on Automatic Control}, 61(1):246--251, 2016.

\bibitem{Anfinsen2022}
H.~Anfinsen and O.~M. Aamo.
\newblock Leak detection, size estimation and localization in branched pipe
  flows.
\newblock {\em Automatica}, 140:110213, 2022.

\bibitem{Aarsnes2019b}
U.~J.~F. Aarsnes and N.~van~de Wouw.
\newblock Axial and torsional self-excited vibrations of a distributed
  drill-string.
\newblock {\em Journal of Sound and Vibration}, 444:127--151, 2019.

\bibitem{Sagert2013}
C.~Sagert, F.~Di~Meglio, M.~Krstic, and P.~Rouchon.
\newblock Backstepping and flatness approaches for stabilization of the
  stick-slip phenomenon for drilling.
\newblock {\em IFAC Proceedings Volumes}, 46(2):779--784, 2013.

\bibitem{Aarsnes2019}
U.~J.~F. Aarsnes, J.~Auriol, F.~Di~Meglio, and R.~Shor.
\newblock Estimating friction factors while drilling.
\newblock {\em Journal of Petroleum Science and Engineering}, 179:80--91, 2019.

\bibitem{Auriol2021}
J.~Auriol, N.~Kazemi, and S.-I. Niculescu.
\newblock Sensing and computational frameworks for improving drill-string
  dynamics estimation.
\newblock {\em Mechanical Systems and Signal Processing}, 160:107836, 2021.

\bibitem{frihauf2010leader}
P.~Frihauf and M.~Krstic.
\newblock Leader-enabled deployment onto planar curves: A {PDE}-based approach.
\newblock {\em IEEE Transactions on Automatic Control}, 56(8):1791--1806, 2010.

\bibitem{meurer2011finite}
T.~Meurer and M.~Krstic.
\newblock Finite-time multi-agent deployment: {A} nonlinear {PDE} motion
  planning approach.
\newblock {\em Automatica}, 47(11):2534--2542, 2011.

\bibitem{zhang2024multi}
J.~Zhang, R.~Vazquez, J.~Qi, and M.~Krstic.
\newblock Multi-agent deployment in 3-{D} via reaction--diffusion system with
  radially-varying reaction.
\newblock {\em Automatica}, 161:111491, 2024.

\bibitem{qi2017wave}
J.~Qi, J.~Zhang, and Y.~Ding.
\newblock Wave equation-based time-varying formation control of multiagent
  systems.
\newblock {\em IEEE Transactions on Control Systems Technology},
  26(5):1578--1591, 2017.

\bibitem{qi2019control}
J.~Qi, S.~Wang, J.~Fang, and M.~Diagne.
\newblock Control of multi-agent systems with input delay via {PDE}-based
  method.
\newblock {\em Automatica}, 106:91--100, 2019.

\bibitem{freudenthaler2020experiments}
G.~Freudenthaler and T.~Meurer.
\newblock {PDE}-based multi-agent formation control using flatness and
  backstepping: {A}nalysis, design and robot experiments.
\newblock {\em Automatica}, 115:108897, 2020.

\bibitem{vazquez2023power}
R.~Vazquez, G.~Chen, J.~Qiao, and M.~Krstic.
\newblock The power series method to compute backstepping kernel gains: theory
  and practice.
\newblock In {\em 2023 62nd IEEE Conference on Decision and Control (CDC)},
  pages 8162--8169. IEEE, 2023.

\bibitem{lin2024towards}
X.~Lin, R.~Vazquez, and M.~Krstic.
\newblock Towards a {MATLAB} toolbox to compute backstepping kernels using the
  power series method.
\newblock {\em arXiv preprint arXiv:2403.16070}, 2024.

\bibitem{Vazquez2014a}
R.~Vazquez and M.~Krstic.
\newblock Marcum {Q}-functions and explicit kernels for stabilization of
  2$\times$ 2 linear hyperbolic systems with constant coefficients.
\newblock {\em Systems \& Control Letters}, 68:33--42, 2014.

\bibitem{marcum1950table}
J.~Marcum.
\newblock Table of {Q}-functions, {US} air force project {RAND} {R}es. {M}emo.
  {M}-339, {ASTIA} document {AD} 1165451.
\newblock {\em RAND Corp}, 1950.

\bibitem{wilhelmsen2022explicit}
N.~Wilhelmsen and O.~M. Aamo.
\newblock Explicit backstepping kernel solutions for leak detection in branched
  pipe flows.
\newblock {\em IEEE Control Systems Letters}, 7:913--918, 2022.

\bibitem{Woittennek2017}
F.~Woittennek, M.~Riesmeier, and S.~Ecklebe.
\newblock On approximation and implementation of transformation based feedback
  laws for distributed parameter systems.
\newblock {\em IFAC-PapersOnLine}, 50(1):6786--6792, 2017.
\newblock 20th IFAC World Congress.

\bibitem{ecklebe2017approximation}
S.~Ecklebe, M.~Riesmeier, and F.~Woittennek.
\newblock Approximation and implementation of transformation based feedback
  laws for distributed parameter systems.
\newblock {\em Proceedings in Applied Mathematics and Mechanics},
  17(1):785--786, 2017.

\bibitem{Auriol2019c}
J.~Auriol, K.~Morris, and F.~Di~Meglio.
\newblock Late-lumping backstepping control of {P}artial {D}ifferential
  {E}quations.
\newblock {\em Automatica}, 100:247--259, 2019.

\bibitem{yu2021reinforcement}
H.~Yu, S.~Park, A.~Bayen, S.~Moura, and M.~Krstic.
\newblock Reinforcement learning versus {PDE} backstepping and {PI} control for
  congested freeway traffic.
\newblock {\em IEEE Transactions on Control Systems Technology},
  30(4):1595--1611, 2021.

\bibitem{lee2025traffic}
J.~W. Lee, H.~Wang, K.~Jang, N.~Lichtl{\'e}, A.~Hayat, M.~Bunting, A.~Alanqary,
  W.~Barbour, Z.~Fu, X.~Gong, et~al.
\newblock Traffic control via connected and automated vehicles ({C}{A}{V}s): An
  open-road field experiment with 100 {CAV}s.
\newblock {\em IEEE Control Systems}, 45(1):28--60, 2025.

\bibitem{Bhan2023}
L.~Bhan, Y.~Shi, and M.~Krstic.
\newblock Neural operators for bypassing gain and control computations in {PDE}
  backstepping.
\newblock {\em IEEE Transactions on Automatic Control}, 69(8):5310--5325, 2024.

\bibitem{vazquez2024gain}
R.~Vazquez and M.~Krstic.
\newblock Gain-only neural operators for {P}{D}{E} backstepping.
\newblock {\em arXiv preprint arXiv:2403.19344}, 2024.

\bibitem{krstic2023neural}
M.~Krstic, L.~Bhan, and Y.~Shi.
\newblock Neural operators of backstepping controller and observer gain
  functions for reaction-diffusion {PDE}s.
\newblock {\em Automatica}, 164:111649, 2024.

\bibitem{redaud2024physics}
J.~Redaud, M.~Darrin, N.~Kazemi, and J.~Auriol.
\newblock Physics-informed dual architecture neural networks for enhanced
  estimation of drilling dynamics.
\newblock In {\em IGARSS 2024-2024 IEEE International Geoscience and Remote
  Sensing Symposium}, pages 1038--1042. IEEE, 2024.

\bibitem{wang2023backstepping}
S.~Wang, M.~Diagne, and M.~Krstic.
\newblock Backstepping neural operators for $2\times 2$ hyperbolic {PDEs}.
\newblock {\em arXiv preprint arXiv:2312.16762}, 2023.

\bibitem{lamarque2024adaptive}
M.~Lamarque, L.~Bhan, Y.~Shi, and M.~Krstic.
\newblock Adaptive neural-operator backstepping control of a benchmark
  hyperbolic {PDE}.
\newblock {\em arXiv preprint arXiv:2401.07862}, 2024.

\bibitem{lamarque2024gain}
M.~Lamarque, L.~Bhan, R.~Vazquez, and M.~Krstic.
\newblock Gain scheduling with a neural operator for a transport {PDE} with
  nonlinear recirculation.
\newblock {\em IEEE Transactions on Automatic Control}, 2025.

\bibitem{Auriol2020c}
J.~Auriol, F.~Bribiesca-Argomedo, and D.~Bresch-Pietri.
\newblock Stabilization of an underactuated 1+2 linear hyperbolic system with a
  proper control.
\newblock In {\em 2020 59th IEEE Conference on Decision and Control (CDC)},
  2020.

\bibitem{Tsubakino2012InDomain}
D.~Tsubakino, M.~Krstic, and S.~Hara.
\newblock Backstepping control for parabolic {PDEs} with in-domain actuation.
\newblock In {\em 2012 American Control Conference (ACC)}, pages 2226--2231,
  2012.

\bibitem{WANG2013indomain}
S.~Wang and F.~Woittennek.
\newblock Backstepping-method for parabolic systems with in-domain actuation.
\newblock {\em IFAC Proceedings Volumes}, 46(26):43--48, 2013.
\newblock 1st IFAC Workshop on Control of Systems Governed by Partial
  Differential Equations.

\bibitem{TSUBAKINO2015observerweighted}
D.~Tsubakino and S.~Hara.
\newblock Backstepping observer design for parabolic {PDEs} with measurement of
  weighted spatial averages.
\newblock {\em Automatica}, 53:179--187, 2015.

\bibitem{Zhu2022distributed}
J.~Zhu, X.~Liu, and Z.~Zhou.
\newblock Stabilization of a wave equation with a distributed control.
\newblock {\em Asian Journal of Control}, 24(5):2796--2805, 2022.

\bibitem{REDAUD2023indomain}
J.~Redaud and J.~Auriol.
\newblock Backstepping stabilization of a clamped string with actuation inside
  the domain.
\newblock {\em IFAC-PapersOnLine}, 56(2):9936--9941, 2023.
\newblock 22nd IFAC World Congress.

\bibitem{zhang2023adaptive}
Z.~Zhang, Q.~Fang, and X.~Chen.
\newblock Adaptive stabilization of uncertain linear system with stochastic
  delay by {PDE} full-state feedback.
\newblock {\em IEEE Transactions on Automatic Control}, 2023.

\bibitem{Auriol2023b}
J.~Auriol, M.~Pereira, and B.~Kulcsar.
\newblock Mean-square exponential stabilization of coupled hyperbolic systems
  with random parameters.
\newblock In {\em IFAC World Congress 2023}, 2023.

\bibitem{guan2024robustness}
D.~Guan, J.~Qi, and M.~Diagne.
\newblock Robustness of reaction--diffusion {PDE}s predictor-feedback to
  stochastic delay perturbations.
\newblock {\em Automatica}, 167:111784, 2024.

\bibitem{wu2022disturbance}
Z.-H. Wu, H.-C. Zhou, F.~Deng, and B.-Z. Guo.
\newblock Disturbance observer-based boundary control for an anti-stable
  stochastic heat equation with unknown disturbance.
\newblock {\em IEEE Transactions on Automatic Control}, 2022.

\bibitem{velho2024stabilization}
G.~Velho, J.~Auriol, R.~Bonalli, and I.~Boussaada.
\newblock Stabilization and optimal control of interconnected {SDE}-scalar
  {PDE} system.
\newblock {\em IEEE Control Systems Letters}, 2024.

\bibitem{auriol2020robust}
J.~Auriol, U.~J.~F. Aarsnes, F.~Di~Meglio, and R.~Shor.
\newblock Robust control design of underactuated 2$\times$2 {PDE}-{ODE}-{PDE}
  systems.
\newblock {\em IEEE Control Systems Letters}, 5(2):469--474, 2020.

\bibitem{redaud2024domain}
J.~Redaud, J.~Auriol, and Y.~Le~Gorrec.
\newblock In domain dissipation assignment of boundary controlled
  {P}ort-{H}amiltonian systems using backstepping.
\newblock {\em Systems \& Control Letters}, 185:105722, 2024.

\bibitem{le2005dirac}
Y.~Le~Gorrec, H.~Zwart, and B.~Maschke.
\newblock Dirac structures and boundary control systems associated with
  skew-symmetric differential operators.
\newblock {\em SIAM Journal on Control and Optimization}, 44(5):1864--1892,
  2005.

\bibitem{ramirez2017backstepping}
H.~Ramirez, H.~Zwart, Y.~Le~Gorrec, and A.~Macchelli.
\newblock On backstepping boundary control for a class of linear
  port-{H}amiltonian systems.
\newblock In {\em 2017 IEEE 56th Annual Conference on Decision and Control},
  pages 658--663. IEEE, 2017.

\bibitem{vazquez2025backstepping}
R.~Vazquez.
\newblock Backstepping control laws for higher-dimensional {P}{D}{E}s: Spatial
  invariance and domain extension methods.
\newblock {\em arXiv preprint arXiv:2503.00225}, 2025.

\end{thebibliography}

 %\newpage
 \par\noindent
\parbox[t]{\linewidth}{
\noindent\parpic{\includegraphics[height=1.2in,width=1in,clip,keepaspectratio]{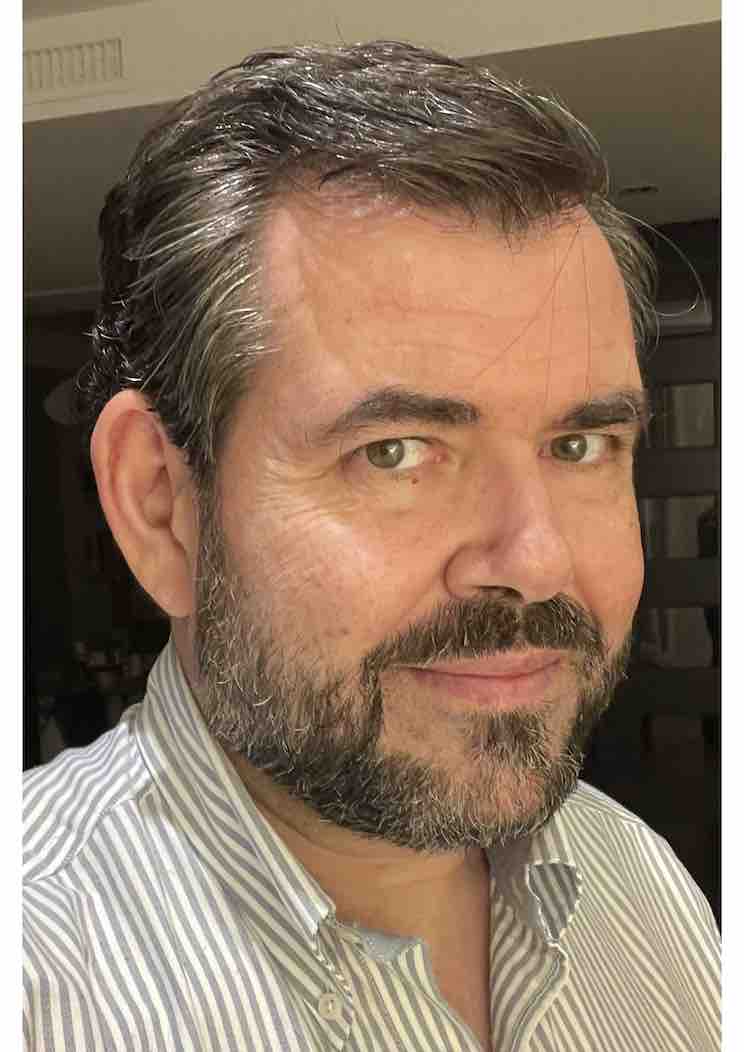}}
	\noindent {\bf Rafael Vazquez}  received the electrical engineering and mathematics degrees from the University of Seville, Seville, Spain, and the M.S. and Ph.D. degrees in aerospace engineering from the University of California, San Diego, La Jolla, CA, USA. He is currently Professor in the Aerospace Engineering Department and Director of the Space Surveillance  Chair at the University of Seville, Spain. His research interests include control theory, estimation, and optimization with applications to distributed parameter systems, spacecraft and aircraft guidance, navigation and control, and space surveillance and awareness. He is coauthor of more than 150 journal and conference publications and the book Control of Turbulent and Magnetohydrodynamic Channel Flows (Birkhauser). Prof. Vazquez currently serves as Associate Editor for Automatica and IEEE Control Systems Letters (L-CSS).}
 
\par\noindent
\parbox[t]{\linewidth}{
	\noindent\parpic{\includegraphics[height=1.2in,width=1in,clip,keepaspectratio]{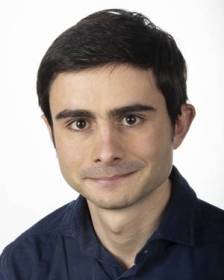}}
	\noindent {\bf Jean Auriol}   is a Researcher (Charg\'e de Recherches) at CNRS, Universit\'e Paris-Saclay, Centrale Supelec, Laboratoire des Signaux et Syst\`emes (L2S), Gif-sur-Yvette, France. He received his Ph.D. degree from Ecole des Mines Paris, PSL Research University in 2018, and obtained his Habilitation \`a Diriger des Recherches from Universit\'e Paris-Saclay in 2024. His research interests include robust control of hyperbolic systems, neutral systems, under-actuated networks and interconnected systems. Past and current applications of interest include Oil \& Gas and geothermal drilling, and the analysis of neural fields.}
%\vspace{4\baselineskip}
\par\noindent
\parbox[t]{\linewidth}{
	\noindent\parpic{\includegraphics[height=1.2in,width=1in,clip,keepaspectratio]{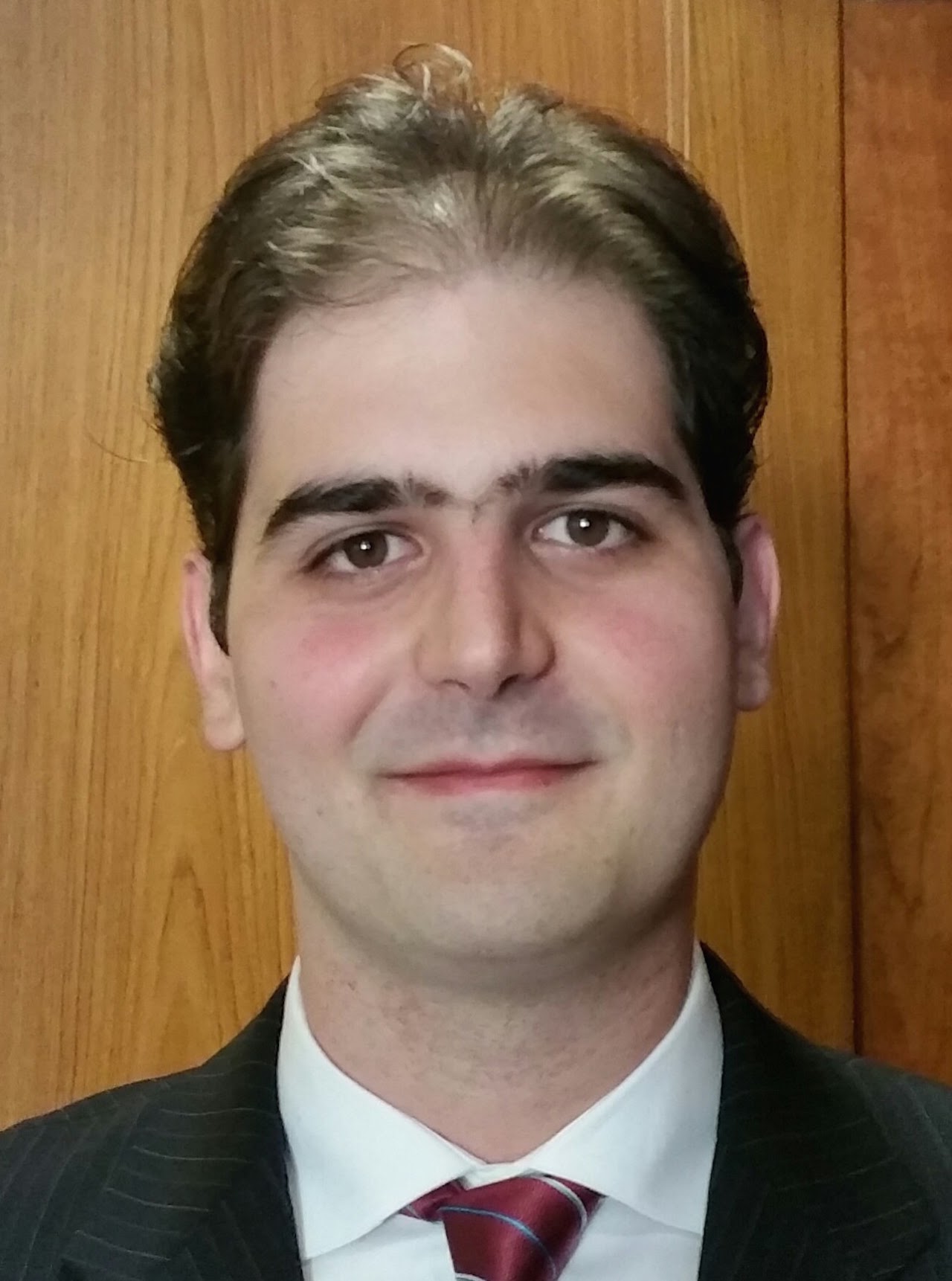}}
	\noindent {\bf Federico Bribiesca-Argomedo} received the Mechatronics Engineering degree from Tecnológico de Monterrey, Monterrey, Mexico., the M.Sc. in Control Systems from Grenoble INP, and the Ph.D. in Control Systems from the Université de Grenoble, in Grenoble, France. He is currently an Associate Professor (Maître de Conférences) at INSA Lyon, Laboratoire Ampère, Villeurbanne, France. His research interest include control and observer design, as well as parameter estimation and optimization for energy systems, including electrochemical energy storage systems and tokamak plasmas,
    with a particular focus on developing control-theoretical tools for systems comprised of parabolic and hyperbolic PDEs.}
%\vspace{4\baselineskip}
\par\noindent
\parbox[t]{\linewidth}{
	\noindent\parpic{\includegraphics[height=1.2in,width=1in,clip,keepaspectratio]{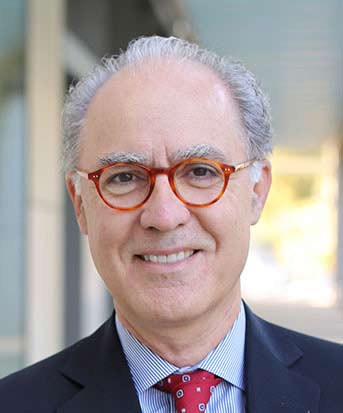}}
	\noindent {\bf Miroslav Krstic} is Distinguished Professor of Mechanical and Aerospace Engineering, holds the Alspach endowed chair, and is the founding director of the Center for Control Systems and Dynamics at UC San Diego. He also serves as Senior Associate Vice Chancellor for Research at UCSD. As a graduate student, Krstic won the UC Santa Barbara best dissertation award and student best paper awards at CDC and ACC. Krstic has been elected Fellow of seven scientific societies---IEEE, IFAC, ASME, SIAM, AAAS, IET (UK), and AIAA (Assoc. Fellow)---and as a foreign member of the Serbian Academy of Sciences and Arts and of the Academy of Engineering of Serbia. He has received the Richard E. Bellman Control Heritage Award, Bode Lecture Prize, SIAM Reid Prize, ASME Oldenburger Medal, Nyquist Lecture Prize, Paynter Outstanding Investigator Award, Ragazzini Education Award, IFAC Nonlinear Control Systems Award, IFAC Ruth Curtain Distributed Parameter Systems Award, IFAC Adaptive and Learning Systems Award, Chestnut textbook prize, AV Balakrishnan Award for the Mathematics of Systems, Control Systems Society Distinguished Member Award, the PECASE, NSF Career, and ONR Young Investigator awards, the Schuck (’96 and ’19) and Axelby paper prizes, and the first UCSD Research Award given to an engineer. Krstic has also been awarded the Springer Visiting Professorship at UC Berkeley, the Distinguished Visiting Fellowship of the Royal Academy of Engineering, the Invitation Fellowship of the Japan Society for the Promotion of Science, and four honorary professorships outside of the United States. He serves as Editor-in-Chief of Systems \& Control Letters and has been serving as Senior Editor in Automatica and IEEE Transactions on Automatic Control, as editor of two Springer book series, and has served as Vice President for Technical Activities of the IEEE Control Systems Society and as chair of the IEEE CSS Fellow Committee. Krstic has coauthored eighteen books on adaptive, nonlinear, and stochastic control, extremum seeking, control of PDE systems including turbulent flows, and control of delay systems.}
\end{document}